\documentclass[twoside,11pt,reqno]{amsart}
\usepackage{amsmath,amssymb,amscd,mathrsfs}

\makeatletter

\@addtoreset{equation}{section}

\def\thesubsection{\S\thesection-\alph{subsection}}
\def\Point#1{\addtocounter{subsection}{1}\vspace{2mm}
\noindent\thesubsection. {\bf #1.}\def\@currentlabel{\thesubsection}}

\newtheorem{Lemma}[equation]{Lemma}
\newtheorem{Theorem}[equation]{Theorem}
\newtheorem{Corollary}[equation]{Corollary}

\newtheorem{Remark}[equation]{Remark}

\newtheorem{Example}[equation]{Example}

\newtheorem{Conjecture}[equation]{Conjecture}

\newtheorem{Procedure}[equation]{Procedure}

\def\LL{{\mathtt{L}}}
\def\RR{{\mathtt{R}}}

\def\C{{\mathbb C}}
\def\Q{{\mathbb Q}}
\def\Z{{\mathbb Z}}
\def\N{{\mathbb N}}

\def\0{{\bar 0}}
\def\1{{\bar 1}}
\def\pr{{\operatorname{pr}}}

\def\H{{\mathscr H}}

\def\V{\mathscr{V}}
\def\W{\mathscr{W}}
\def\T{{\mathscr T}}
\def\TT{{_{\!\mathscr T}}}
\def\E{{\mathscr E}}
\def\EE{{_{\!\mathscr E}}}
\def\U{{\mathscr U}}

\def\A{{\mathscr A}}

\def\stab{{\operatorname{stab}}}

\def\sgn{{\operatorname{sgn}}}
\def\hom{{\operatorname{Hom}}}

\def\ext{{\operatorname{Ext}}}
\def\End{{\operatorname{End}}}

\def\soc{{\operatorname{soc}}}
\def\cosoc{{\operatorname{cosoc}}}
\def\ch{{\operatorname{ch}\:}}
\def\wt{{\operatorname{wt}}}

\def\underbar{\mathpalette\@underbar}
\def\@underbar#1#2{\settowidth{\@tempdimb}{$#1#2$}\@tempdimb=0.8\@tempdimb
                   \ooalign{$#1#2$\crcr%
                         \hfil\rule[-.5mm]{\@tempdimb}{.4pt}\hfil}}
\def\bs{\backslash}

\def\bid{{\text{\boldmath$1$}}}
\def\eps{{\varepsilon}}
\def\phi{{\varphi}}
\def\emptyset{{\varnothing}}

\def\la{{\lambda}}

\def\De{{\Delta}}

\newdimen\hoogte    \hoogte=12pt    
\newdimen\breedte   \breedte=14pt  
\newdimen\dikte     \dikte=0.5pt 

\newenvironment{Young}{\begingroup
       \def\vr{\vrule height0.89\hoogte width\dikte depth 0.2\hoogte}
       \def\fbox##1{\vbox{\offinterlineskip
                    \hrule height\dikte
                    \hbox to \breedte{\vr\hfill##1\hfill\vr}
                    \hrule height\dikte}}
       \vbox\bgroup \offinterlineskip \tabskip=-\dikte \lineskip=-\dikte
            \halign\bgroup &\fbox{##\unskip}\unskip  \crcr }
       {\egroup\egroup\endgroup}
\def\diagram#1{\relax\ifmmode\vcenter{\,\begin{Young}#1\end{Young}\,}\else%
              $\vcenter{\,\begin{Young}#1\end{Young}\,}$\fi}

\begin{document}
\title[Kazhdan-Lusztig polynomials]
{\boldmath
Kazhdan-Lusztig polynomials and character formulae for the Lie superalgebra
$\mathfrak{gl}(m|n)$}
\author{{\sc Jonathan Brundan}}
\thanks{Research partially supported by the NSF (grant no. DMS-0139019).}
\thanks{
{\em 2000 Subject Classification}: 17B10.}
\address
{Department of Mathematics\\ University of Oregon\\
Eugene\\ OR~97403, USA}
\email{brundan@darkwing.uoregon.edu}

\begin{abstract}
We compute the characters of the finite dimensional irreducible
representations of the Lie superalgebra $\mathfrak{gl}(m|n)$,
and determine $\ext$'s between simple modules
in the category of finite dimensional representations. We formulate
conjectures for the analogous results in category $\mathcal O$.
The combinatorics parallels the combinatorics
of certain canonical bases over the Lie algebra 
$\mathfrak{gl}(\infty)$.
\end{abstract}
\maketitle

\section{Introduction}

The problem of computing the characters of the 
finite dimensional irreducible representations of the
Lie superalgebra $\mathfrak{gl}(m|n)$ over $\C$ was raised originally
by V. Kac in 1977 \cite{Kac2,Kacnote}. Kac proved that
the finite dimensional universal highest weight modules, known 
nowadays as {\em Kac modules}, are irreducible for so-called 
{\em typical} highest weights, and gave
a formula for their characters.
After that, there were several conjectures and partial
results dealing with atypical highest weights 
\cite{BL, BR, sergeev0, HKJT, HKJT2, HKJ, KW, PS}, before
the complete solution to the problem was given by V. Serganova 
\cite{Serg0,Serg} in 1995 using a mixture of algebraic and geometric 
techniques.

In this article, we present a different, purely algebraic solution of the 
problem. One consequence is a 
proof of a conjecture made by van der Jeugt and
Zhang \cite{JZ}, which is apparently closely related to the conjectures 
made in \cite{HKJ}. In particular 
the composition multiplicities
of the Kac modules are all either $0$ or $1$, a fact which
does not seem to follow easily from Serganova's
formula since that involves certain alternating sums.
We also formulate for the first time a conjecture 
for the characters of the infinite dimensional
irreducible representations in the analogue of 
category $\mathcal O$ for the Lie superalgebra $\mathfrak{gl}(m|n)$.

Inspired by ideas of Lascoux, Leclerc and Thibon \cite{LLT},
our approach is to relate 
the finite dimensional representation theory of $\mathfrak{gl}(m|n)$
to the structure of the module
$$
\mathscr E^{m|n} := {\bigwedge}^m \V^* \otimes {\bigwedge}^n\V,
$$
where
$\V$ denotes the natural representation of the quantized
enveloping algebra $U_q(\mathfrak{gl}_\infty)$.
By work of Lusztig \cite[ch. 27]{Lubook}, the module $\E^{m|n}$
possesses a {canonical basis} $\{U_\la\}$ 
and a dual canonical basis $\{L_\la\}$, see
Theorems~\ref{thmB} and \ref{dthm}, which for the purpose
of this introduction we parametrize via (\ref{shift}) by
the set 
$X^+(m|n)$ of dominant integral weights
for $\mathfrak{gl}(m|n)$.
The entries of the transition matrices between these bases and the natural
monomial basis $\{K_\la\}$ of $\E^{m|n}$ define polynomials
$u_{\mu,\la}(q)$ and $l_{\mu,\la}(q)$ for each $\mu,\la \in X^+(m|n)$:
$$
U_\la = \sum_{\mu \in X^+(m|n)} u_{\mu,\la}(q) K_\mu,
\qquad
L_\la = \sum_{\mu \in X^+(m|n)} l_{\mu,\la}(q) K_\mu.
$$
The remarkable thing is that it is 
quite easy to compute these polynomials explicitly, 
because all the $\mathfrak{sl}_2$-strings 
in the crystal graph underlying the module $\E^{m|n}$ are of length $\leq 2$,
see Corollary~\ref{clearer} 
for the explicit formulae.
The main result of the article shows that
the polynomials $u_{\mu,\la}(q)$ evaluated 
at $q = 1$ compute the composition multiplicities of the Kac modules,
see Theorem~\ref{MT} and (\ref{lie}).
Moreover, the polynomials $l_{\mu,\la}(-q^{-1})$ 
coincide with the Kazhdan-Lusztig polynomials $K_{\la,\mu}(q)$ defined 
originally by Serganova \cite{Serg0, Serg}, so can be used to compute
$\ext$'s between Kac modules and irreducible modules, see
Theorem \ref{parcond} and Corollary \ref{extLL}. 

The module $\E^{m|n}$
is a summand of the tensor space
$$
\T^{m|n} := {\bigotimes}^m \V^* \otimes {\bigotimes}^n \V.
$$
The latter also
possesses a canonical basis $\{T_\la\}$
and a dual canonical basis $\{L_\la\}$, see Theorems~\ref{thmb} and
\ref{thmd}, which we parametrize via (\ref{shift})
by the set $X(m|n)$ of {\em all} integral weights for $\mathfrak{gl}(m|n)$.
For $\la \in X^+(m|n) \subset X(m|n)$, the elements $L_\la$ here
coincide with the images of the elements with the same name in the 
previous paragraph under the embedding $\E^{m|n} \hookrightarrow \T^{m|n}$.
The entries of the transition matrices between these bases
and the natural monomial basis $\{M_\la\}$ of $\T^{m|n}$
give us polynomials 
$t_{\mu,\la}(q)$ and $l_{\mu,\la}(q)$ for each $\mu,\la \in X(m|n)$:
$$
T_\la = \sum_{\mu \in X(m|n)} t_{\mu,\la}(q) M_\mu,
\qquad
L_\la = \sum_{\mu \in X(m|n)} l_{\mu,\la}(q) M_\mu.
$$
These 
should be viewed as the true combinatorial analogues for $\mathfrak{gl}(m|n)$
of the Kazhdan-Lusztig polynomials of
\cite{KL, Deo}.
We explain an explicit algorithm to compute $t_{\mu,\la}(q)$
in \ref{s5}, and conjecture based on calculations with
this algorithm that our polynomials 
share the positivity properties
of the usual Kazhdan-Lusztig polynomials, see Conjecture~\ref{pc}.
We conjecture moreover that the polynomials $t_{\mu,\la}(q)$ 
evaluated at $q = 1$
compute the composition multiplicities of the Verma modules
in category $\mathcal O$, see Conjecture~\ref{mc} and (\ref{mmult}).
This conjecture is true in the case $m = 0$
by the original Kazhdan-Lusztig conjecture \cite{KL}
for the Lie algebra $\mathfrak{gl}(n)$ proved in \cite{BB, BrK},
see Theorem~\ref{typthm}(i).

Returning to finite dimensional representations, 
let us now formulate the conjecture of van der Jeugt and Zhang proved
here precisely, to give the flavor of the combinatorics that arises. 
So let 
$\mathfrak g$  denote the Lie superalgebra $\mathfrak{gl}(m|n)$ over $\C$, 
labeling
rows and columns of matrices in $\mathfrak g$ by the ordered index set 
$I(m|n) = \{-m,\dots,-1,1,\dots,n\}$. 
We work always with the standard choices $\mathfrak h$
and $\mathfrak b$ of Cartan and Borel
subalgebras, consisting of diagonal and upper triangular matrices
respectively. 
For $i \in I(m|n)$, let $\delta_i \in \mathfrak h^*$ denote the function
picking out the $i$th diagonal entry of a diagonal matrix.
Put a symmetric bilinear form $(.|.)$ on $\mathfrak h^*$ by setting
$(\delta_i|\delta_j) = 1$ if $i=j>0$, $-1$ if $i=j < 0$, and 0 otherwise.
Let $W \cong S_m \times S_n$ denote the Weyl group
associated to $\mathfrak{g}_{\0}$, acting naturally
on $\mathfrak h^*$. We also 
need the dot action of $W$ on $\mathfrak{h}^*$: $w \cdot \la = w(\la+\rho) - \rho$
where $\rho = - \sum_{i \in I(m|n)} i \delta_i$.

Let $X(m|n) \subset \mathfrak h^*$ denote the set of all 
$\Z$-linear combinations of the weights $\{\delta_i\}_{i\in I(m|n)}$, and let
$X^+(m|n) \subset X(m|n)$ denote the {dominant integral 
weights}, namely,
the weights $\la = \sum_{i \in I(m|n)} \la_i \delta_i \in X(m|n)$
with $\la_{-m} \geq \dots \geq \la_{-1},
\la_1 \geq \dots \geq \la_n$.
Associated to $\la \in X^+(m|n)$, we have the 
Kac module $K(\la)$ and its unique irreducible quotient $L(\la)$,
see \ref{s9}.
We should note that there is no loss of generality in restricting our
attention to integral weights, since an arbitrary finite dimensional
irreducible 
representation of $\mathfrak{gl}(m|n)$ is either typical or can be obtained
from $L(\la)$ for some $\la \in X^+(m|n)$ by tensoring with a one dimensional
representation.

\vspace{2mm}\noindent
{\bf Main Theorem.} {\em
Let $\la \in X^+(m|n)$.
Let $r$ be maximal such that there exist
$-m \leq i_1 < \dots < i_r < 0 < j_r < \dots < j_1 \leq n$ with
$(\la+\rho|\delta_{i_s}-\delta_{j_s}) = 0$
for each $s = 1,\dots,r$.
Let $(k_1,\dots,k_r)$ be the lexicographically smallest tuple
of strictly positive integers such that
for all $\theta= (\theta_1,\dots,\theta_r) \in \{0,1\}^r$,
$\la + \sum_{s=1}^r \theta_s k_s (\delta_{i_s} - \delta_{j_s})$
is conjugate under the dot action of $W$ to a dominant weight,
denoted $\RR_{\theta}(\la) \in X^+(m|n)$. Then, for each $\mu \in X^+(m|n)$,
$$
[K(\mu):L(\la)] = \left\{
\begin{array}{ll}
1&\hbox{if $\mu = \RR_{\theta}(\la)$ 
for some $\theta = (\theta_1,\dots,\theta_r)
\in \{0,1\}^r$,}\\
0&\hbox{otherwise.}
\end{array}
\right.
$$
}
\vspace{1mm}

To prove the Main Theorem, we work with a different family of modules
$\{U(\la)\}_{\la \in X^+(m|n)}$
called {\em indecomposable tilting modules}, following the general
framework developed by Soergel \cite{so2} and extended to Lie
superalgebras in \cite{Btilt}. 
The problem of computing the multiplicities of Kac modules in 
indecomposable tilting modules is
roughly speaking transpose to the problem of 
computing the composition multiplicities of 
Kac modules, see (\ref{kdt}) for the precise relationship 
(a twisted BGG reciprocity).
The main step in the proof 
gives an explicit inductive construction
of the $U(\la)$'s starting from the typical case, when $U(\la) = K(\la)$,
and applying certain special translation functors that arise from tensoring
with the natural module and its dual.
Actually, we see eventually that 
the indecomposable tilting modules in this finite dimensional setting
{\em coincide} with the indecomposable projectives (also injectives),
but they are parametrized by highest weight rather than by
their irreducible quotients. 
Though one could just as well choose to work with the
latter more familiar labeling, the alternate
parameterization seems
to be the one that emerges naturally when considering
canonical bases.
There are also indecomposable tilting modules
denoted $\{T(\la)\}_{\la \in X(m|n)}$ in category $\mathcal O$, where
again they seem to correspond most directly to the canonical basis.

We now explain how the remainder of the article is organized.
In sections 2 and 3, we give the construction and properties 
of the canonical bases of the modules $\T^{m|n}$ and $\E^{m|n}$
from a purely combinatorial standpoint. 
Then in section 4 we describe the representation theory of 
$\mathfrak{gl}(m|n)$, working in two natural categories
$\mathcal O_{m|n}$ and $\mathcal F_{m|n}$ whose Grothendieck groups
are identified with the spaces
$\T^{m|n}$ and $\E^{m|n}$ respectively. 
In sections 2 and 3 we work exclusively
in a $\rho$-shifted notation which is more convenient for the combinatorics, 
replacing the set $X(m|n)$ of weights with the
set $\Z^{m|n}$ of functions $I(m|n)\rightarrow \Z$.
See (\ref{shift}) for the rule to translate between the two notations.

\vspace{2mm}
\noindent{\bf Acknowledgements. }
I would like to thank Arkady Berenstein and Jon Kujawa for helpful
conversations, Barrie Hughes for pointing out  the conjectures
made in \cite{HKJ, JZ}, and
Nathan Geer for asking questions that pushed me to
think about the results in \ref{klpol}.

\section{Tensor algebra}

In this section, we 
define 
and study the canonical basis of the tensor space $\T^{m|n}$.
We will work throughout 
over the field $\Q(q)$ of rational functions, where
$q$ is an indeterminate.

\Point{Combinatorial notation}\label{s1}
For $m,n \geq 0$,
let $S_{m|n}$ denote the symmetric group $S_m \times S_n$ acting
on the set $I(m|n) = \{-m,\dots,-1,1,\dots,n\}$ so that 
$S_m$ permutes $\{-m,\dots,-1\}$ and $S_n$ permutes $\{1,\dots,n\}$.
Thus $S_{m|n}$ is generated by the basic transpositions
$$s_{-m+1} = (-\!m\:\:-\!m\!+\!1),\dots, s_{-1} = (-\!2\:\:-\!1),s_1= (1\:\:2),\dots,s_{n-1} = 
(n\!-\!1\:\:n).
$$
Let $\Z^{m|n}$ be the set of all functions
$I(m|n) \rightarrow \Z$.
We call $f \in \Z^{m|n}$ {\em antidominant} if
$f(-m) \geq \dots \geq f(-1), f(1) \leq \dots \leq f(n)$.
Note $S_{m|n}$ acts on the right on $\Z^{m|n}$ by 
composition of functions, and every $f \in \Z^{m|n}$ is conjugate
under this action to a unique antidominant function.
We also have the `flip' $\omega:\Z^{m|n} \rightarrow \Z^{n|m}$,
where $\omega(f)$ is the function $I(n|m) \rightarrow \Z,
i \mapsto f(-i)$.

Let $P$ denote the free abelian group on basis
$\{\eps_a\:|\:a \in \Z\}$ endowed with a symmetric bilinear
form $(.,.)$ for which the $\eps_a$ form an orthonormal basis.
We view $P$ as the integral weight lattice associated
to the Lie algebra $\mathfrak{gl}_{\infty}$. The 
{\em simple roots} are the elements
$\eps_a - \eps_{a+1}\in P$ for $a \in\Z$.
The {\em dominance ordering} on $P$ is defined by
$\mu \leq \nu$ if $(\nu - \mu)$ is an
$\N$-linear combination of simple roots 
(here and later $\N = \{0,1,2,\dots\}$).
Equivalently, $\mu \leq \nu$ if
\begin{equation}\label{dim}
\sum_{b \leq a} (\mu, \eps_b) \leq \sum_{b \leq a} (\nu, \eps_b)
\end{equation}
for all $a \in \Z$ with equality for $a \gg 0$.

For $f \in \Z^{m|n}$ and $j \in I(m|n)$, define
\begin{equation}\label{wd}
\wt(f) := \sum_{i \in I(m|n)}\sgn(i) \eps_{f(i)},\qquad
\wt_{j} (f) := \sum_{j \leq i \in I(m|n)} \sgn(i) \eps_{f(i)},
\end{equation}
where $\sgn(i) \in \{\pm 1\}$ denotes the sign of $i$.
The {\em degree of atypicality} of $f \in \Z^{m|n}$ is defined to be
\begin{equation}\label{dat}
\#f := \frac{1}{2}\left(m+n - \sum_{a \in \Z} \big|(\wt(f), \eps_a)\big|\right).
\end{equation}
If $\#f = 0$, then $f$ is called {\em typical}.
So $f$ is typical if and only if
$$
\{f(-m),\dots,f(-1)\}\cap\{f(1),\dots,f(n)\} = \emptyset.
$$

\Point{Bruhat ordering}\label{bord}
Introduce a partial ordering on $\Z^{m|n}$ by declaring that
$g \preceq f$ if $\wt(g) = \wt(f)$
and
$\wt_{j}(g) \leq \wt_{j}(f)$ for all $j \in I(m|n)$.
It is immediate that if $g \preceq f$ then $\# g = \# f$.
Using (\ref{dim}), we see that
$g \preceq f$ if and only if
\begin{equation}\label{bruhat}
\sum_{
\stackrel
{\scriptstyle j \leq i \in I(m|n)} 
{g(i) \leq a}
} \sgn(i) \leq
\sum_{
\stackrel{\scriptstyle j \leq  i \in I(m|n)} 
{f(i) \leq a}
} \sgn(i)
\end{equation}
for all $a \in \Z$ and 
$j \in I(m|n)$, with equality if either $a \gg 0$ or $j = -m$.
From this, one gets in particular that $g \preceq f$ if and only if
$\omega(g) \preceq \omega(f)$.
In proofs, it will be convenient to have a shorthand for the sums appearing
in the inequality (\ref{bruhat}), 
so for $f \in \Z^{m|n}, a \in \Z$ and $j \in I(m|n)$
we abbreviate
\begin{equation*}
\#(f,a,j) = \sum_{\substack{j \leq i \in I(m|n) \\f(i)\leq a}}
\sgn(i).
\end{equation*}
Thus, $g \preceq f$ if and only if $\#(g,a,j) \leq \#(f,a,j)$ for all
$a \in \Z$ and $j \in I(m|n)$, with equality for $a \gg 0$ or
$j = -m$.

Here is another description of the partial order.
Let $d_i \in \Z^{m|n}$ be the function
$j \mapsto \sgn(i) \delta_{i,j},$ for each $i \in I(m|n)$. 
Write $f \downarrow g$ if one of the following holds:
\begin{itemize}
\item[(1)] $g = f - d_i + d_j$ for some $-m \leq i \leq -1, 1 \leq j \leq n$
such that $f(i) = f(j)$;
\item[(2)] $g = f \cdot (i\:\:j)$ for some $1 \leq i < j \leq n$
such that $f(i) > f(j)$;
\item[(3)] $g = f \cdot (i\:\:j)$ for some $-m \leq i < j \leq -1$
such that $f(i) < f(j)$.
\end{itemize}
Then:

\begin{Lemma}
$f \succeq g$ if and only if
there is a sequence $h_1, \dots, h_r \in \Z^{m|n}$ such that
$f = h_1 \downarrow \dots \downarrow h_r = g$.
\end{Lemma}

\begin{proof}
($\Leftarrow$) Obvious.

($\Rightarrow$) We show by induction on $(m+n)$ that if
$f \succ g$ are neigbors in the ordering, then $f \downarrow g$.
The case $m+n=0$ is vacuous, so suppose $m+n > 0$.
Replacing $f,g$ by $\omega(f),\omega(g)$ if necessary, we may assume 
in fact that $n > 0$.
If $f(n) = g(n)$, then we are done by induction, so
we may assume that $a = f(n) < g(n) = b$. We consider two cases.

\noindent
{\em Case one: }there exists $0 < i < n$ with $a < f(i) \leq b$.
Pick the greatest such $i$, so each $f(j)$ for $j=i+1,\dots,n$
is either $\leq a$ or $> b$, and set $c = f(i)$.
We claim that $f \succ f \cdot (i\:\:n) \succeq g$,
whence $f \downarrow g$ as required since $f$ and $g$ are neighbors.
For $i < j$ and $a \leq d < c$, we have that
$\#(f\cdot (i\:\:n),d,j) = \#(f,d,j) - 1$, while 
$\#(f \cdot (i\:\:n),d,j) = \#(f,d,j)$ for all other $j,d$.
Therefore to prove the claim, we just need to show that
$\#(f,d,j) > \#(g,d,j)$ for each $i < j$ and each $a \leq d < c$.
But by the choice of $i$, we have that
$\#(f,d,j) = \#(f,b,j) \geq \#(g,b,j) > \#(g,d,j)$ since $g(n) = b$.

\noindent
{\em Case two: }each $f(j)$ for $j = 1,\dots,n$ is either
$\leq a$ or $> b$.
From $(\wt(g), \eps_b) = (\wt(f),\eps_b) \leq 0$, we deduce that 
there must exist $-m \leq i < 0$ with $g(i) = b$. Take the greatest
such $i$. Now we claim that $f \succeq g + d_i - d_n \succ g$,
so again $f \downarrow g$ as they are neighbors.
To prove the claim, 
note that $\#(g+d_i-d_n,d,j) = \#(g,d,j)$ unless $j > i$ and $d = b-1$,
while $\#(g+d_i-d_n,b-1,j) = \#(g,b-1,j)+1$ for $j > i$.
Therefore we need to show that
$\#(f,b-1,j) > \#(g,b-1,j)$ for each $j > i$.
Now observe that
$\#(f,b-1,j) \geq \#(f,b,j) \geq \#(g,b,j) > \#(g,b-1,j)$.
\end{proof}

For example, writing elements of $\Z^{2|2}$ as tuples,
$$
(1,2|2,1) \downarrow (1,2|1,2) \downarrow (1,3|1,3) \downarrow 
(3,1|1,3).
$$
It is worth pointing out that $f \in \Z^{m|n}$ is minimal
with respect to the ordering just defined if and only if $f$ is typical
and antidominant.

\Point{\boldmath The quantum group}
Recall that the quantum integer associated to $n \geq 0$ is
$[n] := {(q^n-q^{-n})}/{(q-q^{-1})}$
and the quantum factorial is $[n]! := [n][n-1]\dots[2][1].$
Let $-:\Q(q) \rightarrow \Q(q)$ be the 
field automorphism induced by $q \mapsto q^{-1}$.
We will call an additive map
$f:V \rightarrow W$ between $\Q(q)$-vector spaces
{\em antilinear} if $f(c v) = \overline{c} f(v)$
for all $c \in \Q(q), v \in V$.

Let $\U$ denote the quantum group $U_q(\mathfrak{gl}_{\infty})$.
By definition, this is the $\Q(q)$-algebra on generators
$E_a, F_a, K_a, K_a^{-1}\:(a \in \Z)$
subject to relations
\begin{align*}
K_aK_a^{-1} &= K_a^{-1}K_a = 1,\\
K_aK_b &= K_bK_a,\\
K_aE_bK_a^{-1} &= q^{(\eps_a, \eps_b - \eps_{b+1})} E_b,\\
K_aF_bK_a^{-1} &= q^{(\eps_a, \eps_{b+1}-\eps_b)} 
F_b,
\end{align*}\begin{align*}
E_aF_b - F_bE_a &= \delta_{a, b} \frac{K_{a, a+1} - K_{a+1, a}}{q - q^{-1}},\\
E_aE_b &= E_bE_a&\hbox{if $|a-b|>1$},\\
E_a^2E_b  + E_bE_a^2 &= (q+q^{-1})E_aE_bE_a&\hbox{if $|a-b|=1$},\\
F_aF_b &= F_bF_a&\hbox{if $|a-b| > 1$,}\\
\qquad \qquad \qquad F_a^2F_b+ F_bF_a^2 &=   (q+q^{-1})F_aF_bF_a&\hbox{if $|a-b|=1$.}
\end{align*}
Here, for any $a, b \in\Z$, $K_{a, b}$ denotes
$K_aK_b^{-1}$. 
Also introduce the {\em divided powers} $F_a^{(r)} := F_a^r / [r]!$ and
$E_a^{(r)} := E_a^r / [r]!$.
We have the bar involution on $\U$, namely, the unique
antilinear automorphism such that $\overline{E_a} = E_a,
\overline{F_a} = F_a, \overline{K_a} = K_a^{-1}.$

We regard $\U$ as a Hopf algebra with comultiplication
$\Delta:\U \rightarrow \U \otimes \U$ defined on generators by
\begin{align*}
\Delta(E_a) &= 1 \otimes E_a + E_a \otimes K_{a+1, a},\\
\Delta(F_a) &= K_{a, a+1} \otimes F_a + F_a \otimes 1,\\
\Delta(K_a) &= K_a \otimes K_a.
\end{align*}
This is the comultiplication from Kashiwara \cite{Ka}, 
and is different from the one in Lusztig's book \cite{Lubook}.
The counit $\eps$ is defined by $\eps(E_a) = \eps(F_a) = 0$,
$\eps(K_a) = 1$, the antipode $S$
by
$S(E_a) = -E_a K_{a,a+1},
S(F_a) = -K_{a+1,a}F_a,
S(K_a) = K_a^{-1}.$

\Point{\boldmath The space $\T^{m|n}$}\label{s2}
Let $\V$ be the natural $\U$-module, 
with basis $\{v_a\}_{a\in \Z}$ and
action defined by
$$
K_a v_b = q^{\delta_{a,b}} v_b,
\qquad
E_a v_b = \delta_{a+1,b} v_a,
\qquad
F_a v_b = \delta_{a,b} v_{a+1}.
$$
Let $\W = \V^*$ be the dual $\U$-module, with
basis $\{w_a\}_{a\in\Z}$ related to the basis of $\V$ by
$\langle w_a, v_b \rangle = (-q)^{-a}\delta_{a,b}$.
The action of $\U$ on $\W$ satisfies
$$
K_a w_b = q^{-\delta_{a,b}} w_b,
\qquad
E_a w_b = \delta_{a,b}  w_{a+1},
\qquad
F_a w_b = \delta_{a+1,b}w_{a}.
$$
Let
$\T^{m|n} := 
\W^{\otimes m} \otimes \V^{\otimes n},$
viewed as a $\U$-module in the natural way.
Recall that $\Z^{m|n}$ denotes the set of all functions
$I(m|n)\rightarrow \Z$.
For $f \in \Z^{m|n}$, we let
$$
M_f =
w_{f(-m)} \otimes\dots\otimes w_{f(-1)} \otimes
v_{f(1)} \otimes \dots\otimes v_{f(n)}.
$$
The vectors 
$\{M_f\}_{f \in \Z^{m|n}}$ give a basis for $\T^{m|n}$.
A vector $v$ in a $\U$-module $M$
is said to be of {\em weight} $\nu \in P$ if
$K_a v = q^{(\nu,\eps_a)} v$
for all $a \in \Z$.
The weight of the vector $M_f$ is $\wt(f)$, as defined in (\ref{wd}).

We will often work with a 
completion $\widehat \T^{m|n}$ of $\T^{m|n}$.
To define this formally, let
$\Z^{m|n}_{\leq d}$ denote the set of all $f \in \Z^{m|n}$
with $f(i) \leq d$ for all $i \in I(m|n)$.
Let $\T^{m|n}_{\leq d}$ denote the subspace of $\T^{m|n}$
spanned by $\{M_f\}_{f \in \Z^{m|n}_{\leq d}}$, 
and let $\pi_{\leq d}:\T^{m|n} \rightarrow \T^{m|n}_{\leq d}$
denote projection along the basis.
The filtration $(\ker \pi_{\leq d})_{d \in \Z}$ induces
a topology on the abelian group
$\T^{m|n}$, see \cite[ch.III, $\S$2.5]{BouCA}.
Let 
$$
\widehat \T^{m|n} = \varprojlim \T^{m|n}_{\leq d}
$$ 
denote the corresponding completion, and identify $\T^{m|n}$
with its image in $\widehat{\T}^{m|n}$. 
The projections 
$\pi_{\leq d}$ extend by continuity to give maps 
${\pi}_{\leq d}:\widehat \T^{m|n} \rightarrow \T^{m|n}_{\leq d}$. 
As usual, 
we will view elements of $\widehat{\T}^{m|n}$ as infinite
$\Q(q)$-linear combinations of the basis elements $\{M_f\}_{f \in \Z^{m|n}}$ 
whose projections onto each $\T^{m|n}_{\leq d}$ are finite sums.
A homomorphism $\theta: \T^{m|n} \rightarrow \widehat{\T}^{m|n}$
of abelian groups satisfying the compatibility condition
$$
\hbox{
\em $\pi_{\leq d} (u) = 0$ implies ${\pi}_{\leq d}(\theta (u)) = 0$
for all $u \in \T^{m|n}$ and all $d\gg 0$}
$$
is automatically continuous, 
hence extends uniquely to a continuous endomorphism of $\widehat{\T}^{m|n}$.
In particular, the action of $\U$ lifts uniquely to a continuous action on 
$\widehat{\T}^{m|n}$, 
since $E_a, F_a$ and $K_a$ commute with $\pi_{\leq d}$
for all $d > a$.

\Point{The Iwahori-Hecke algebra}
Associated to the symmetric group $S_{m|n}$ we have the
Iwahori-Hecke algebra 
$\H_{m|n}$. This is defined as the
$\Q(q)$-algebra on generators
$H_{-m+1},\dots,H_{-1},H_1,\dots,H_{n-1}$ subject to relations
\begin{align*}
&H_i^2 = 1 - (q-q^{-1}) H_i,\\
&H_i H_{i+1} H_i = H_{i+1}H_i H_{i+1},\\
&H_i H_j = H_j H_i\quad\hbox{if }|i-j| > 1.
\end{align*}
For $x \in S_{m|n}$, we have the corresponding element 
$H_x \in \H_{m|n}$, 
where $H_x = H_{i_1} \dots H_{i_r}$ if $x = s_{i_1} \dots s_{i_r}$ 
is a reduced expression for $x$.
The bar involution on $\H_{m|n}$ is the 
unique antilinear automorphism such that $
\overline{H_x}
= H_{x^{-1}}^{-1},
$
in particular
$\overline{H_i} = H_i + (q-q^{-1})$.

We define a linear right action of $\H_{m|n}$ on $\T^{m|n}$ by the formulae
$$
M_f H_i = 
\left\{
\begin{array}{ll}
M_{f \cdot s_i}
& \hbox{if $f \prec f \cdot s_i$,}\\
q^{-1} M_f&\hbox{if $f = f \cdot s_i$,}\\
M_{f \cdot s_i} - (q-q^{-1})M_f&\hbox{if $f \succ f \cdot s_i$.}
\end{array}
\right.
$$
Since the action of $\H_{m|n}$ commutes with all $\pi_{\leq d}$,
it lifts by continuity to $\widehat{\T}^{m|n}$.
As is well-known, see e.g. \cite{Du}, the actions of $\U$ and $\H_{m|n}$ on $\T^{m|n}$
commute with one another, hence the actions on the completion
$\widehat{\T}^{m|n}$ also commute.

\Point{Some (anti)automorphisms}
Let $\sigma,\tau:\U \rightarrow \U$
be the antiautomorphisms
and $\omega:\U \rightarrow \U$
be the automorphism defined by
\begin{align*}
\sigma(E_a) = E_{-1-a}, \qquad
&\sigma(F_a) = F_{-1-a}, \qquad\quad\!
\sigma(K_a) = K_{-a},\\
\tau(E_a) = q^{-1} K_{a+1,a} F_a, \qquad
&\tau(F_a) = q E_a K_{a,a+1}, \quad
\tau(K_a) = K_{a},\\
\omega(E_a) = F_a, \qquad
&\omega(F_a) = E_a, \quad
\qquad\quad\:\omega(K_a) = K_{a}^{-1}.
\end{align*}
Let $\tau:\H_{m|n} \rightarrow \H_{m|n}$ 
be the antiautomorphism and
$\omega:\H_{m|n} \rightarrow \H_{n|m}$
be the isomorphism defined by 
$\tau(H_i) = H_i$ and
$\omega(H_i) = H_{-i}$
for $i \in I(m-1|n-1)$. 
Introduce the linear map
\begin{equation}\label{odef}
\omega:\T^{m|n} \rightarrow \T^{n|m}, \qquad M_f \mapsto 
M_{\omega(f)},
\end{equation}
where $\omega(f)$ is as in \ref{s1}.
Note $\omega$ extends by continuity to a linear map
$\widehat{\T}^{m|n} \rightarrow \widehat{\T}^{n|m}$.
Next let $(.,.)_{\TT}$ be the symmetric bilinear form on
$\T^{m|n}$ defined by
\begin{equation}\label{fdef}
(M_f, M_g)_\TT = \delta_{f,g}
\end{equation}
for $f,g \in \Z^{m|n}$.
Finally, define an antilinear map 
\begin{equation}\label{sdef}
\sigma:\T^{m|n} \rightarrow \T^{m|n},
\qquad
M_f \mapsto M_{-f}.
\end{equation}
The form $(.,.)_{\TT}$ and the map $\sigma$ do not extend to the completion.

\begin{Lemma}{\label{omega}}
\begin{itemize}
\item[(i)] $\omega(X u H) = \omega(X) \omega(u) \omega(H)$
for all $X \in \U, H \in \H_{m|n}$ and $u \in \widehat{\T}^{m|n}$.
\item[(ii)]
$(X u H, v)_\TT = (u, \tau(X) v \tau(H))_\TT$
for all $X \in \U, H \in \H_{m|n}$ and $u,v \in \T^{m|n}$.
\item[(iii)]
$\sigma(X u H) = \tau(\overline{\sigma(X)}) \sigma(u) \overline{H}$
for all $X \in \U, H \in \H_{m|n}$ and $u \in \T^{m|n}$.
\end{itemize}
\end{Lemma}

\begin{proof}
These are all checked directly for $\H_{m|n}$.
To prove them for $\U$, one first checks that
$\tau$ and $- \circ \sigma$ are coalgebra automorphisms and
$\omega$ is a coalgebra antiautomorphism of $\U$.
Hence it suffices to check (i)--(iii) when
$m+n=1$.
\end{proof}

\Point{Generation}
We proceed to prove that 
$\widehat{\T}^{m|n}$ is generated as a topological $\U$-module
by the vectors $M_f$ for {\em typical}
$f \in \Z^{m|n}$.

\begin{Lemma}\label{enasty}
Suppose that $f \in \Z^{m|n}$ and
$1 \leq i_1 < \dots < i_r \leq n$ are such that
$f(i_1) = \dots = f(i_r) = a+1$
and $f(j) \neq a, a+1$ for
all $j \in \{i_1,i_1+1,\dots,n\} - \{i_1,\dots,i_r\}$.
Let $f'$ be the function with
$f'(i_1) = \dots = f'(i_r) = a$ and $f'(j) = f(j)$
for all $j \neq i_1,\dots,i_r$.
Then, for any $g \preceq f$,
$$
E_a^{(r)} M_g \in \delta_{f, g} M_{f'} + 
\sum_{g' \prec f'} \Z[q,q^{-1}] M_{g'}.
$$
\end{Lemma}

\begin{proof}
Take $g \preceq f$.
Recall the definition of $d_j \in \Z^{m|n}$ from \ref{s1}.
Note $E_a^{(r)} M_g$ is a linear combination of $M_{g'}$'s where
$g' = g - d_{j_1} - \dots - d_{j_r}$
for $j_1 < \dots <j_r \in I(m|n)$ such that
$$
g(j_s) = \left\{\begin{array}{ll}a&\hbox{if $j_s < 0$},\\a+1&\hbox{if $j_s > 0$.}\end{array}\right.
$$
Let us show that for such a $g'$, we have that 
$g' \preceq f'$.
By (\ref{bruhat}), we need to show that
$\#(g',b,j) \leq \#(f',b,j)$
for all $b \in \Z$ and $j \in I(m|n)$. Since $g \preceq f$, we know that
$\#(g,b,j) \leq
\#(f,b,j)$.
So we are done except
possibly for $b = a$. Suppose then that 
$\#(g',a,j) > \#(f',a,j)$
for some $j$.
Say $i_1,\dots,i_s < j \leq  i_{s+1}, \dots, i_r$
and $j_1, \dots, j_t < j \leq j_{t+1}, \dots, j_r$.
Then, 
\begin{align*}
\#(f',a,j)
&<
\#(g',a,j)
=
\#(g,a,j) + (r-t)
\leq
\#(f,a,j)
 + (r-t)\\
&=
\#(f',a,j)
-(r-s)+(r-t)
=
\#(f',a,j)
+s-t.
\end{align*}
Hence, we must have that $s > t$.
This implies in particular that $j > 0$, and using this we get that
$$
\#(g,a+1,j)
\geq
\#(g',a,j)
>
\#(f',a,j)
=
\#(f,a+1,j), 
$$
which is a contradiction. So indeed we must have that $g' \preceq f'$.
Finally suppose that $g' = f'$.
The assumption that $f(j) \neq a,a+1$ for $j \in \{i_1, i_1+1,\dots,n\}-\{i_1,\dots,i_r\}$ means that we must have $j_1 \leq i_1,\dots,j_r \leq i_r$.
Hence, $f \preceq g$.
Since we started with the assumption that $g \preceq f$, we therefore
have $g = f$ which completes the proof.
\end{proof}

Twisting with $\omega$ using Lemma~\ref{omega}(i), 
we also have  the analogous statement for $F_a^{(r)}$:

\begin{Lemma}\label{fnasty}
Suppose that $f \in \Z^{m|n}$ and
$-m \leq i_r < \dots < i_1 \leq -1$ are such that
$f(i_1) = \dots = f(i_r) = a+1$
and $f(j) \neq a, a+1$ for
all $j \in \{-m,1-m,\dots,i_1\} - \{i_1,\dots,i_r\}$.
Let $f'$ be the function with
$f'(i_1) = \dots = f'(i_r) = a$ and $f'(j) = f(j)$
for all $j \neq i_1,\dots,i_r$.
Then, for any $g \preceq f$,
$$
F_a^{(r)} M_g \in \delta_{f, g} M_{f'} + 
\sum_{g' \prec f'} \Z[q,q^{-1}] M_{g'}.
$$
\end{Lemma}

\begin{Theorem}\label{triv}
We can write
each $M_f$ as a (possibly infinite) $\Z[q,q^{-1}]$-linear combination
of terms of the form $F_{a_1}^{(r_1)} \dots F_{a_s}^{(r_s)} M_g$
for $a_1,\dots,a_s \in \Z, r_1,\dots,r_s \geq 1$ and
{\em typical} $g \in \Z^{m|n}$.
\end{Theorem}

\begin{proof}
To prove the theorem, 
it suffices to show for each $d \in \Z$ and $f \in \Z^{m|n}$ that
we can write $M_f$ as a finite linear combination
of terms of the form
$F_{a_1}^{(r_1)} \dots F_{a_s}^{(r_s)} M_g$
for typical $g \in \Z^{m|n}$ modulo $\ker \pi_{\leq d}$.
So fix $d \in \Z$ and $f \in \Z^{m|n}$.
There
are only finitely many $g \preceq  f$ with $\pi_{\leq d} M_g \neq 0$,
So proceeding by induction on the dominance ordering, 
we may assume that
every $M_g$ with $g \prec f$ can be expressed
as a finite linear combination
of terms of the form
$F_{a_1}^{(r_1)} \dots F_{a_s}^{(r_s)} M_g$
for typical $g \in \Z^{m|n}$ modulo $\ker \pi_{\leq d}$.

Let $\{a_1 < a_2 < \dots < a_s\} = \{f(-m), \dots, f(-1)\}$
and let $r_t = \#\{i \in I(m|0)\:|\:f(i) = a_t\}$
for each $t = 1,\dots,s$.
Choose $k \gg 0$ so that every element of the set
$\{f(-m)+k, \dots, f(-1)+k\}$ exceeds every element of the
set $\{f(1), \dots, f(n)\}$.
Define $g \in \Z^{m|n}$ by
$$
g(i) = 
\left\{
\begin{array}{ll}
f(i)&\hbox{if $i > 0$,}\\
f(i)+k&\hbox{if $i < 0$.}
\end{array}\right.
$$
Note $g$ is typical by the choice of $k$. Now consider
$$
F_{a_s}^{(r_s)} \dots F_{a_s+k-1}^{(r_s)} 
F_{a_{s-1}}^{(r_{s-1})} \dots
F_{a_{s-1}+k-1}^{(r_{s-1})} 
\dots 
F_{a_1}^{(r_1)} \dots F_{a_1+k-1}^{(r_1)} M_g.
$$
One checks using Lemma~\ref{fnasty}
that this equals $M_f$ plus a $\Z[q,q^{-1}]$-linear combination
of $M_h$'s with $h \prec f$. So we are done by the induction hypothesis.
\end{proof}

\begin{Corollary}\label{swr}
Suppose $\theta:\widehat{\T}^{m|n} \rightarrow \widehat{\T}^{m|n}$
is a continuous $\U, \H_{m|n}$-bimodule endomorphism fixing 
$M_f$ for all typical antidominant $f \in \Z^{m|n}$.
Then $\theta$ is the identity map.
\end{Corollary}

\begin{proof}
If $f$ is antidominant, then
$M_f H_x = M_{f\cdot x}$ for all $x \in S_{m|n}$.
So for typical antidominant $f$ we have that
$\theta(M_{f \cdot x}) =  \theta(M_f H_x) = \theta(M_f) H_x = M_f H_x = M_{f \cdot x}$.
This shows that $\theta$ fixes $M_g$ for  {\em all} 
typical $g \in \Z^{m|n}$.
Now using the continuity of $\theta$ and Theorem~\ref{triv},
we get that $\theta$ fixes all $M_f$.
Hence by continuity again, $\theta$ is the identity map.
\end{proof}

\Point{\boldmath Canonical bases}\label{s4}
We now follow ideas of Lusztig \cite[ch. 27]{Lubook} to define a canonical
topological basis of $\widehat\T^{m|n}$.
We should note that in {\em loc. cit.}, 
Lusztig only considers finite dimensional
quantum groups, but the techniques generalize to our situation
on passing to the completion.
The first step in the construction is to introduce 
a bar involution on the space $\widehat\T^{m|n}$
that is compatible with the bar involutions on $\U$ and on $\H_{m|n}$.
The definition of this in Lusztig's work involves the quasi-$R$-matrix
associated to $\U$.
One gets from \cite[$\S$27.3]{Lubook} a bar involution
$-:\widehat\T^{m|n} \rightarrow \widehat\T^{m|n}$ that satisfies 
property (iv), hence (i), 
in the theorem below and that is compatible with the
bar involution on $\U$. One then checks easily using Lusztig's
definition that it is also compatible with the bar involution on $\H_{m|n}$,
giving the existence half of the proof of the theorem.
We will sketch a direct construction of the bar
involution on $\widehat\T^{m|n}$ below, independent of Lusztig's work.

\begin{Theorem}\label{thma}
There exists a unique continuous, antilinear map
 $-:\widehat{\T}^{m|n} \rightarrow
\widehat{\T}^{m|n}$
such that
\begin{itemize}
\item[(i)]
$\overline{M_f} = M_f$ for all typical antidominant
$f \in \Z^{m|n}$;
\item[(ii)]
$\overline{X u H} = \overline{X}\, \overline{u} \,\overline{H}$
for all $X \in \U, u \in \widehat{\T}^{m|n}$, $H \in \H_{m|n}$.
\end{itemize}
Moreover, 
\begin{itemize}
\item[(iii)] bar is an involution;
\item[(iv)] $\overline{M_f} = M_f + (*)$ where $(*)$ is
a (possibly infinite) $\Z[q,q^{-1}]$-linear combination
of $M_g$'s for $g \prec f$;
\item[(v)] $\overline{\omega (u)} = \omega(\overline{u})$ 
for all $u \in\widehat{\T}^{m|n}$.
\end{itemize}
\end{Theorem}

\begin{proof}
Let us first explain how to prove the uniqueness, and also the fact that
bar is necessarily an involution.
Suppose that we are given two continuous antilinear maps
$-,\sim:\widehat{\T}^{m|n} \rightarrow \widehat{\T}^{m|n}$
satisfying properties (i) and (ii).
Then we can consider the composite map
$$
\phi: \widehat{\T}^{m|n} \rightarrow \widehat{\T}^{m|n},
u \mapsto \widetilde{\overline{u}}.
$$
This is a continuous $\U, \H_{m|n}$-bimodule 
endomorphism of $\widehat{\T}^{m|n}$ fixing $M_f$ for all typical, antidominant
$f$.
Hence, $\phi$ is the identity map by Corollary~\ref{swr}.
In particular, this gives that
$\overline{\overline{u}} = u$ and $\widetilde{\widetilde{u}} = u$
for each $u \in \widehat{\T}^{m|n}$.
Finally, applying $\sim$ to both sides of the equation
$\widetilde{\overline{u}} = u$
gives that $\overline{u} = \widetilde{u}$ for all $u$, whence uniqueness.

To get existence without appealing directly to Lusztig's work, 
we need a little more notation.
For a partition $\la = (\la_1 \geq \la_2 \geq \dots)$, 
let $|\la|$ denote the sum of its parts,
$r(\lambda)$ denote the total number of non-zero parts, and
$r_s(\lambda)$ denote the total number of parts equal to $s$.
Also let
$$
p(\lambda) = (1-q^2)^{r(\lambda)} (-q)^{-|\lambda|}
\prod_{s \geq 1}
q^{r_s(\la)(r_s(\la)-1)/2}[r_s(\la)]!.
$$
Finally, for $a \in \Z$ and a partition $\lambda$,
let $a_\lambda^{m|n} \in \Z^{m|n}$ 
denote the function 
$i \mapsto a+\lambda_{|i|}$.
Let $D_\lambda^{m|n}$ 
denote 
the set of minimal length $\stab_{S_{m|n}}(a^{m|n}_\lambda) \bs S_{m|n}$-coset 
representatives.
For example, if $\la = (2,1,1)$, then
\begin{align*}
M_{0^{3|2}_\lambda} &= w_1\otimes 
w_1 \otimes w_2 \otimes v_2 \otimes v_1,\\
D_\lambda^{3|2} &= \{1, s_{-1}, s_{-1}s_{-2},
s_1, s_{-1}s_1, s_{-1}s_{-2}s_1\}.
\end{align*}
Now define
$-:\T^{m|n} \rightarrow \widehat{\T}^{m|n}$ to be the unique antilinear map 
satisfying the following properties:
\begin{itemize}
\item[(B1)] For $a \in \Z$,
$\displaystyle
\overline{w_a^{\otimes m} \otimes v_a^{\otimes n}}
=
\sum_{
\stackrel{\scriptstyle \lambda\text{ with}}{r(\lambda) \leq m,n}
}
p(\lambda)\bigg[\sum_{x \in D^{m|n}_\lambda}
q^{\ell(x)} M_{a^{m|n}_\lambda \cdot x}\bigg].$
\item[(B2)]
If $x \in S_{m|n}$ is the unique element of minimal length
such that
$f\cdot x$ is antidominant then 
$\overline{M_f} = \overline{M_{f\cdot x}} H_x^{-1}$.
\item[(B3)]
If all elements of
$\{f(-m),\dots,f(-k-1), f(l+1),\dots,f(n)\}$
are strictly greater than all elements of
$\{f(-k),\dots,f(-1),f(1),\dots,f(l)\}$ for some
$1 \leq k \leq m, 1 \leq l \leq n$, then
$$
\overline{M_f} =  
\sum_{
\stackrel{\scriptstyle g\in I(m-k|n-l)}{h\in I(k|l)}}
a_gb_h
w_{g(k-m)} \otimes\dots\otimes w_{g(-1)}
\otimes M_h  \otimes v_{g(1)} \otimes\dots\otimes
v_{g(n-l)}
$$
 where
\begin{align*}
\overline{w_{f(-m)} \otimes\dots\otimes w_{f(-k-1)}
\otimes v_{f(l+1)}\otimes\dots\otimes v_{f(n)}}
&=
\sum_{g\in I(m-k|n-l)} a_g M_g,\\
\overline{w_{f(-k)} \otimes\dots\otimes w_{f(-1)}
\otimes v_{f(1)}\otimes\dots\otimes v_{f(l)}}
&=
\sum_{h\in I(k|l)} b_h M_h.
\end{align*}
\end{itemize}
The following property is a consequence of (B1), and is useful
in inductive arguments. We omit the proof.
\vspace{1mm}
\begin{itemize}
\item[(B1)$'$]
For $m,n > 0$ and $a \in \Z$,
\begin{multline*}
\overline{w_a^{\otimes m} \otimes v_a^{\otimes n}}
= 
w_a \otimes \overline{w_a^{\otimes (m-1)} \otimes v_a^{\otimes n}}
+\\
q^{m-n}(q-q^{-1})
\sum_{l \geq 1} (-q)^{1-l}
w_{a+l} \otimes \overline{w_a^{\otimes(m-1)}\otimes v_a^{\otimes(n-1)}} 
\otimes v_{a+l}
X
\end{multline*}
where
$X = 1 + q H_{n-1}^{-1} + \dots + q^{n-1} H_{n-1}^{-1} \dots H_1^{-1}.$
\end{itemize}
For example:
\begin{align*}
\overline{w_a^{\otimes m}} &= w_a^{\otimes m},\qquad\qquad
\overline{v_a^{\otimes n}} = v_a^{\otimes n},\\
\overline{w_a \otimes v_a} &= 
w_a \otimes v_a + (q-q^{-1}) \sum_{l \geq 1} (-q)^{1-l}
w_{a+l} \otimes v_{a+l}.
\end{align*}

Now one checks easily that the map $-:\T^{m|n} \rightarrow \widehat\T^{m|n}$ 
just defined satisfies (iv), hence (i), and (v).
In particular, (iv) implies that bar is continuous,
so it extends uniquely to a continuous antilinear
map $-:\widehat{\T}^{m|n} \rightarrow \widehat{\T}^{m|n}$.
One finally needs to show that it satisfies (ii). 
This is done by a lengthy -- but elementary -- verification, using (B1), (B1)$'$, (B2) and (B3) directly.
We omit the details.
\end{proof}

Now we appeal to the following general lemma originating in 
\cite{KL}, also used implicitly in
\cite[$\S$27.3]{Lubook}. See \cite[1.2]{Du2} for a short proof.

\begin{Lemma}\label{dlem}
Let $(I,\preceq)$ be a partially ordered set with the property that 
$\{j \in I\:|\: j \preceq i\}$ is finite for all $i \in I$. 
Suppose that  $M$ is a $\Q(q)$-vector space with basis $\{u_i\}_{i \in I}$
equipped with an antilinear involution 
$-:M \rightarrow M$ such that
$\overline{u_i} = u_i + (*)$
for each $i \in I$,
where $(*)$ is a $\Z[q,q^{-1}]$-linear combination of 
$u_j$'s for $j \prec i$.
Then there exist unique bases
$\{x_i\}_{i \in I}$, $\{y_i\}_{i \in I}$
for $M$ such that
\begin{itemize}
\item[(i)] $\overline{x_i} = x_i$
and $\overline{y_i} = y_i$;
\item[(ii)]
$x_i \in u_i + \sum_{j \in I} q \Z[q] u_j$ and
$y_i \in u_i + \sum_{j \in I} q^{-1} \Z[q^{-1}] u_j,$
\end{itemize}
for each $i \in I$.
Moreover,
the coefficient of $u_j$ in $x_i$
(resp. $y_i$) is zero unless $j \preceq i$.
\end{Lemma}

Applying this to the space $\T^{m|n}_{\leq d}$ for fixed $d \in \Z$, the
basis $\{M_f\}_{f \in \Z^{m|n}_{\leq d}}$ and 
the antilinear involution 
$\pi_{\leq d} \circ - : \T^{m|n}_{\leq d} \rightarrow \T^{m|n}_{\leq d}$, 
we deduce:

\begin{Lemma}\label{lemmb}
There exist unique bases
$\{T_f^{(d)}\}_{f \in \Z_{\leq d}^{m|n}}$, 
$\{L_f^{(d)}\}_{f \in \Z_{\leq d}^{m|n}}$ for
${\T}_{\leq d}^{m|n}$ such that
\begin{itemize}
\item[(i)]
$\pi_{\leq d} (\overline{T_f^{(d)}}) = T_f^{(d)}$
and 
$\pi_{\leq d} (\overline{L_f^{(d)}}) = L_f^{(d)}$;
\item[(ii)]
$T_f^{(d)} \in M_f + {\sum}_{g \in \Z_{\leq d}^{m|n}} 
q \Z[q] M_g$
and
$L_f^{(d)} \in M_f + {\sum}_{g \in \Z_{\leq d}^{m|n}} q^{-1} 
\Z[q^{-1}] M_g$.
\end{itemize}
Moreover, the coefficient of $M_g$ in $T_f^{(d)}$ (resp. $L_f^{(d)}$)
is zero unless $g \preceq f$.
\end{Lemma}

Passing to the completion gives us bases
for $\widehat\T^{m|n}$:

\begin{Theorem}\label{thmb}
There exist unique topological bases
$\{T_f\}_{f \in \Z^{m|n}}$,
$\{L_f\}_{f \in \Z^{m|n}}$ for
$\widehat{\T}^{m|n}$ such that
\begin{itemize}
\item[(i)]
$\overline{T_f} = T_f$
and 
$\overline{L_f} = L_f$;
\item[(ii)]
$T_f \in M_f + \widehat{\sum}_{g \in \Z^{m|n}} q \Z[q] M_g$
and
$L_f \in M_f + \widehat{\sum}_{g \in \Z^{m|n}} q^{-1} 
\Z[q^{-1}] M_g$.
\end{itemize}
Moreover, 
\begin{itemize}
\item[(iii)]
$T_f = M_f + (*)$ and $L_f = M_f + (**)$ where
$(*)$ and $(**)$ are (possibly infinite) linear combinations of 
$M_g$'s for $g \prec f$;
\item[(iv)] $\omega(T_f) = T_{\omega(f)}$ and
$\omega(L_f) = L_{\omega(f)}$.
\end{itemize}
\end{Theorem}

\begin{proof}
Take $e \geq d$ and $f \in \Z^{m|n}_{\leq d}$.
Consider the elements 
$T_f^{(e)} \in \T^{m|n}_{\leq e}$ and
$T_f^{(d)} \in \T^{m|n}_{\leq d}$
given by Lemma~\ref{lemmb}.
We know that if $\pi_{\leq d} (u) = 0$ then
$\pi_{\leq d} (\overline{u}) = 0$.
Hence we have that
$\pi_{\leq d}(\overline{\pi_{\leq d}(u) - u}) = 0$
for all $u \in \widehat{\T}^{m|n}$.
Applying this to $u = T^{(e)}_f$ 
we deduce that
$
\pi_{\leq d} 
(\overline{
\pi_{\leq d}(T^{(e)}_f)
}
)
=
\pi_{\leq d}(T^{(e)}_f).
$
Hence by the uniqueness in Lemma~\ref{lemmb}, we have that
$\pi_{\leq d} (T^{(e)}_f) = T^{(d)}_f.$
Similarly,
$\pi_{\leq d} (L_f^{(e)}) = {L_f^{(d)}}.$
Hence, for all $f \in \Z^{m|n}$
there exist unique elements  $T_f, L_f \in \widehat{\T}^{m|n}$ such that
$$
\pi_{\leq d} (T_f) = T_f^{(d)},
\qquad
\pi_{\leq d} (L_f) = {L_f^{(d)}}
$$
for all $d \in \Z$ and all $f \in \Z^{m|n}_{\leq d}$.
Using the lemma for each $d \in \Z$,
one now easily checks that these satisfy (i)--(iii), while
(iv) follows from Theorem~\ref{thma}(v).
\end{proof}

We call the basis $\{T_f\}_{f\in \Z^{m|n}}$
the {\em canonical basis} of $\widehat\T^{m|n}$, and
$\{L_f\}_{f\in \Z^{m|n}}$ is the {\em dual canonical basis}.
Let us introduce notation for the coefficients: let
\begin{equation}\label{tde}
T_f = \sum_{g \in \Z^{m|n}} t_{g,f}(q) M_g,
\qquad
L_f = \sum_{g \in \Z^{m|n}} l_{g,f}(q) M_g
\end{equation}
for polynomials $t_{g,f}(q) \in \Z[q]$ and 
$l_{g,f}(q) \in \Z[q^{-1}]$.
We know that $t_{g,f}(q) = l_{g,f}(q) = 0$ unless
$g \preceq f$, and that
$t_{f,f}(q) = l_{f,f}(q) = 1$.

\begin{Example}\rm
If $m = n = 1$, the bases $\{T_f\}_{f \in \Z^{m|n}}$
and
$\{L_f\}_{f \in \Z^{m|n}}$ are
\begin{align*}
\{w_a \otimes v_b\:|\:a,b \in \Z, a \neq b\}
&\cup\{w_a \otimes v_a + q w_{a+1} \otimes v_{a+1}\:|\:a \in \Z\},\\
\{w_a \otimes v_b\:|\:a,b \in \Z, a \neq b\}
&\cup
\{w_a \otimes v_a - q^{-1} w_{a+1} \otimes v_{a+1}
\\&\qquad\qquad\:\:\:\:\,+ q^{-2} w_{a+2} \otimes v_{a+2} - \dots\:|\:a \in \Z\}
\end{align*}
respectively.
\end{Example}

\Point{Duality}\label{dsec}
We wish next to explain the relationship
between the bases $\{T_f\}_{f \in \Z^{m|n}}$ 
and $\{L_f\}_{f \in \Z^{m|n}}$.
Recall the definitions (\ref{fdef}) and (\ref{sdef}).
Define a new bilinear form $\langle.,.\rangle_\TT$ on
$\widehat{\T}^{m|n}$ by
\begin{equation}\label{od}
\langle u, v \rangle_\TT = 
(u, \sigma(\overline{v}))_\TT
\end{equation}
for $u, v \in \widehat{\T}^{m|n}$. 
Note this makes sense, even though the expression 
$\sigma(\overline{v})$ may not make sense in its own right.
Indeed, it is clear that on expanding $u$ and
$\overline{v}$ in terms of the basis
$\{M_f\}$, there are only finitely many $f$
such that $M_f$ is involved in $u$ and $M_{-f}$ is involved in
$\overline{v}$.
So we can interpret $(u, \sigma(\overline{v}))_\TT$ as
\begin{equation*}\label{interp}
\sum_{f \in \Z^{m|n}} (u, M_f)_\TT \overline{(M_{-f},\overline v)_\TT},
\end{equation*}
all but finitely many terms in the sum being zero.

\begin{Lemma}\label{sf} $\langle X u H, v \rangle_\TT
= \langle u, \sigma(X) v \tau(H) \rangle_\TT$
for all $X \in \U, H \in \H_{m|n}$ and $u,v \in \widehat{\T}^{m|n}$.
\end{Lemma}

\begin{proof}
According to Lemma~\ref{omega}, we have that
\begin{align*}
\langle X u H, v \rangle_\TT &= 
(X u H, \sigma(\overline{v}))_\TT
=
(u, \tau(X) \sigma(\overline{v}) \tau(H))_\TT\\
&=
(u, \tau(\overline{\sigma\overline{(\sigma(X))}}) 
\sigma(\overline{v}) \tau(H))_\TT
=
(u, 
\sigma(\overline{\sigma(X)}\overline{v} \overline{\tau(H)}))_\TT\\
&=
(u, \sigma(\overline{\sigma(X) v \tau(H)}))_\TT
=
\langle u, \sigma(X) v \tau(H) \rangle_\TT.
\end{align*}
\end{proof}

\begin{Lemma}\label{sym} 
The bilinear form $\langle.,. \rangle_\TT$ is symmetric.
\end{Lemma}

\begin{proof}
Let us first show that $\langle u, M_f \rangle_\TT 
= \langle M_f, u \rangle_\TT$ for all $u \in \widehat{\T}^{m|n}$
and typical antidominant $f \in \Z^{m|n}$.
We need to show that
$(u, M_{-f})_\TT = \overline{(\overline{u}, M_{-f})_\TT},$
for which it suffices to consider the special case $u = M_g$.
Then, $(M_g, M_{-f})_\TT = \delta_{g, -f}$.
Consider $(\overline{M_g}, M_{-f})_\TT$. By Theorem~\ref{thma}(iv), 
it is zero unless $\wt(g)  = \wt(-f)$.
So since $-f$ is typical, $g$ must be too.
By Theorem~\ref{thma}(iv),
$$
\overline{M_g} = M_g + (\hbox{a linear combination of
$M_h$'s with $h \prec g$}).
$$
Since $f$ is antidominant, 
we deduce 
that $(\overline{M_g}, M_{-f})_\TT = (M_g, M_{-f})_\TT = \delta_{g, -f}$.

Now we show that 
$\langle u, M_f \rangle_\TT 
= \langle M_f, u \rangle_\TT$ for all $u \in \widehat{\T}^{m|n}$
and all typical $f \in \Z^{m|n}$.
Let $x \in S_{m|n}$ be of minimal length such that $f \cdot x$ is 
antidominant. Then, by the previous paragraph and Lemma~\ref{sf}, we have that
\begin{align*}
\langle u, M_f \rangle_\TT
&= \langle u, M_{f \cdot x} T_{x^{-1}} \rangle_\TT
= \langle u T_x, M_{f \cdot x} \rangle_\TT
\\&= 
\langle M_{f \cdot x}, u T_x \rangle_\TT
= \langle M_{f\cdot x} T_{x^{-1}}, u \rangle_\TT
=
\langle M_f, u \rangle_\TT,
\end{align*}
as required.

Now let us consider the general case. 
In view of Theorem~\ref{triv}, we may assume that
$v = X M_g$ for typical $g$ and
$X = F_{a_1}^{(r_1)} \dots F_{a_s}^{(r_s)} \in \U$.
Then, by the previous paragraph and Lemma~\ref{sf}, we have that
\begin{align*}
\langle u, v \rangle_\TT &= 
\langle u, X M_g \rangle_\TT
= 
\langle \sigma(X) u, M_g \rangle_\TT\\ &=
\langle M_g, \sigma(X) u \rangle_\TT = \langle X M_g, u \rangle_\TT
= \langle v, u \rangle_\TT.
\end{align*}
This completes the proof.
\end{proof}

The following 
theorem characterizes $\{L_f\}_{f \in \Z^{m|n}}$
uniquely as the dual basis to $\{T_{-g}\}_{g \in \Z^{m|n}}$
under the bilinear form $\langle.,.\rangle_\TT$.

\begin{Theorem}\label{thmd}
For $f, g \in \Z^{m|n}$, 
$\langle L_f, T_{-g} \rangle_\TT = \delta_{f,g}$.
\end{Theorem}

\begin{proof}
Consider
$\langle L_f, T_{-g} \rangle_\TT = (L_f, \sigma(T_{-g}))_\TT.$
We observe that $L_f$ is equal to $M_f$ plus a $q^{-1}\Z[q^{-1}]$-linear combination
of $M_h$'s with $h \prec f$.
Also $\sigma(T_{-g})$ equals $M_g$ plus a $q^{-1} \Z[q^{-1}]$-linear
combination of $M_h$'s with $h \succ g$.
Hence, $\langle L_f, T_{-g} \rangle_\TT$ is zero unless
$f \succeq g$, it is $1$ if $f = g$ and it is in
$q^{-1} \Z[q^{-1}]$ if $f \succ g$.

On the other hand, by Lemma~\ref{sym}, 
$\langle L_f, T_{-g} \rangle_\TT = 
(\sigma(L_f), T_{-g})_\TT$.
Hence, arguing as in the previous paragraph, 
$\langle L_f, T_{-g} \rangle_\TT$ is zero unless
$f \succeq g$, it is $1$ if $f = g$ and it is in
$q \Z[q]$ if $f \succ g$.
Since $q \Z[q] \cap q^{-1} \Z[q^{-1}] = \{0\}$, this completes the proof.
\end{proof}

\begin{Corollary}\label{lcornew} For $f \in \Z^{m|n}$,
$$
M_f = \sum_{g \in \Z^{m|n}} t_{-f,-g}(q^{-1}) L_g
= \sum_{g \in \Z^{m|n}} l_{-f,-g}(q^{-1}) T_g.
$$
\end{Corollary}

\begin{proof}
By the theorem, we can write
$M_f = \sum_{g \in \Z^{m|n}} \langle M_f, T_{-g} \rangle_\TT L_g$.
Now a computation from the definition 
(\ref{od})
of the form $\langle.,.\rangle_\TT$
gives that  $\langle M_f, T_{-g} \rangle_\TT = t_{-f,-g}(q^{-1})$.
The second equality is proved similarly.
\end{proof}

\Point{\boldmath An algorithm}\label{s5}
The goal in this subsection is to explain an
algorithm to compute 
$T_f^{(d)}$ (cf. Lemma~\ref{lemmb})
for each $d \in \Z$ and $f \in\Z^{m|n}_{\leq d}$. 
Assuming a certain positivity conjecture which 
ensures that the $T_f^{(d)}$ converge
to $T_f$ in finitely many steps, the algorithm can be modified to
actually compute the canonical basis elements $T_f$ themselves.

The algorithm proceeds by induction on the degree of atypicality 
$\# f$ of $f$.
To begin with, we describe the base of the induction by 
explaining how to compute $T_f$ (hence all $T_f^{(d)}$) 
for {\em typical} $f$.
If $f$ is typical and 
antidominant, then we have that $T_f = M_f$ and we are done.
Otherwise, we can find $i \in I(m-1|n-1)$ such that
$f \cdot s_i \prec f$.
We may assume by induction that $T_{f \cdot s_i}$ is already known,
and consider the bar invariant element
$T_{f \cdot s_i} (H_i + q),$
which we view as a first approximation to $T_{f}$.
It equals $M_f$ plus a sum of terms $p_{g,f}(q) M_g$
for polynomials $p_{g,f}(q) \in \Z[q]$ and
$g$ with $g \prec f$ (there being only finitely many such $g$'s since $f$ is
typical).
For each such $g$ with $p_{g,f}(0) \neq 0$,
we make a correction by subtracting $p_{g,f}(0) T_g$
from our first approximation.
The result is a bar invariant 
expression that equals $M_f$ plus a $q\Z[q]$-linear combination of
$M_g$'s.
This must be $T_f$ by the uniqueness in Theorem~\ref{thmb}.

We have just described the usual algorithm to compute
the parabolic Kazhdan-Lusztig polynomials associated to 
the Hecke algebra $\H_{m|n}$ \cite{KL,Deo}.
To make this precise, let us recall the definition of the latter, following
\cite{soergel}.
Let $f \in \Z^{m|n}$ be antidominant.
Then, $S_f := \stab_{S_{m|n}} (f)$ is a parabolic subgroup of 
$S_{m|n}$. Let $\H_f$ be the corresponding parabolic subalgebra of
$\H_{m|n}$.
Let $\bid_{\H_f}$ denote the one dimensional right $\H_f$-module
on basis $1$ with action $1 H_i = q^{-1} 1$ for each 
$H_i \in \H_f$.
We consider the induced module
$\bid_{\H_f} \otimes_{\H_f} \H_{m|n}.$
This has a basis given by the elements $M_x^{(f)} := 1 \otimes H_x$
as $x$ runs over the set $D_f$ of minimal length
$S_f \bs S_{m|n}$-coset representatives.
The bar involution on $\bid_{\H_f} \otimes_{\H_f} \H_{m|n}$
is the antilinear map defined by
$\overline{1 \otimes H_x} := 1 \otimes \overline{H_x}.$
The Kazhdan-Lusztig basis 
of  $\bid_{\H_f} \otimes_{\H_f} \H_{m|n}$
is the unique
bar invariant basis $\{C^{(f)}_x\}_{x \in D_f}$ 
such that
$$
C_x^{(f)} \in M_x^{(f)} + \sum_{y \in D_f} q \Z[q] M_y^{(f)}.
$$
The corresponding parabolic Kazhdan-Lusztig polynomials are defined from
$$
C_x^{(f)} = \sum_{y \in D_f} m_{y,x}^{(f)}(q) M_y^{(f)}.
$$
The relationship to our situation is as follows:

\begin{Lemma}\label{tc}
Suppose that $f$ is typical and antidominant.
Then, for $x \in D_f$ and $g \in \Z^{m|n}$,
$$
\displaystyle
t_{g, f \cdot x}(q)
= 
\left\{
\begin{array}{ll}
m_{y,x}^{(f)}(q)&\hbox{if $g = f \cdot y$ for some $y \in D_f$},\\
0&\hbox{otherwise.}
\end{array}
\right.
$$
\end{Lemma}

\begin{proof}
Note that $M_f H_i = q^{-1} M_f$ for all $H_i \in \H_f$.
So we get from Frobenius reciprocity a unique $\H_{m|n}$-module homomorphism
$\theta: \bid_{\H_f} \otimes_{\H_f} \H_{m|n}
\rightarrow \T^{m|n}$
under which $M_x^{(f)}$ maps to $M_{f \cdot x}$.
By (B2) from the proof of Theorem~\ref{thma}, $\overline{M_{f \cdot x}} = M_f \overline{H_{x}}$.
So $\theta(\overline{u}) = \overline{\theta(u)}$
for each $u \in \bid_{\H_f} \otimes_{\H_f} \H_{m|n}$.
Therefore, $\theta(C_x^{(f)})$ is bar invariant, and it equals
$M_{f \cdot x}$ plus a $q \Z[q]$-linear combination of other
$M_g$'s.
Hence, $\theta(C_x^{(f)}) = T_{f \cdot x}$ by the uniqueness.
This shows that for each $x \in D_f$,
$$
T_{f \cdot x} = \sum_{y \in D_f} m_{y,x}^{(f)}(q) M_{f \cdot y}. 
$$
The lemma follows.
\end{proof}

\begin{Remark}\rm
In the case $m = 0$, this lemma shows that in
type $A$ the parabolic Kazhdan-Lusztig
polynomials 
coincide with the coefficients of the canonical basis of 
$\V^{\otimes n}$. This is a well-known observation, 
see for example \cite{FKK}.
\end{Remark}

Now we describe the algorithm to compute $T_f^{(d)}$ for atypical $f$.
We assume therefore that we are given $d \in \Z$ and
$f \in \Z^{m|n}_{\leq d}$ with $\# f > 0$, and that we have already
constructed an algorithm to compute $T_g^{(e)}$ for each $e \in \Z$
and $g$ with $\# g < \# f$.
Let us write $f_- = (f(-m), \dots, f(-1))$ and
$f_+ = (f(1),\dots,f(n))$.
Define $a_1$
to be the greatest integer that appears in both the tuples $f_-$ and $f_+$.
Now we iterate a certain {\em bumping procedure}:

Let $n_{1,1}$ 
be the number of entries equal to $a_1$ appearing in the tuple $f_+$,
and label all such entries.
If there are no entries equal to $(a_1+1)$ appearing to the right of
labeled $a_1$'s, move on to the next paragraph.
Otherwise, let $n_{1,2}$ be the number of entries equal to
$(a_1+1)$ appearing to the right of labeled $a_1$'s, and label
all such $(a_1+1)$'s. Next, if there are no $(a_1+2)$'s 
to the right of labeled $(a_1+1)$'s, move on to the next paragraph.
Otherwise let $n_{1,3}$ be the number of $(a_1+2)$'s to the right of labeled
$(a_1+1)$'s, and label all such $(a_1+2)$'s.
Continue in this way.

When the process just described terminates, we are left with
a sequence $n_{1,1}, n_{1,2}, \dots, n_{1,k_1}$ for
some $k_1 \geq 1$, where there are 
$n_{1,i}$ labeled $(a_1+i-1)$'s in the tuple $f_+$.
We define $X_1 := E_{a_1+k_1-1}^{(n_{1,k_1})} 
\dots E_{a_1+1}^{(n_{1,2})} E_{a_1}^{(n_{1,1})}$
and $a_2 := a_1+k_1$.
If there are no entries equal to $a_2$ in the tuple $f_-$, the
bumping procedure is finished. 
Otherwise, we need to repeat
the bumping procedure but applied to $f_-$ instead, as follows.
Let $n_{2,1}$ 
be the number of entries equal to $a_2$ appearing in the tuple $f_-$,
and label all such entries.
If there are no entries equal to $(a_2+1)$ appearing to the left of
labeled $a_2$'s, move on to the next paragraph.
Otherwise, let $n_{2,2}$ be the number of entries equal to
$(a_2+1)$ appearing to the left of labeled $a_2$'s, and label
all such $(a_2+1)$'s. Continue in this way until the process terminates.

We are left with
a sequence $n_{2,1}, n_{2,2}, \dots, n_{2,k_1}$ for
some $k_2 \geq 1$, where there are 
$n_{2,i}$ labeled $(a_2+i-1)$'s in the tuple $f_-$.
Let $X_2 := F_{a_2+k_2-1}^{(n_{2,k_2})} 
\dots F_{a_2+1}^{(n_{2,2})} F_{a_2}^{(n_{2,1})}$
and $a_3 := a_2+k_2$.
This time if there are no entries equal to $a_3$ in the tuple $f_+$, the
bumping procedure is finished. 
Otherwise, we repeat the whole process once more from the
beginning, but using $a_3$ in place of $a_1$, to construct $X_3, a_4, X_4, \dots$ and so on.

When the bumping procedure finally ends, we are left with a sequence
of monomials $X_1, \dots, X_N$ and integers
$a_1<a_2<\dots<a_{N+1}$.
Increase all labeled entries in the tuples $f_-, f_+$ by $1$ and let
$h \in \Z^{m|n}$ be the corresponding function.
Note that $\#h < \# f$. So by induction, we can compute $T_h^{(e)}$,
where $e = \max(d, a_{N+1})$.
Now consider the bar invariant element
$\pi_{\leq d} (X_N \dots X_1 T_h^{(e)}) \in \T^{m|n}_{\leq d},$
which is our first approximation to $T_f^{(d)}$.
By Lemmas~\ref{enasty} and \ref{fnasty}, 
it equals $M_f$ plus a finite linear combination of terms
$p_{g, f}(q) M_g$ for polynomials 
$p_{g,f}(q) \in \Z[q,q^{-1}]$ and
$g \in \Z^{m|n}_{\leq d}$ with $g \prec f$.
Now we make corrections to the first approximation.
Let $g \prec f$ be maximal such that
$p_{g,f}(q) \notin q \Z[q]$.
Let $p_{g,f}'(q)$ be the unique bar invariant element of $\Z[q,q^{-1}]$
such that $p_{g,f}'(q) \equiv p_{g, f}(q) \pmod{q \Z[q]}$.
Proceeding by induction on the ordering on $\Z^{m|n}_{\leq d}$, 
we may assume that $T_g^{(d)}$ is already known.
Subtract $p_{g,f}'(q) T_g^{(d)}$ from the first approximation, 
to obtain a second approximation to $T_f^{(d)}$.
Repeating the correction procedure, we reduce in finitely many steps
to a bar invariant expression 
that equals $M_f$ plus a $q \Z[q]$-linear combination of $M_g$'s.
This must be $T_f^{(d)}$ by the uniqueness.
We are done.

\begin{Example}\rm
We explain how to compute $T_{(0,4,1|0,2,3)}^{(4)}$ using the algorithm.
The bumping procedure proceeds as follows:
\begin{align*}
(0,4,1|0,2,3) 
&\stackrel{a_1 = 0}{\longrightarrow} (0,4,1|\underline{0},2,3)
\stackrel{a_2 = 1}{\longrightarrow} (0,4,\underline{1}|\underline{0},2,3)\\&
\stackrel{a_3 = 2}{\longrightarrow} (0,4,\underline{1}|\underline{0},\underline{2},\underline{3})
\stackrel{a_4 = 4}{\longrightarrow} (0,\underline{4},\underline{1}|\underline{0},\underline{2},\underline{3}).
\end{align*}
Now, $(0,5,2|1,3,4)$ is typical, so we can compute 
$\pi_{\leq 4} (F_4 E_3 E_2 F_1 E_0 T^{(5)}_{(0,5,2|1,3,4))})$
using the Kazhdan-Lusztig algorithm.
It turns out that this equals $M_{(0,4,1|0,2,3)} + 
M_{(1,4,1|1,2,3)}$ plus a $q \Z[q]$-linear combination of lower terms.
Now one computes $T^{(4)}_{(1,4,1|1,2,3)}$ needed for the correction
procedure
by repeating the algorithm
(which is rather lengthy). Finally one obtains
\begin{multline*}
T^{(4)}_{(0,4,1|0,2,3)} = 
M_{(0,4,1|0,2,3)} + q M_{(1,4,0|0,2,3)} + q M_{(4,0,1|0,2,3)}\\
+ q^2 M_{(1,4,1|1,2,3)} + q^2 M_{(4,1,0|0,2,3)} + q^3 M_{(4,1,1|1,2,3)}.
\end{multline*}
Note there is no reason why we chose to start
the bumping procedure with $f_+$ in describing the algorithm. 
One could also start the bumping procedure with $f_-$, increasing
all entries in $f_-$ equal to $a_1$ by $1$ and so on \dots.
In practice, one should always choose to start with 
the side for which the resulting 
word $X_N \dots X_1 \in \U$ is as short as possible.
In the present example, it is better to start the bumping procedure
with $f_-$, since then there is only one step:
$$
(0,4,1|0,2,3) \stackrel{a_1 = 0}{\longrightarrow} (\underline{0},4,1|0,2,3).
$$
Thus, we need to compute
$\pi_{\leq 4} (F_0 T^{(4)}_{(1,4,1|0,2,3)})$
instead, which is much quicker as only one generator of $\U$ needs to be 
applied. It turns out that this equals $T^{(4)}_{(0,4,1|1,2,3)}$ directly
(indeed it already equals $T_{(0,4,1|1,2,3)}$), 
with no corrections needed.
\end{Example}

Computer calculations using the above algorithm support the following
positivity conjecture:

\begin{Conjecture}\label{pc} Let $f \in \Z^{m|n}$.
\begin{itemize}
\item[(i)]
The coefficients $t_{g, f}(q)$ of $T_f$ 
when expanded in the basis 
$\{M_g\}_{g \in \Z^{m|n}}$
belong to $\N[q]$.
\item[(ii)]
The coefficients $l_{g, f}(q)$ of $L_f$ 
when expanded in the basis 
$\{M_g\}_{g \in \Z^{m|n}}$
belong to $\N[-q^{-1}]$.
\item[(iii)]
For each $a \in \Z$, $f \in \Z^{m|n}$ and $r \geq 1$,
the coefficients of $E_a^{(r)} T_f$ and $F_a^{(r)} T_f$
when expanded in the basis $\{T_g\}_{g \in \Z^{m|n}}$
belong to
$\N[q,q^{-1}]$.
\item[(iv)]
For each $a \in \Z$, $f \in \Z^{m|n}$ and $r \geq 1$,
the coefficients of $E_a^{(r)} L_f$ and $F_a^{(r)} L_f$
when expanded in the basis $\{L_g\}_{g \in \Z^{m|n}}$
belong to
$\N[q,q^{-1}]$.
\end{itemize}
\end{Conjecture}

If this positivity conjecture is true, it follows in particular that
each $T_f$ belongs to $\T^{m|n}$ rather than the completion
$\widehat{\T}^{m|n}$, i.e. each $T_f$ is a finite linear combination 
of $M_g$'s.
To see this, we modify the above algorithm
to obtain an algorithm that computes $T_f$ itself (not just the
$T_f^{(d)}$'s) in finitely many steps, as follows.
To start with, one follows the bumping procedure to obtain
$h$ and the elements $X_1, \dots, X_N \in \U$ exactly as above.
Since $\# h < \# f$, we may assume that $T_h$ is known 
inductively and is a finite sum of $M_g$'s.
Consider $X_N \dots X_1 T_h \in \T^{m|n}$, 
and choose $e$ to be minimal so that
$X_N \dots X_1 T_h \in \T^{m|n}_{\leq e}$.
In view of Conjecture~\ref{pc}(iii), 
$X_N \dots X_1 T_h$ equals $T_f$ plus a 
$\N[q,q^{-1}]$-linear combination of $T_g$'s.
So by Conjecture~\ref{pc}(i), 
we must have that that $T_f \in \T^{m|n}_{\leq e}$, hence
$T_f = T_f^{(e)}$.
Now follow the above algorithm to compute $T_f^{(e)}$.

\Point{Crystal structures}\label{s6}
Finally in this section, we 
review some results of Kashiwara, see e.g. \cite{Ka} for the basic
language used here.
Let $\A$ be the subring of $\Q(q)$ consisting of rational functions
having no pole at $q = 0$.
Evaluation at $q = 0$ induces an isomorphism $\A / q \A \rightarrow \Q$.

Let $\V_\A$ be the $\A$-lattice in $\V$ spanned by the $v_a$'s,
and let $\W_\A$ be the $\A$-lattice in $\W$ spanned by the $w_a$'s.
Then, $\V_\A$ together with the basis of the $\Q$-vector space
$\V_\A / q \V_\A$
given by the images of the $v_a$'s is a lower {crystal basis}
for $\V$ at $q = 0$ in the sense of \cite[4.1]{Ka}.
Similarly, $\W_\A$ together with the basis for $\W_\A / q \W_\A$
given by the images of the $w_a$'s is a lower crystal basis for $\W$ at $q = 0$.
Let
$\T^{m|n}_\A = \W_\A^{\otimes m} \otimes_\A \V_\A^{\otimes n}$
be the $\A$-lattice in $\T^{m|n}$ spanned by the $M_f$'s.
Then, by \cite[Theorem 4.1]{Ka}, $\T^{m|n}_\A$
together with the basis for
$\T^{m|n}_\A / q \T^{m|n}_\A$
given by the images of the $M_f$ for
$f \in \Z^{m|n}$ is a lower crystal basis for $\T^{m|n}$ at $q = 0$.
Moreover, we can easily describe the associated crystal graph
using Kashiwara's tensor product rule.

To do this, let us identify
the set 
$\{M_f + q \T^{m|n}_\A\}_{f \in \Z^{m|n}}$ underlying the crystal basis
with the set $\Z^{m|n}$ in the obvious way.
Then, Kashiwara's crystal operators induce maps
$\tilde E'_a, \tilde F'_a: 
\Z^{m|n} \rightarrow \Z^{m|n} \sqcup \{\varnothing\}$.
(We are using $\tilde E_a', \tilde F_a'$ because $\tilde E_a, \tilde F_a$ 
are used for something else later on.)
Fix $a \in \Z$ and $f \in \Z^{m|n}$.
The {\em $a$-signature}
$(\sigma_{-m}, \dots, \sigma_{-1},\sigma_1, \dots, \sigma_n)$
of $f$ is defined by
\begin{equation}\label{sigdef}
\sigma_i = \left\{
\begin{array}{ll}
+&\hbox{if $i > 0$ and $f(i) = a$, or if $i < 0$ and $f(i) = a+1$,}\\
-&\hbox{if $i > 0$ and $f(i) = a+1$, or if $i < 0$ and $f(i) = a$,}\\
0&\hbox{otherwise.}
\end{array}
\right.
\end{equation}
From this, we form the {\em reduced $a$-signature}
by successively replacing 
subsequences of the form $+-$ (possibly separated by $0$'s)
in the signature 
with $0 0$
until no $-$ appears to the right of a $+$.
Recall the definition of $d_j \in \Z^{m|n}$ from \ref{s1}.
We define
$$
\tilde E'_a (f) = \left\{
\begin{array}{ll}
\varnothing&\hbox{if there are no $-$'s in the reduced $a$-signature,}\\
f-d_j&\hbox{if the rightmost $-$ is in position $j \in I(m|n)$,}
\end{array}
\right.
$$
and
$$
\tilde F'_a (f) = \left\{
\begin{array}{ll}
\varnothing&\hbox{if there are no $+$'s in the reduced $a$-signature,}\\
f+d_j&\hbox{if the leftmost $+$ is in position $j \in I(m|n)$.}
\end{array}
\right.
$$
Also let
\begin{align*}
\eps'_a(f) &= \max\{r \geq 0\:|\:(\tilde E_a')^r(f) \neq 0\}\\
&=
\hbox{the total number of $-$'s in the reduced $a$-signature},\\
\phi'_a(f) &= \max\{r \geq 0\:|\:(\tilde F_a')^r(f) \neq 0\}\\
&=
\hbox{the total number of $+$'s in the reduced $a$-signature}.
\end{align*}
Then, the datum $(\Z^{m|n}, \tilde E'_a, \tilde F'_a, \eps'_a, \phi'_a, \wt)$
is the crystal associated to the module $\T^{m|n}$.

\begin{Example}\rm\label{eg1}
Consider the function
$f = (3,6,2,0,2,1|3,2,1) \in \Z^{6|3}$.
The $2$-signature is
$(+,0,-,0,-,0|-,+,0)$.
Cancelling off $+-$ pairs, we deduce that the reduced $2$-signature is
$(0,0,0,0,-,0|-,+,0)$.
Hence, the $2$-string through
$f$ in the crystal graph is
\begin{align*}
(3,6,2,0,3,1|2,2,1)
\stackrel{\tilde F'_2}{\longrightarrow}
(3,6,2,0,2,1|2,2,1)&\stackrel{\tilde F'_2}{\longrightarrow}
(3,6,2,0,2,1|3,2,1)\\
&\stackrel{\tilde F'_2}{\longrightarrow}
(3,6,2,0,2,1|3,3,1).
\end{align*}
\end{Example}

\begin{Theorem}\label{thmc}
Let $f \in \Z^{m|n}$ and $a \in \Z$.
\begin{itemize}
\item[(i)]
$E_a T_f = 
[\phi_a'(f)+1] T_{\tilde E_a'(f)} + 
\widehat{\sum}_{g \in \Z^{m|n}}
u_{g,f}^a T_g$
where the coefficient $u_{g,f}^a$ 
belongs to $q^{2 - \phi'_a(g)} \Z[q]$
and is zero unless
$\eps'_b(g) \geq \eps'_b(f)$ for all $b \in \Z$.
\item[(ii)]
$F_a T_f = [\eps'_a(f)+1] T_{\tilde F'_a(f)} + 
\widehat{\sum}_{g \in \Z^{m|n}}
v_{g,f}^a T_g$
where the coefficient $v_{g,f}^a$ belongs to $q^{2 - \eps'_a(g)} \Z[q]$
and is zero unless
$\phi'_b(g) \geq \phi'_b(f)$ for all $b \in \Z$.
\end{itemize}
(In (i) resp. (ii), the first term on the right hand side 
should be omitted if
$\tilde E'_a(f)$ resp. 
$\tilde F'_a(f)$ equals $\varnothing$.)
\end{Theorem}

\begin{proof}
Fix $d \in \Z$, and consider
$\T = \T^{m|n}_{\leq d}$, which is an integrable 
module in the sense of \cite[1.3]{KaG} 
with respect to the subalgebra of $\U$
generated by all $E_a, F_a, K_b^{\pm 1}$ 
for $a < d, b \leq d$.
Let $\T_0$ (resp. $\T_{\Q}$) be the $\A$ (resp.
$\Q[q,q^{-1}]$) -lattice in $\T$ spanned
by the basis elements $\{M_f\}_{f \in \Z^{m|n}_{\leq d}}$. 
Let $\T_{\infty} = \pi_{\leq d} (\overline{\T_0})$, 
an $\overline{\A}$-lattice in $\T$.
The canonical map
$\T_{\Q} \cap \T_{0} \cap \T_{\infty} 
\rightarrow \T_{0} / q \T_0$
is an isomorphism; this follows at once from
Lemma~\ref{lemmb} since that shows that 
all three lattices are generated by the elements 
$\{T_f^{(d)}\}_{f \in \Z^{m|n}_{\leq d}}$.
The preimage of 
the crystal basis element 
$M_f + q \T_0$ is $T_f^{(d)}$,
for each $f \in \Z_{\leq d}^{m|n}$.
In the language of Kashiwara \cite{KaG}, this shows that
$(\T_\Q, \T_0, \T_{\infty})$ 
is a balanced triple, and that
$\{T_f^{(d)}\}_{f \in \Z^{m|n}_{\leq d}}$ 
is a lower global crystal basis 
for $\T$ at $q = 0$.

Now we get from
\cite[Proposition 5.3.1]{KaG} (which is about upper global crystal bases)
and an argument involving duality \cite[$\S$3.2]{KaG}, that
$$
E_a T_f^{(d)} = 
[\phi'_a(f)+1] T_{\tilde E'_a(f)}^{(d)} + 
\sum_{g \in \Z_{\leq d}^{m|n}}
u_{g,f}^a T_g^{(d)}
$$
where the coefficient $u_{g,f}^a$ belongs to $q^{2 - \phi'_a(g)} \Z[q]$
and is zero unless
$\eps'_b(g) \geq \eps'_b(f)$ for all $b < d$.
Taking the limit as $d \rightarrow \infty$, we get (i).
The proof of (ii) is similar.
\end{proof}

We will also meet certain dual crystal operators on $\Z^{m|n}$.
Define
\begin{align}\label{dc1}
\tilde E_a^*(f) := -\tilde F'_{-1-a} (-f),\qquad
&\tilde F_a^*(f) := -\tilde E'_{-1-a} (-f),\\
\eps_a^*(f) :=  \phi'_{-1-a}(-f),
\qquad
&\phi_a^*(f) :=  \eps'_{-1-a}(-f).\label{dc2}
\end{align}
These can be described explicitly in a similar way to the above:
for fixed  $a \in \Z$ and $f \in \Z^{m|n}$, let
$(\sigma_{-m}, \dots, \sigma_{-1}, \sigma_1, \dots, \sigma_n)$ 
be the $a$-signature
as defined in (\ref{sigdef}).
Form the {\em dual reduced $a$-signature} 
by successively replacing 
sequences of the form
$-+$ (possibly separated by $0$'s) with $00$
until no $-$ appears to the left of a $+$.
Then:
$$
\tilde E^*_a (f) = \left\{
\begin{array}{ll}
0&\hbox{if there are no $-$'s in the dual reduced $a$-signature,}\\
f-d_j&\hbox{if the leftmost $-$ is in position $j \in I(m|n)$,}
\end{array}
\right.
$$
and
$$
\tilde F_a^* (f) = \left\{
\begin{array}{ll}
0&\hbox{if there are no $+$'s in the dual reduced $a$-signature,}\\
f+d_j&\hbox{if the rightmost $+$ is in position $j \in I(m|n)$.}
\end{array}
\right.
$$
Also
\begin{align*}
\eps^*_a(f) &=
\hbox{the total number of $-$'s in the dual reduced $a$-signature},\\
\phi_a^*(f) &=
\hbox{the total number of $+$'s in the dual reduced $a$-signature}.
\end{align*}
In this way, we obtain the dual crystal structure
$(\Z^{m|n}, \tilde E^*_a, \tilde F^*_a, \eps_a^*, 
\phi_a^*, \wt)$
on the underlying set $\Z^{m|n}$.

\begin{Theorem}\label{cthm}
Let $f \in \Z^{m|n}$ and $a \in \Z$.
\begin{itemize}
\item[(i)]
$E_a L_f = [\eps_a^*(f)] L_{\tilde E_a^* (f)} + 
\widehat{\sum}_{g \in \Z^{m|n}}
w_{g,f}^a L_g$
where the coefficient $w_{g,f}^a$ belongs to $q^{2 - \eps_a^*(f)} \Z[q]$
and is zero unless
$\phi_b^*(g) \leq \phi_b^*(f)$ for all $b\in\Z$.
\item[(ii)]
$F_a L_f = [\phi_a^*(f)] L_{\tilde F_a^* (f)} + 
\widehat{\sum}_{g \in \Z^{m|n}}
x_{g,f}^a L_g$
where the coefficient $x_{g,f}^a$ belongs to $q^{2 - \phi_a^*(f)} \Z[q]$
and is zero unless
$\eps_b^*(g) \leq \eps_b^*(f)$ for all $b\in\Z$.
\end{itemize}
\end{Theorem}

\begin{proof}
Dualize Theorem~\ref{thmc} using Theorem~\ref{thmd} and
Lemma~\ref{sf}. 
\end{proof}

\section{Exterior algebra}

Now we descend from the tensor space $\T^{m|n}$ to 
$\E^{m|n}$. We continue with the same notation as in section 2.

\Point{\boldmath The space $\E^{m|n}$}\label{thespace}
Let $w_0$ denote the longest element of $S_{m|n}$. Let
\begin{equation}
H_0 := 
\sum_{x \in S_{m|n}} (-q)^{\ell(x) - \ell(w_0)} H_x \in \H_{m|n}.
\end{equation}
The first lemma summarizes some elementary properties.

\begin{Lemma}\label{hprops}
The following properties hold:
\begin{itemize}
\item[(i)] $H_i H_0 = -q H_0 = H_0 H_i$ for any $i \in I(m-1|n-1)$;
\item[(ii)] $\overline{H_0} = H_0$;
\item[(iii)] $H_0^2 = -[m]![n]! H_0$;
\item[(iv)] $H_0 = \tau(H_0)$;
\item[(v)] the map $\omega:\H_{m|n} \rightarrow \H_{n|m}$ maps
$H_0  \in \H_{m|n}$ to $H_0 \in \H_{n|m}$.
\end{itemize}
\end{Lemma}

\begin{proof}
Part (i) is an easy exercise.
For (ii), use \cite[Proposition 2.9]{soergel}
and apply the map $dia$ there.
For (iii), one gets at once using (i) that
$$
H_0^2 = \sum_{x \in S_{m|n}} (-q)^{\ell(x) - \ell(w_0)} H_x H_0
=
\sum_{x \in S_{m|n}} (-q)^{2 \ell(x) - \ell(w_0)} H_0.
$$
Now use the well-known formula for the Poincar\'e polynomial
of $S_{m|n}$ to rewrite the sum.
Finally, (iv) and (v) are obvious.
\end{proof}

Let
$\E^{m|n} := \T^{m|n} H_0,$
a $\U$-submodule of $\T^{m|n}$.
Note $\E^{m|n}$ is the $q$-analogue of the exterior power
${\bigwedge}^m \W \otimes {\bigwedge}^n \V.$
Form the completion $\widehat{\E}^{m|n} = \widehat{\T}^{m|n} H_0$
as in \ref{s2}.
By Lemmas~\ref{omega}(i) and \ref{hprops}(v), the restriction of the map 
$\omega:\widehat{\T}^{m|n}\rightarrow\widehat{\T}^{n|m}$ 
is an isomorphism
$\omega:\widehat{\E}^{m|n}\rightarrow\widehat{\E}^{n|m}$.

We will call $f \in \Z^{m|n}$ {\em dominant} if
$f(-m) < \dots < f(-1), f(1) > \dots > f(n).$
We warn the reader that the inequality signs here are strict, unlike in the
earlier definition of antidominant!
Let $\Z^{m|n}_+$ denote the set of all dominant $f \in \Z^{m|n}$.
For $f \in \Z_+^{m|n}$, let
\begin{equation}\label{kf}
K_f := M_{f \cdot w_0} H_0 \in \E^{m|n}.
\end{equation}
The following lemma implies that 
the $\{K_f\}_{f \in \Z_+^{m|n}}$
form a basis for $\E^{m|n}$.

\begin{Lemma}\label{kl}
Suppose $f \in \Z^{m|n}$ and let 
$x \in S_{m|n}$ be the unique element of minimal length such that
$f \cdot x$ is antidominant.
Then,
$$
M_f H_0 = \left\{
\begin{array}{ll}
(-q)^{\ell(x)}K_{f\cdot x w_0}
&\hbox{if $f \cdot xw_0$ is dominant,}\\
0&\hbox{otherwise.}
\end{array}
\right.
$$
\end{Lemma}

\begin{proof}
We have that $M_f = M_{f \cdot x} H_{x^{-1}}$.
So applying Lemma~\ref{hprops}(i),
$M_f H_0 = M_{f \cdot x} H_{x^{-1}} H_0
= (-q)^{\ell(x)} M_{f \cdot x} H_0$.
Finally, note that if $f \cdot xw_0$ is not dominant,
then $M_{f \cdot x} H_0 = 0$.
\end{proof}

\Point{\boldmath Canonical bases}
Since $\overline{H_0} = H_0$ by Lemma~\ref{hprops}(ii), 
the bar involution on $\widehat\T^{m|n}$ leaves
$\widehat{\E}^{m|n}$ invariant.
Moreover, for dominant $f$, 
$\overline{K_f} = \overline{M_{f \cdot w_0}} H_0$.
So using Lemma~\ref{kl} and 
the explicit description of $\overline{M_{f \cdot w_0}}$
given by (B1) and (B3) in the proof of Theorem~\ref{thma}, we see that
$\overline{K_f} = K_f + (*)$
where $(*)$ is a (possibly infinite) $\Z[q,q^{-1}]$-linear
combination of $K_g$'s for $g \in \Z^{m|n}_+$ with $g \prec f$.
Moreover, for typical dominant $f$, 
we have that $\overline{K_f} = K_f$.
As in Theorem~\ref{thma}, these properties uniquely
characterize the bar involution on $\widehat\E^{m|n}$:

\begin{Theorem}
There exists a unique continuous, antilinear map
 $-:\widehat{\E}^{m|n} \rightarrow
\widehat{\E}^{m|n}$
such that
\begin{itemize}
\item[(i)]
$\overline{K_f} = K_f$ for all typical
$f \in \Z_+^{m|n}$;
\item[(ii)]
$\overline{X u} = \overline{X}\, \overline{u}$
for all $X \in \U$ and $u \in \widehat{\E}^{m|n}$.
\end{itemize}
Moreover, 
\begin{itemize}
\item[(iii)] bar is an involution;
\item[(iv)] $\overline{K_f} = K_f + (*)$ where $(*)$ is
a (possibly infinite) $\Z[q,q^{-1}]$-linear combination
of $K_g$'s for dominant $g \prec f$;
\item[(v)] $\overline{\omega (u)} = \omega(\overline{u})$ 
for all $u \in\widehat{\E}^{m|n}$.
\end{itemize}
\end{Theorem}

\begin{proof}
We have already proved existence above. For uniqueness, note on applying
$H_0$ to the conclusion of Theorem~\ref{triv} that we can write each
$K_f$ as a possibly infinite linear combination of $X_g K_g$'s 
for $X_g \in \U$ and typical $g \in \Z^{m|n}_+$.
Hence, as in Corollary~\ref{swr}, the only continuous 
$\U$-endomorphism of $\widehat{\E}^{m|n}$ that fixes $K_f$
for all typical $f \in \Z^{m|n}_+$ is the identity map.
Now using this one gets uniqueness by exactly the same argument as in the
proof of Theorem~\ref{thma}.
\end{proof}

Now applying the general principles 
used in the proof of Theorem~\ref{thmb}, 
we deduce:

\begin{Theorem}\label{thmB}
There exist unique topological bases
$\{U_f\}_{f \in \Z_+^{m|n}}$,
$\{L_f\}_{f \in \Z_+^{m|n}}$ for
$\widehat{\E}^{m|n}$ such that
\begin{itemize}
\item[(i)]
$\overline{U_f} = U_f$
and 
$\overline{L_f} = L_f$;
\item[(ii)]
$U_f \in K_f + \widehat{\sum}_{g \in \Z_+^{m|n}} q \Z[q] K_g$
and
$L_f \in K_f + \widehat{\sum}_{g \in \Z_+^{m|n}} q^{-1} 
\Z[q^{-1}] K_g$.
\end{itemize}
Moreover, 
\begin{itemize}
\item[(iii)]
$U_f = K_f + (*)$ and $L_f = K_f + (**)$ where
$(*)$ and $(**)$ are (possibly infinite) linear combinations of 
$K_g$'s for dominant $g \prec f$;
\item[(iv)] $\omega(U_f) = U_{\omega(f)}$ and
$\omega(L_f) = L_{\omega(f)}$.
\end{itemize}
\end{Theorem}

We use the following notation for the coefficients:
\begin{equation}\label{tde2}
L_f = \sum_{g \in \Z^{m|n}_+} l_{g, f}(q) K_g,
\qquad
U_f = \sum_{g \in \Z^{m|n}_+} u_{g, f}(q) K_g,
\end{equation}
for polynomials $l_{g, f}(q) \in \Z[q^{-1}], u_{g, f}(q) \in \Z[q]$.
We know that $l_{g, f}(q) = u_{g, f}(q) = 0$ unless $g \preceq f$, and
$l_{f, f}(q) = u_{f, f}(q) = 1$.

Note that $K_f = M_f+(*)$ where $(*)$ is a $q^{-1} \Z[q^{-1}]$-linear
combination of $M_g$'s.
So the element $L_f$ defined in Theorem~\ref{thmB} is bar invariant
and equals $M_f$ plus a $q^{-1} \Z[q^{-1}]$-linear combination of $M_g$'s.
So by the uniqueness in Theorem~\ref{thmb}, the elements $L_f$ and the
polynomials $l_{g, f}(q)$ defined
here are {\em the same} as the ones defined in \ref{s4}, for
dominant $g,f$. Thus our notation is consistent with the earlier notation.
The relationship between the elements $U_f$ here and the $T_f$'s from
before is given by:

\begin{Lemma}\label{ud} For $f \in \Z^{m|n}_+$, 
$U_f = T_{f \cdot w_0} H_0$.
\end{Lemma}

\begin{proof}
Note that $T_{f \cdot w_0} H_0$ is a bar invariant element of $\widehat{\E}^{m|n}$.
Recall that $T_{f \cdot w_0}$ equals $M_{f \cdot w_0}$ plus a $q \Z[q]$-linear
combination of $M_g$'s.
So applying Lemma~\ref{kl}, $T_{f \cdot w_0} H_0$ 
equals $K_f$ plus a  $q \Z[q]$-linear combination of $K_g$'s.
Hence $T_{f \cdot w_0} H_0 = U_f$ by the uniqueness in Theorem~\ref{thmB}.
\end{proof}

\Point{Duality}
Recall the antilinear involution $\sigma:\T^{m|n}\rightarrow \T^{m|n}$
defined in (\ref{sdef}).
In view of Lemmas~\ref{omega}(iii) and \ref{hprops}(ii), 
this leaves the subspace $\E^{m|n}$ invariant. 
Indeed, by Lemma~\ref{kl}, we have that
\begin{equation}\label{srel}
\sigma(K_f) =  (-q)^{\ell(w_0)}K_{-f \cdot w_0}
\end{equation}
for each $f \in \Z^{m|n}_+$.
Let $(.,.)_\EE$ be the bilinear form on $\E^{m|n}$ 
defined so that the elements $\{K_f\}_{f \in \Z^{m|n}_+}$ are
orthonormal.
Note by Lemmas~\ref{omega}(ii) and \ref{hprops}(iii) that
\begin{align*}
(K_f, K_g)_{\TT} &= (M_{f \cdot w_0} H_0, M_{g \cdot w_0} H_0)_{\TT}
=
-[m]![n]! (M_{f \cdot w_0} H_0, M_{g \cdot w_0})_{\TT}\\
&= 
-(-q)^{-\ell(w_0)} [m]![n]! \delta_{f,g}.
\end{align*} 
Hence,
\begin{equation}\label{frel}
(u,v)_\EE = - \frac{(-q)^{\ell(w_0)}}{[m]![n]!} (u,v)_{\TT}
\end{equation}
for all $u,v \in \E^{m|n}$.
Finally, define a bilinear form $\langle.,.\rangle_\EE$
on $\widehat{\E}^{m|n}$ by setting
\begin{equation}
\langle u, v \rangle_\EE := (-q)^{-\ell(w_0)}(u, \sigma(\overline{v}))_\EE
\end{equation}
for all $u,v \in \widehat{\E}^{m|n}$. 
Comparing with the definition of the form $\langle.,.\rangle_\TT$
from (\ref{od}) and using (\ref{frel}), one sees
immediately that
\begin{equation}
\langle u, v \rangle_\EE = 
-
\frac{1}{[m]![n]!} \langle u,v \rangle_\TT
\end{equation}
for all $u,v \in \widehat{\E}^{m|n}$. Hence in particular
we get from Lemma~\ref{sym} that the form $\langle.,.\rangle_\EE$ is symmetric.

\begin{Theorem}\label{dthm} For $f, g \in \Z^{m|n}_+$, 
$\langle L_f, U_{-g \cdot w_0} \rangle_\EE = \delta_{f,g}$
\end{Theorem}

\begin{proof}
Since $L_f \in \widehat{\E}^{m|n}$, we have by Lemma~\ref{hprops}(iii) that
$L_f H_0 = -[m]![n]! L_f$.
So applying Theorem~\ref{thmd} and Lemmas~\ref{ud}
 and \ref{sf}, we have that
\begin{align*}
\langle L_f, U_{-g \cdot w_0} \rangle_\EE &= 
-\frac{1}{[m]![n]!}
\langle L_f, T_{-g} H_0 \rangle_{\TT}
=
-\frac{1}{[m]![n]!}
\langle L_f H_0 , T_{-g} \rangle_{\TT}\\
&= \langle L_f, T_{-g} \rangle_{\TT}
= 
\delta_{f, g}.
\end{align*}
\end{proof}

By the theorem and the argument used to prove
Corollary~\ref{lcornew}, we get:

\begin{Corollary}\label{lcornew2} For $f \in \Z_+^{m|n}$,
$$
K_f = \sum_{g \in \Z_+^{m|n}} u_{-f\cdot w_0,-g\cdot w_0}(q^{-1}) L_g
=\sum_{g \in \Z_+^{m|n}} l_{-f\cdot w_0,-g\cdot w_0}(q^{-1}) U_g.
$$
\end{Corollary}

\Point{Crystal structures}\label{xs}
Next we describe the crystal structure on $\widehat\E^{m|n}$,
following  the same language as in \ref{s6}.
Recalling (\ref{kf}), let
\begin{align}\label{nc1}
\tilde E_a(f) := (\tilde E'_a (f \cdot w_0) ) \cdot w_0,\qquad
&\tilde F_a(f) := (\tilde F'_{a} (f\cdot w_0))\cdot w_0,\\
\eps_a(f) := \eps'_a(f \cdot w_0),\qquad&
\phi_a(f) := \phi'_a(f \cdot w_0)\label{nc2}
\end{align}
for $f \in \Z^{m|n}_+$. 
Then, 
$(\Z^{m|n}_+, \tilde E_a, \tilde F_a, \eps_a, \phi_a, \wt)$ 
is the crystal associated to the $\U$-module $\E^{m|n}$.
Actually, the crystal structure is so simple in this case, 
that we can actually list all the possibilities explicitly.
There are ten 
possible configurations for edges of color $a$ in the crystal graph, listed below. Here, $\dots$
denotes entries different from $a, a+1$.
\begin{itemize}
\item[(1)]
$(\dots,  a, a+1, \dots | \dots, a+1, a, \dots)$;
\item[(2)]
$(\dots,  a, a+1, \dots | \dots, a, \dots)
\stackrel{\tilde F_a}{\longrightarrow}(\dots,  a, a+1, \dots | \dots, a+1, \dots)$;
\item[(3)]
$(\dots,  a+1, \dots | \dots, a+1,a, \dots)
\stackrel{\tilde F_a}{\longrightarrow}(\dots,  a, \dots | \dots, a+1, a,\dots)$;
\item[(4)]
$(\dots,  a+1, \dots | \dots, a, \dots)
\stackrel{\tilde F_a}{\longrightarrow}
(\dots,  a, \dots | \dots, a,\dots)\\
\phantom{hello}\hspace{49mm}
\stackrel{\tilde F_a}{\longrightarrow}
(\dots,  a, \dots | \dots, a+1,\dots)$;
\item[(5)]
$(\dots,  a+1, \dots | \dots, a+1, \dots)$;
\item[(6)]
$(\dots,  a, a+1, \dots | \dots)$;
\item[(7)]
$(\dots|\dots,  a+1, a, \dots)$;
\item[(8)]
$(\dots|\dots,a,\dots) 
\stackrel{\tilde F_a}{\longrightarrow}
(\dots|\dots,a+1,\dots)$;
\item[(9)]
$(\dots,a+1,\dots|\dots) 
\stackrel{\tilde F_a}{\longrightarrow}
(\dots,a,\dots|\dots)$;
\item[(10)]
$(\dots| \dots)$.
\end{itemize}
We also have the {\em dual crystal} 
$(\Z^{m|n}_+, \tilde E_a^*, \tilde F_a^*, \eps_a^*, \phi_a^*, \wt)$,
where $\tilde E_a^*, \tilde F_a^*, \eps_a^*$ and $\phi_a^*$
are the restrictions of the functions from
\ref{s6} to $\Z^{m|n}_+$. Combining (\ref{dc1}) and (\ref{nc1}),
we have that
\begin{equation}\label{ndc1}
\tilde E_a^*(-f \cdot w_0) = - \tilde F_{-1-a}(f) \cdot w_0,
\qquad
\tilde F_a^*(-f \cdot w_0) = - \tilde E_{-1-a}(f) \cdot w_0.
\end{equation}
Again, there are ten possible configurations for the edges in the
corresponding dual crystal graph, all of which 
are exactly the same as (1)--(10) above (replacing
$\tilde F_a$ with $\tilde F_a^*$)
with the exception of (4) and (5) which change to
\begin{itemize}
\item[(4$^*$)]
$(\dots,  a+1, \dots | \dots, a, \dots)
\stackrel{\tilde F_a^*}{\longrightarrow}
(\dots,  a+1, \dots | \dots, a+1,\dots)\\
\phantom{hello}\hspace{49mm}
\stackrel{\tilde F_a^*}{\longrightarrow}
(\dots,  a, \dots | \dots, a+1,\dots)$;
\item[(5$^*$)]
$(\dots,  a, \dots | \dots, a, \dots)$.
\end{itemize}

\begin{Remark}\label{lrem}\rm
In (\ref{LLdef})
and Lemma~\ref{opprop}(v) below we will define 
mutually inverse bijections $\LL, \RR:\Z^{m|n}_+ \rightarrow \Z^{m|n}_+$.
By considering all the above cases (1)--(10) one by one, it is
not hard to check that $\LL$ satisfies, indeed is characterized
uniquely by, the following properties
\begin{itemize}
\item[(1)] if $f \in \Z^{m|n}_+$ is typical then
$\LL(f) = f$;
\item[(2)] 
for every $a \in \Z$ and $f \in \Z^{m|n}_+$, 
$\tilde E_a^* \LL(f) = \LL(\tilde E_af)$ and
 $\tilde F_a^* \LL (f) = \LL(\tilde F_af)$;
\item[(3)]  for every $f \in \Z^{m|n}_+$, $\wt(\LL (f)) = \wt(f)$.
\end{itemize}
Hence, 
$\LL:
(\Z^{m|n}_+, \tilde E_a, \tilde F_a, \eps_a, \phi_a, \wt)
\rightarrow
(\Z^{m|n}_+, \tilde E_a^*, \tilde F_a^*, \eps_a^*, \phi_a^*, \wt)$
is an isomorphism of crystals with inverse $\RR$.
\end{Remark}

The crucial observation to be made from the above description
of the crystal graph is that {\em all $a$-strings are of length
at most $2$}. The following lemma is a consequence of this particularly
simple structure.

\begin{Lemma}\label{thme} Let $f \in \Z^{m|n}_+$ and $a \in \Z$.
\begin{itemize}
\item[(i)] If $\eps_a(f) > 0$, then
$E_a U_f = [\phi_a(f)+1] U_{\tilde E_a (f)}$.
\item[(ii)]
If $\phi_a(f) > 0$, then
$F_a U_f = [\eps_a(f)+1] U_{\tilde F_a (f)}$.
\end{itemize}
\end{Lemma}

\begin{proof} We prove (i), (ii) being similar.
Dualizing Theorem~\ref{cthm}(i) using Theorem~\ref{dthm} and Lemma~\ref{sf}
(or by considering the effect of the 
Kashiwara operators directly and arguing as in the proof of
Theorem~\ref{thmc}) gives us that
$E_a U_f = 
[\phi_a(f)+1] U_{\tilde E_a(f)} + 
\widehat{\sum}_{g \in \Z_+^{m|n}}
y_{g,f}^a U_g$
where $y_{g,f}^a$ 
belongs to $q^{2 - \phi_a(g)} \Z[q]$
and is zero unless
$\eps_b(g) \geq \eps_b(f)$ for all $b \in \Z$.
Suppose that $y_{g,f}^a \neq 0$ for some $g$.
By assumption, $\eps_a(g) \geq \eps_a(f) \geq 1$,
so $\phi_a(g) \leq 1$ since all $a$-strings are of length $\leq 2$.
So $0 \neq y_{g,f}^a \in q \Z[q]$.
But $y_{g,f}^a$ is bar invariant, so this is a contradiction.
\end{proof}

\Point{\boldmath Two algorithms}\label{ta}
In this subsection, we describe algorithms to compute the
canonical basis $\{U_f\}$ and
the coefficients $\{l_{g,f}(q)\}$ of the
dual canonical basis of $\E^{m|n}$ explicitly.
The first algorithm computes $U_f$, and is similar to the algorithm for
computing the $T_f$'s explained in \ref{s5} --- but it is {\em much}
simpler since no corrections are needed thanks to Lemma~\ref{thme}.

\begin{Procedure}\label{alg2}\rm
Suppose we are given $f \in \Z^{m|n}_+$ with $\#f > 0$.
Compute $h \in \Z^{m|n}_+$ and operators
$X_a, Y_a
\in \{E_a, F_a\}_{a\in\Z}$ by following the instructions below starting
at step (0).
\begin{itemize}
\item[(0)\phantom{$'$}]
Choose the largest $i \in \{-m,\dots,-1\}$ such that $f(i) = f(j)$
for some $j \in \{1,\dots,n\}$. Go to step (1).
\item[(1)\phantom{$'$}] 
If $i < -1$ and $f(i+1) = f(i)+1$, replace $i$ by $(i+1)$
and repeat step (1). Otherwise, go to step (2).
\item[(2)\phantom{$'$}] If $f(i) + 1 = f(j)$ for some (necessarily unique)
$j \in \{1,\dots,n\}$
go to step (1)$'$. 
Otherwise, set $X_a = F_{f(i)}, Y_a = E_{f(i)}$ and 
$h = f - d_i$. Stop.
\item[(1)$'$]
If $j > 1$ and $f(j-1) = f(j)+1$, replace $j$ by $(j-1)$
and repeat step (1)$'$. Otherwise, go to step (2)$'$.
\item[(2)$'$] If $f(j)+1 = f(i)$ for some (necessarily unique)
$i \in \{-m,\dots,-1\}$
go to step (1). 
Otherwise, set $X_a = E_{f(j)}, Y_a = F_{f(j)}$ and 
$h = f + d_j$. Stop.
\end{itemize}
\end{Procedure}

The following lemma follows immediately from the nature of the
above procedure and Lemma~\ref{thme}.

\begin{Lemma}\label{props}
Take $f \in \Z^{m|n}_+$
with $\#f > 0$.
Define $h$ and operators $X_a,Y_a \in \{E_a,F_a\}_{a \in \Z}$
according to Procedure~\ref{alg2}.
Then, one of the following holds:
\begin{itemize}
\item[(i)] $\#h = \#f$.
In this case, the $a$-string through $f$ is
$h \stackrel{\tilde X_a}{\longrightarrow} f$, of length $1$.
Moreover, $X_a U_h = U_f$, $Y_a U_f = U_h$ and $X_a K_h = K_f$.
\item[(ii)] $\#h = \#f - 1$.
In this case, the $a$-string through $f$ is
$h \stackrel{\tilde X_a}{\longrightarrow} f \stackrel{\tilde X_a}{\longrightarrow} g$,
of length $2$.
Moreover $X_a U_h = U_f$, $Y_a U_f = [2] U_h$ and 
$X_a K_h = K_f + q K_{\tilde X_a^*(h)}$.
\end{itemize}
Case (ii) (when the atypicality gets strictly smaller) must occur
after at most $(m+n-1)$ repetitions of the procedure.
Hence after finitely many recursions, the 
procedure reduces $f$ to a {typical} weight.
\end{Lemma}

Lemma~\ref{props}
implies the following algorithm for computing $U_f$.
If $f \in \Z^{m|n}_+$ is typical then $U_f = K_f$,
since such $f$'s are minimal in the ordering $\preceq$ in $\Z^{m|n}_+$.
Otherwise, apply Procedure~\ref{alg2} to get $h \in \Z^{m|n}_+$ and $X_a
\in \{E_a, F_a\}_{a\in\Z}$.
Since the procedure always reduces $f$ to a typical weight in 
finitely many steps, we may assume $U_h$ is known recursively. Then
$U_f = X_a U_h$.

\begin{Example}\label{oldex}\rm Applying the algorithm repeatedly, 
we get that
\begin{align*}
U_{(0,1,3,4|2,1,0)} &= 
F_4 F_3 E_2 F_1 F_5 F_4 E_3 F_2 E_1 F_0
K_{(1,3,5,6|4,2,0)}\\
&= K_{(0,1,3,4|2,1,0)} + q K_{(1,3,4,6|6,2,1)} + q K_{(0,3,4,5|5,2,0)}
+ q^2 K_{(3,4,5,6|6,5,2)}.
\end{align*}
\end{Example}

In the next subsection, we will apply the above algorithm to derive
a closed formula for $U_f$.
We turn now to describing the second algorithm,
which computes the polynomials $l_{g,f}(q)$.
It will not be needed until \ref{klpol} below.
First we state a variation on Procedure~\ref{alg2}.

\begin{Procedure}\label{alg3}\rm
Suppose we are given $g \in \Z^{m|n}_+$ with $\#g > 0$.
Compute $h \in \Z^{m|n}_+$ and operators
$X_a, Y_a
\in \{E_a, F_a\}_{a\in\Z}$ by following the instructions below starting
at step (0).
\begin{itemize}
\item[(0)\phantom{$'$}]
Choose the smallest $i \in \{-m,\dots,-1\}$ such that $g(i) = g(j)$
for some $j \in \{1,\dots,n\}$. Go to step (1).
\item[(1)\phantom{$'$}] 
If $i > -m$ and $g(i-1) = g(i)-1$, replace $i$ by $(i-1)$
and repeat step (1). Otherwise, go to step (2).
\item[(2)\phantom{$'$}] If $g(i) - 1 = g(j)$ for some (necessarily unique)
$j \in \{1,\dots,n\}$
go to step (1)$'$. 
Otherwise, set $h = g + d_i$, $X_a = E_{h(i)}$ and $Y_a = F_{h(i)}$.
Stop.
\item[(1)$'$]
If $j < n$ and $g(j+1) = g(j)-1$, replace $j$ by $(j+1)$
and repeat step (1)$'$. Otherwise, go to step (2)$'$.
\item[(2)$'$] If $g(j)-1 = g(i)$ for some (necessarily unique)
$i \in \{-m,\dots,-1\}$
go to step (1). 
Otherwise, set $h = g - d_j$, $X_a = F_{h(j)}$ and $Y_a = E_{h(j)}$. Stop.
\end{itemize}
\end{Procedure}

\begin{Lemma}\label{secondalg} Suppose $g, f \in \Z^{m|n}_+$ with $\#g > 0$.
Define $h$ and operators $X_a,Y_a \in \{E_a,F_a\}_{a \in \Z}$
according to Procedure~\ref{alg3}.
Then,
$$
l_{g,f}(-q^{-1}) = 
\left\{
\begin{array}{ll}
l_{h, \tilde Y_a^*(f)}(-q^{-1}) &\hbox{if $\#h = \#g$,}\\
l_{h, \tilde Y_a^*(f)}(-q^{-1}) +q l_{\tilde X_a(h), f} (-q^{-1})
&\hbox{if $\#h = \#g-1$,}
\end{array}\right.
$$
interpreting $l_{h, \tilde Y^*_a (f)}(-q^{-1})$ as $0$ if
$\tilde Y^*_a (f) = \varnothing$.
\end{Lemma}

\begin{proof}
Let $f, g \in \Z^{m|n}_+$ with $\#f > 0$.
Apply Procedure~\ref{alg2} to construct $h$ and 
operators $X_a, Y_a \in \{E_a, F_a\}_{a \in \Z}$. 
Apply the operator $X_a$ to both sides of the equation
$$
K_h = \sum_{k \in \Z^{m|n}_+} l_{-h\cdot w_0, -k \cdot w_0}(q^{-1}) U_k
$$
from Corollary~\ref{lcornew2}. In the case
that $\#h = \#f$, we know by Lemma~\ref{props} that
$h$, hence also all $k \in \Z^{m|n}_+$ 
with the same weight as $h$, is at one end of an $a$-string of length 1.
So by Lemma~\ref{thme}, 
$X_a U_k = U_{\tilde X_a (k)}$ for all $k$ with $l_{-h\cdot w_0, -k\cdot w_0}(q^{-1}) \neq 0$.
Hence, 
$$
K_f = X_a K_h = 
\sum_{k \in \Z^{m|n}_+} l_{-h\cdot w_0, -k\cdot w_0}(q^{-1}) 
U_{\tilde X_a(k)}.
$$
On the other hand, if $\#h =\#f-1$, then $h$, hence also all $k$
of the same weight as $h$, lies at one end of an
$a$-string of length 2, so
$X_a U_k = U_{\tilde X_a (k)}$.
We also know from Lemma~\ref{props} that $X_a K_h = K_f + 
q K_{\tilde X_a^* (h)}$.
So in this case,
$$
K_f = \sum_{k \in \Z^{m|n}_+} l_{-h\cdot w_0, -k\cdot w_0}(q^{-1})
 U_{\tilde X_a(k)}
- q \sum_{h \in \Z^{m|n}_+} l_{-\tilde X_a^*(h) \cdot w_0, -g\cdot w_0}(q^{-1})
 U_g.
$$
Now compute the coefficient of $U_g$ in the above
expressions for $K_f$ using Corollary~\ref{lcornew2} again, to get
$$
l_{-f\cdot w_0,-g\cdot w_0}(q^{-1}) = 
\left\{
\begin{array}{ll}
l_{-h\cdot w_0, -\tilde Y_a (g)\cdot w_0}(q^{-1})
&\hbox{if $\#h = \#f$,}\\
l_{-h \cdot w_0, -\tilde Y_a (g)\cdot w_0}(q^{-1})
&\hbox{if $\#h = \#f-1$,}\\
\qquad- q l_{-\tilde X_a^*(h) \cdot w_0, -g \cdot w_0}(q^{-1})&
\end{array}
\right.
$$
interpreting $l_{-h\cdot w_0, -\tilde Y_a (g)\cdot w_0}(q)$ as $0$ if
$\tilde Y_a (g) = \varnothing$.
The lemma follows from this, replacing $f$ by 
$-g\cdot w_0$ and $g$ by $-f \cdot w_0$ and using (\ref{ndc1}),
since Procedure~\ref{alg3} is just Procedure~\ref{alg2}
twisted by the involution $f \mapsto -f \cdot w_0$.
\end{proof}

Now to compute $l_{g,f}(-q^{-1})$, 
we have that $l_{g,f}(-q^{-1}) = \delta_{g,f}$ if $g$ is typical,
and it is $0$ if $g \not\preceq f$. 
Otherwise, if $\#g > 0$ and $g \preceq f$, apply Procedure~\ref{alg3}
and
Lemma~\ref{secondalg} to write $l_{g,f}(-q^{-1})$ in terms
of $l_{h, \tilde Y_a^*(f)}(-q^{-1})$
and (in case $\#h = \#g-1$) $l_{\tilde X_a(h), f}(-q^{-1})$, and repeat.
This process terminates in finitely many steps,
because $h$ is closer to being typical than $g$
in the sense of Procedure~\ref{alg3}, and $\tilde X_a(h)$ 
is closer than $g$ to failing the condition $g \preceq f$.
Note this algorithm shows in particular that
$l_{g,f}(-q^{-1}) \in \N[q]$, as also follows from the explicit
description given in Corollary~\ref{clearer}(ii) below.

\Point{Combinatorial description of canonical bases}\label{cdes}
We now introduce some combinatorics
to enable us to write down closed formulae for the
canonical basis and dual canonical basis elements.
The material in this subsection was inspired originally by \cite{JZ}.
Suppose $f \in \Z^{m|n}$ is conjugate under the action of $S_{m|n}$
to an element of $\Z^{m|n}_+$.
We will denote this ``dominant conjugate'' of $f$ by $f^+$.
For
$-m \leq i < 0 < j \leq n$ with $f(i) = f(j)$,
let
\begin{align}\label{lij}
\LL_{i,j}(f) &:= f - a(d_i - d_j),\\\intertext{where 
$a$ is the smallest positive integer such that
$f - a(d_i - d_j)$ and all
$\LL_{k,l} (f) - a(d_i-d_j)$
for $i < k < 0 < l < j$ with $f(k) = f(l)$
are conjugate to elements of $\Z^{m|n}_+$. Similarly, let}
\RR_{i,j}(f) &:= f + b(d_i-d_j),\label{rij}
\end{align}
where $b$ is the smallest positive integer such that
$f + b(d_i - d_j)$ and all
$\RR_{k,l} (f) + b(d_i-d_j)$
for $-m \leq k < i, j < l \leq n$ with $f(k) = f(l)$
are conjugate to elements of $\Z^{m|n}_+$.

Now take $f \in \Z^{m|n}_+$. 
Let $r = \# f$ and 
 $-m \leq i_1 < \dots < i_r
< 0 < j_r < \dots < j_1 \leq n$ be the unique integers with
$f(i_s) = f(j_s)$ for each $s = 1,\dots,r$.
For a tuple $\theta = (\theta_1,\dots,\theta_r) \in \N^r$, 
$|\theta|$ denotes $\theta_1+\dots+\theta_r$. Let
\begin{align}\label{one2}
\LL_{\theta}(f) &= 
\left(\LL_{i_r,j_r}^{\theta_r}\circ\LL_{i_{r-1},j_{r-1}}^{\theta_{r-1}}\circ \dots \circ\LL_{i_1, j_1}^{\theta_1}(f)\right)^+,\\
\LL'_{\theta}(f) &= 
\left(\LL_{i_1,j_1}^{\theta_1}\circ\LL_{i_2,j_2}^{\theta_2}\circ \dots \circ\LL_{i_r, j_r}^{\theta_r}(f)\right)^+,\\
\RR_{\theta}(f) &= 
\left(\RR_{i_1,j_1}^{\theta_1}\circ\RR_{i_2,j_2}^{\theta_2}\circ \dots\circ \RR_{i_r, j_r}^{\theta_r}(f)\right)^+,\\
\RR'_{\theta}(f) &= 
\left(\RR_{i_r,j_r}^{\theta_r} \circ\RR_{i_{r-1},j_{r-1}}^{\theta_{r-1}} \circ\dots\circ \RR_{i_1, j_1}^{\theta_1}(f)\right)^+.\label{rrpdef}
\end{align}
Note that $\LL'_{\theta}(f), \LL_{\theta}(f) \preceq f
\preceq\RR_{\theta}(f), \RR'_{\theta}(f)$.
The operators $\LL_{\theta}$ and $\RR_{\theta}$ will only ever be used
for $\theta$ belonging to the set $\{0,1\}^r$.
In the special case that $\theta_1 = \dots = \theta_r = 1$, we let
\begin{equation}
\LL(f) := \LL_{\theta}(f),\label{LLdef}\qquad
\RR(f) := \RR_{\theta}(f).
\end{equation}
The following combinatorial lemma lists some 
elementary properties of the lowering and raising operators,
which follow immediately from the definition.

\begin{Lemma}\label{opprop} Let $f \in \Z^{m|n}_+$ and $r = \#f$.
\begin{itemize}
\item[(i)]
Suppose $\theta \in \N^r$ and let
$\varphi = (\theta_r,\dots,\theta_1)$. Then,
$\RR_{\theta}(-f\cdot w_0) = - \LL_{\varphi}(f) \cdot w_0$ and
$\RR'_{\theta}(-f\cdot w_0) = - \LL'_{\varphi}(f) \cdot w_0.$
\item[(ii)] The sets
$\{\LL_{\theta}(f)\}_{\theta \in \{0,1\}^r}$
and
$\{\RR_{\theta}(f)\}_{\theta \in \{0,1\}^r}$
contain $2^r$
distinct elements.
\item[(iii)]
Suppose $\theta \in \{0,1\}^r$ and 
let $\varphi = (1-\theta_r,\dots,1-\theta_1)$. Then,
$\LL_{\theta}(\RR(f)) = \RR_{{\varphi}}(f)$
and $\RR_{\theta}(\LL(f)) = \LL_{{\varphi}}(f)$.
In particular, taking $\theta_1=\dots=\theta_r = 1$,
the maps $\LL, \RR:\Z^{m|n}_+\rightarrow\Z^{m|n}_+$ are mutually inverse bijections.
\end{itemize}
\end{Lemma}

\begin{Example}\rm 
Take $f = (0,1,3,4|2,1,0)$ as in Example~\ref{oldex},
so $\#f = 2$.
Then we have that $\LL_{(0,0)}(f) = (0,1,3,4|2,1,0)$,
$\LL_{(1,0)}(f) = (1,3,4,6|6,2,1)$, $\LL_{(0,1)}(f) = (0,3,4,5|5,2,0)$
and $\LL_{(1,1)}(f) = (3,4,5,6|6,5,2)$.
Observe these are exactly the $K_g$'s
appearing in the expression for $U_f$ computed 
in Example~\ref{oldex}.
\end{Example}

The main theorem of the subsection is the following.

\begin{Theorem}\label{nm}
For $f \in \Z^{m|n}_+$ and $r = \#f$,
\begin{itemize}
\item[(i)]$\displaystyle
U_f =
\sum_{\theta \in \{0,1\}^{r}}
q^{|\theta|} K_{\LL_{\theta}(f)};$
\item[(ii)]$\displaystyle
K_f =
\sum_{\theta \in \N^{r}}
(-q)^{|\theta|} U_{\LL'_{\theta}(f)}.$
\end{itemize}
\end{Theorem}

\begin{proof}
(i) If $f$ is typical, then $U_f = K_f$ and there is nothing to prove.
So suppose that $\#f > 0$ and define $h$ and $X_a,Y_a \in \{E_a,F_a\}_{a \in \Z}$ according to Procedure~\ref{alg2}.
We may assume by induction that the result has already been established for 
$h$. Recalling Lemma~\ref{props}, we need to consider two cases.
In the first case $\#h = \#f$, we know that
$U_h = \sum_{\theta \in \{0,1\}^r} q^{|\theta|} K_{\LL_{\theta}(h)}.$
Applying $X_a$ to both sides, noting that
$X_a U_h = U_f$ and
that $X_a K_{\LL_{\theta}(h)} = K_{\LL_{\theta}(f)}$
for each $\theta$, gives the desired
conclusion. In the second case $\#h = \#f - 1$.
This time, we know that
$U_h = \sum_{\theta \in \{0,1\}^{r-1}} q^{|\theta|} K_{\LL_{\theta}(h)}.$
For each $\theta$, $\LL_{\theta}(h)$ here has the form
$(\dots, a+1, \dots|\dots, a,\dots)$, so
$$
X_a K_{\LL_{\theta}(h)} = K_{(\dots,a,\dots|\dots,a,\dots)}+q
K_{(\dots,a+1,\dots|\dots,a+1,\dots)}
= K_{\LL_{{\theta\cup 0}}(f)} + q K_{\LL_{{\theta\cup 1}}(f)},
$$
where $\theta\cup x$ denotes $(\theta_1,\dots,\theta_{r-1},x) \in \{0,1\}^r$.
So again we see on applying $X_a$ to both sides that
$U_f = X_a U_h = 
\sum_{\theta \in \{0,1\}^{r}}
 q^{|\theta|}
K_{\LL_{\theta}(f)}$.

(ii) To deduce this from (i), we will work 
in the free $\Z[q,q^{-1}]$-module
$\mathcal M^{m|n}$ 
on basis $\{[f]\}_{f \in \Z_+^{m|n}\cdot S_{m|n}}$, completed to a topological
$\Z[q,q^{-1}]$-module $\widehat{\mathcal M}^{m|n}$ exactly as in
\ref{s2} so that expressions of the form $[f] + $(a possibly infinite 
linear combination
of $[g]'s$ with $g \prec f$) make sense.
We define continuous linear maps
$U, K:\widehat{\mathcal M}^{m|n} \rightarrow \widehat{\E}^{m|n}$
by letting $U([f]) = U_{f^+}, K([f]) = K_{f^+}$.
These maps have the right inverses 
$U^{-1}, K^{-1}:\widehat{\E}^{m|n} \rightarrow \widehat{\mathcal M}^{m|n}$
with $U^{-1}(U_f) = [f], K^{-1}(K_f) = [f]$ for each $f \in \Z^{m|n}_+$.
Finally, define continuous linear operators
$\lambda_{i,j}:\widehat{\mathcal M}^{m|n} 
\rightarrow \widehat{\mathcal M}^{m|n}$
for each $-m \leq i < 0 < j \leq n$ by 
$$
\lambda_{i,j}([f]) = \left\{\begin{array}{ll}
[\LL_{i,j}(f)]&\hbox{if $f(i)=f(j)$,}\\
0&\hbox{if $f(i)\neq f(j)$,}
\end{array}\right.
$$
for each $f \in \Z^{m|n}_+\cdot S_{m|n}$.
Now consider the maps
\begin{align*}
P := K \circ \left(\stackrel{\rightarrow}{\prod_{-m \leq i < 0 <j\leq n}}
(1+q \lambda_{i,j}) \right)\circ U^{-1}&:\widehat{\E}^{m|n}\rightarrow\widehat{\E}^{m|n},\\
Q := U \circ \left(
\stackrel{\leftarrow}{\prod_{-m \leq i < 0 <j\leq n}}
\frac{1}{1+q \lambda_{i,j}}\right) \circ K^{-1}&:\widehat{\E}^{m|n} \rightarrow\widehat{\E}^{m|n},
\end{align*}
where $\stackrel{\rightarrow}{\prod}$ is taken in some ordering
with $i$'s decreasing and $j$'s increasing from left to right
and $\stackrel{\leftarrow}{\prod}$ is taken in the opposite ordering,
and $\frac{1}{1+q \la_{i,j}}$ denotes 
$(1 - q\la_{i,j}+q^2 \la_{i,j}- \dots)$.
By (i) and the definition of the operator $\LL_{i,j}$,
the map $P$ sends $U_f$ to $U_f$, so $P = \operatorname{id}$.
On the other hand the result we are trying to prove
is equivalent to the statement that $Q$ sends $K_f$ to $K_f$. 
Therefore we will be done if we can show that
$P \circ Q = \operatorname{id}$, i.e. that for every $f \in \Z^{m|n}_+$,
$$
K 
\circ\left(\stackrel{\rightarrow}{\prod_{-m \leq i < 0 < j \leq n}}
(1+q\la_{i,j})\right) \circ (U^{-1}\circ U) \circ
\left(\stackrel{\leftarrow}{\prod_{-m \leq i < 0 < j \leq n}}
\frac{1}{1+q\la_{i,j}}\right) ([f]) = K_f.
$$
This is obvious if we can show that the 
inside map $(U^{-1}\circ U):[g] \mapsto [g^+]$ on the left hand side
can be omitted. For this, we check that
\begin{equation}\label{wewant}
K \circ
\left(
\stackrel{\rightarrow}{\prod_{-m \leq i < 0 <j\leq n}} (1+q \lambda_{i,j}) 
\right)
([g^+])
=
K \circ\left(
\stackrel{\rightarrow}{\prod_{-m \leq i < 0 <j\leq n}} (1+q \lambda_{i,j}) \right)
([g])
\end{equation}
for every 
$g \in \Z^{m|n}$ such that $[g]$ is involved in
$
\left(\stackrel{\leftarrow}{\prod}
\frac{1}{1+q \lambda_{i,j}} \right)([f])$ 
with non-zero coefficient.
Suppose we have such a $g$.
The crucial observation is that 
whenever there exist $-m \leq i' < i < 0 < j < j' \leq n$
with $g(j) = g(i) < g(i') = g(j')$, one can find
$c$ with $g(i) < c < g(i')$ that does not arise in the tuple $g$.
Given this it is not hard to see that (\ref{wewant}) holds.
\end{proof}

\begin{Corollary}\label{lfomr}
For $f \in \Z^{m|n}_+$,
\begin{itemize}
\item[(i)]
$K_f = \sum_g q^{-|\theta_g|} L_g$ where the
sum is over all $g \in \Z^{m|n}_+$
such that $\RR_{{\theta_g}}(g) = f$
for some (unique) $\theta_g \in \{0,1\}^{\#g}$;
\item[(ii)]
$L_f = \sum_{g, \theta}
(-q)^{-|\theta|} K_g$ where the sum is over all $g \in \Z^{m|n}_+$
and $\theta \in \N^{\#g}$ such that $\RR'_{\theta}(g) = f$.
\end{itemize}
\end{Corollary}

\begin{proof}
(i)
Recall from Corollary~\ref{lcornew2} that the 
coefficient of $L_g$ in $K_f$ is equal to 
$u_{-f \cdot w_0, -g \cdot w_0}(q^{-1})$.
By Theorem~\ref{nm}(i) and Lemma~\ref{opprop}(ii), 
$u_{-f\cdot w_0, -g\cdot w_0}(q^{-1})=
q^{-|\theta|}$ if $-f\cdot w_0= \LL_{\theta}(-g \cdot w_0)$ for some 
(necessarily unique) $\theta \in \{0,1\}^{\#g}$, and is zero otherwise.
Equivalently, invoking Lemma~\ref{opprop}(i), 
$u_{-f\cdot w_0, -g \cdot w_0}(q^{-1})= q^{-|\theta|}$
if $f = \RR_{{\theta}}(g)$ for some $\theta$, and is zero otherwise.

(ii) By Corollary~\ref{lcornew2} again,
$l_{g,f}(q^{-1})$ is equal to the coefficient of
$U_{-f\cdot w_0}$ in $K_{-g \cdot w_0}$.
By Theorem~\ref{nm}(ii), this
equals $\sum_{\theta} (-q)^{|\theta|}$
where the sum is over all $\theta \in \N^{\#g}$ with
$\LL'_{\theta}(-g \cdot w_0) = -f \cdot w_0$, equivalently,
$\RR'_{\theta}(g) = f$.
\end{proof}

\begin{Example}\label{goodie}\rm
Suppose $f = (-m,\dots,-2,-1|-1,-2,\dots,-n)$, so
$r = \#f = \min(m,n)$.
We observe that any $g \preceq f$ can be represented as $\LL_{\theta}'(f)$
for a unique element $\theta \in \N^r$ with
$\theta_1 \leq \dots \leq \theta_r$.
Moreover, this $\theta$ 
is also the unique element of $\N^r$
with the property that $f=\RR_{\theta}'(g)$.
We deduce from Corollary~\ref{lfomr}(ii) that
\begin{equation}
L_f = \sum_{\theta = (\theta_1 \leq \dots \leq \theta_r) \in \N^r}
(-q)^{-|\theta|} K_{\LL_{\theta}'(f)}.
\end{equation}
\end{Example}

Recalling the definitions from 
(\ref{tde2}), we can restate Theorem~\ref{nm}(i)
and Corollary~\ref{lfomr}(ii) as follows:

\begin{Corollary}\label{clearer}
For $g,f \in \Z^{m|n}_+$,
\begin{itemize}
\item[(i)]
$u_{g,f}(q) = 
q^{|\theta|}$ 
if $g = \LL_\theta(f)$ for some $\theta \in \{0,1\}^{\#f}$,
$u_{g,f}(q) = 0$
otherwise;
\item[(ii)]
$l_{g,f}(-q^{-1}) =
\sum_{\theta}
q^{|\theta|}$
summing over all
$\theta \in \N^{\#g}$ with $\RR_\theta'(g) = f$.
\end{itemize}
\end{Corollary}

\begin{Example}\rm
Using Corollary~\ref{clearer}(ii) and arguing by induction on $n$,
one gets that
$l_{g,f}(-q^{-1}) = q^{2}(1+q^{2})^{n-1}$,
for  $f = (0,2,\dots,2n-2|2n-2,\dots,2,0)$ and
$g = (2,4,\dots,2n|2n,\dots,4,2)$.
\end{Example}

\begin{Corollary}\label{dty}
For $f, g \in \Z^{m|n}_+$,
$u_{-g \cdot w_0, -f \cdot w_0}(q) =
q^{\# f} u_{g,\RR(f)}(q^{-1}).$
\end{Corollary}

\begin{proof}
Let $r = \#f$.
By Corollary~\ref{clearer}(i) and Lemma~\ref{opprop}(i), we know that
$u_{-g \cdot w_0, -f \cdot w_0}(q) = 
q^{|\theta|}$
if $g = \RR_{\theta}(f)$
for some $\theta \in \{0,1\}^r$
and is zero otherwise.
Similarly, by Lemma~\ref{opprop}(iii),
$u_{g, \RR(f)}(q)$ is $q^{|\theta|}= q^{r-|\varphi|}$
if $g = \LL_{\theta}(\RR(f))= \RR_{{\varphi}}(f)$
for some $\theta\in\{0,1\}^{r}$ and
$\varphi = (1-\theta_r,\dots,1-\theta_1)$ and is zero otherwise.
\end{proof}

\Point{Length function}\label{lf}
We now consider some further properties of the
polynomials $l_{g,f}(q)$.

\begin{Lemma}\label{moregen}
Let $g, f \in \Z^{m|n}_+$ with $g \preceq f$ and set $r = \#g = \#f$.
There exists a {\em unique} $\theta = \theta(g,f) \in \N^r$
such that
\begin{itemize}
\item[(i)] $f = \RR_\theta'(g)$;
\item[(ii)] if $f = \RR_{\varphi}'(g)$ 
for some $\theta\neq\varphi\in\N^{r}$ 
 then
$|\varphi| < |\theta|$ and $|\varphi| \equiv |\theta| \pmod{2}$.
\end{itemize}
Given in addition $h \in \Z^{m|n}_+$
with $h\preceq g \preceq f$, 
$\theta(h,f) = \theta(h,g)+\theta(g,f)$.
\end{Lemma}

\begin{proof}
We just explain how to construct $\theta$, and leave the rest
of the proof to the reader.
Define
$-m  \leq i_1 < \dots < i_r < 0 < j_r < \dots < j_1 \leq n$
such that $g(i_s) = g(j_s)$ for each $s = 1,\dots,r$,
and
$-m  \leq i_1' < \dots < i_r' < 0 < j_r' < \dots < j_1' \leq n$
such that $f(i_s') = f(j_s')$ for each $s = 1,\dots,r$.
For $0 \leq s \leq r$, let
$$
g_s = g + 
\sum_{t=1}^s (g(i_s) -f(i_s'))(d_{i_s} - d_{j_s}),
$$
so $g_0 = g$ and $(g_r)^+ = f$.
Now for each $s = 1,\dots,r$, 
let $\theta_s$ be the unique non-negative integer such that
$\RR_{i_s,j_s}^{\theta_s}(g_{s-1}) = g_s$, recalling (\ref{rij}),
and take $\theta = \theta(g,f)$ to be the tuple
$(\theta_1,\dots,\theta_r)$.
\end{proof}

Using Lemma~\ref{moregen}, we can introduce a length function
on $\Z^{m|n}_+$.
Suppose to start with that $g \preceq f$.
Let $\ell(g,f) = |\theta(g,f)|$, where $\theta(g,f)$ is the tuple
defined in the lemma.
Notice that if $h \preceq g \preceq f$, then
\begin{equation}\label{lt}
\ell(h,g)+\ell(g,f)=\ell(h,f),
\end{equation} 
as follows from the stronger
fact that $\theta(h,g)+\theta(g,f)=\theta(h,f)$ established
by Lemma~\ref{moregen}.
This allows us to extend the notion of length to arbitrary
$g,f \in \Z^{m|n}_+$
with $\wt(g) = \wt(f)$: pick $h \in \Z^{m|n}_+$
with $h\preceq g$ and $h\preceq f$ and 
set $\ell(g,f) = \ell(h,f) - \ell(h,g)$.
To check that this is well-defined, suppose $h' \in \Z^{m|n}_+$
also satisfies $h' \preceq g$ and $h' \preceq f$.
Choose another $k \in \Z^{m|n}_+$ with
$k \preceq h$ and $k \preceq h'$.
Then using (\ref{lt}),
\begin{align*}
\ell(h,f) - \ell(h,g) &= (\ell(k,h)+\ell(h,f)) - (\ell(k,h)+\ell(h,g))\\
&= \ell(k,f) - \ell(k,g) = \ell(h',f)-\ell(h',g),
\end{align*}
as required. 
So we have now defined $\ell(g,f)$,
the {\em length of $f$ relative to $g$}, for arbitrary $g,f \in \Z^{m|n}_+$
with $\wt(g) = \wt(f)$.
It is immediate from the definition that (\ref{lt}) holds
for all $h,g,f \in \Z^{m|n}_+$ with $\wt(h)=\wt(g)=\wt(f)$.

Finally we can somewhat arbitrarily introduce an absolute notion
of length. For each weight $\gamma$ of $\E^{m|n}$,
we fix a choice of ``origin'' $o_\gamma \in \Z^{m|n}_+$ with
$\wt(o_\gamma) = \gamma$. Then for any $f \in \Z^{m|n}_+$,
we define
\begin{equation}
\ell(f) := \ell(o_\gamma,f)
\end{equation}
where $\gamma = \wt(f)$.
The important thing is that if $\wt(g) = \wt(f)$,
then $\ell(g,f) = \ell(f) - \ell(g)$, so we can recover the
length of $f$ relative to $g$ from the absolute lengths of $f$ and $g$.
In this notation, Lemma~\ref{moregen}(ii) and Corollary~\ref{clearer}(ii)
combine to show:

\begin{Corollary}\label{parcor}
For $g,f \in \Z^{m|n}_+$ with $g \preceq f$,
the polynomial $l_{g,f}(-q^{-1})$ belongs to 
$q^{\ell(f)-\ell(g)} \N[q^{-2}]$,
and the coefficient of $q^{\ell(f) - \ell(g)}$ is $1$.
\end{Corollary}

\section{Representations of $\mathfrak{gl}(m|n)$}

We now relate the combinatorics developed in sections 2 and 3
to two natural categories $\mathcal O_{m|n}$ and $\mathcal F_{m|n}$
of representations of $\mathfrak{gl}(m|n)$.
For basic notions regarding Lie superalgebras, we follow \cite{Kac}.
We denote the parity of a vector $v$ in a vector
superspace by $\bar v \in\Z_2$.
For a Lie superalgebra $\mathfrak g = \mathfrak g_{\0} \oplus
\mathfrak g_{\1}$ and $\mathfrak g$-supermodules $M, N$,  the space
$\hom_{\mathfrak g}(M, N)$
has a canonical $\Z_2$-grading, and the category of all
$\mathfrak g$-supermodules is a superadditive category in the sense of
\cite[ch.3,$\S$2.7]{Man}.
We will use the notation $M \simeq N$ as opposed to the usual
$M \cong N$ to indicate that 
there is an {\em even} isomorphism between $M$ and $N$.
Also $\Pi$ denotes the parity change functor.

\Point{\boldmath Two categories}\label{s9}
From now on, we let $\mathfrak g$ denote the Lie superalgebra
$\mathfrak{gl}(m|n)$. So $\mathfrak g$ consists of
$(m+n)\times(m+n)$ matrices over $\C$, where we label
rows and columns of such matrices by the ordered index set 
$I(m|n) = \{-m,\dots,-1,1,\dots,n\}$ as in the introduction. 
For $i \in I(m|n)$, let $\bar i = \0$ if $i > 0$ and $\1$ if $i < 0$.
Then, the parity
of the $ij$-matrix unit $e_{i,j} \in \mathfrak g$ 
is $\bar i + \bar j$, and the superbracket satisfies
\begin{equation}
[e_{i,j}, e_{k,l}] = \delta_{j, k} e_{i, l} - (-1)^{(\bar i + \bar j)(\bar k + \bar l)} \delta_{i,l} e_{k,j}.
\end{equation}
Note that the subalgebra $\mathfrak g_{\0}$ of $\mathfrak g$
is isomorphic to $\mathfrak{gl}(m) 
\oplus \mathfrak{gl}(n)$.
We will need some other important subalgebras:
let $\mathfrak h$ denote the standard Cartan 
subalgebra of $\mathfrak g$ consisting of all diagonal matrices,
let $\mathfrak b$ be the standard Borel 
subalgebra of all upper triangular matrices,
and let $\mathfrak p = \mathfrak g_\0 + \mathfrak b$.

For $\lambda \in \mathfrak h^*$ and a $\mathfrak g$-supermodule $M$, 
we define the $\lambda$-weight space $M_\lambda$ of $M$ with respect
to $\mathfrak h$ as usual:
$M_\la = \{m \in M\:|\:h m = \la(h) m
\hbox{ for all $h \in \mathfrak h$}\}.$
Given a $\mathfrak g$-supermodule $M$ such that
$M = \bigoplus_{\lambda \in \mathfrak h^*} M_\lambda$,
we can consider the graded dual $M^\star := 
\bigoplus_{\lambda \in \mathfrak h^*}
\hom_{\C}(M_\lambda, \C)$ with the usual $\Z_2$-grading and
$\mathfrak g$-action.
Twisting the $\mathfrak g$-action on $M^\star$ with the
automorphism $X \mapsto - X^{st}$, where 
$st:\mathfrak g \rightarrow \mathfrak g$ is the {\em supertranspose}
$e_{i,j} \mapsto (-1)^{\bar i (\bar i + \bar j)} e_{j,i}$,
we obtain a new $\mathfrak g$-supermodule denoted $M^\tau$.
If all the weight spaces are finite dimensional, then there are natural
isomorphisms $(M^\star)^\star \simeq M$ and $(M^\tau)^\tau \simeq M$.

Let $\{\delta_i\}_{i \in I(m|n)}$ be the basis for $\mathfrak h^*$
dual to the basis $\{e_{i,i}\}_{i \in I(m|n)}$ for $\mathfrak h$.
Define a symmetric bilinear form $(.|.)$ on $\mathfrak h^*$ by declaring
that $(\delta_i | \delta_j) = (-1)^{\bar i} \delta_{i,j}$.
The Weyl group $W$ associated to the reductive Lie algebra
$\mathfrak g_{\0}$ can be identified with the symmetric group $S_{m|n}$
from \ref{s1}. It acts linearly on $\mathfrak h^*$ 
so that $x \delta_i = \delta_{x i}$ for $x \in W, i \in I(m|n)$.
As before, we write $w_0$ for the longest element of $W$.
We will also need the {\em dot action} of $W$ on $X(m|n)$ defined by
$x \cdot \lambda := x (\lambda + \rho) - \rho$
where 
\begin{equation}
\rho = -\sum_{i \in I(m|n)} i \delta_i.
\end{equation}
The {\em root system} of $\mathfrak g$ is the set
$R = \{\delta_i - \delta_j\:|\:i,j \in I(m|n), i \neq j\}.$
We write $R = R_{\0} \cup R_{\1}$ where $R_{\0}$ consists of all
{\em even roots} $\delta_i - \delta_j$ with $\bar i = \bar j$, and
$R_{\1}$ consists of the remaining {\em odd roots}.
Corresponding to the Borel subalgebra $\mathfrak b$,
we have the standard choice of {\em positive roots}
$R^+ = R^+_{\0} \cup R^+_{\1} 
= \{\delta_i - \delta_j\:|\:i,j\in I(m|n), i < j\}$.
The {\em dominance ordering}
on $\mathfrak h^*$ is defined by $\lambda \leq \mu$ if 
$(\mu-\lambda)$ is an $\N$-linear combination of positive roots.

From now on, we will restrict our attention to
the {\em integral weights}, i.e. the weights belonging to 
the subset $X(m|n)$ of $\mathfrak h^*$ consisting of 
all $\Z$-linear combinations of $\{\delta_i\}_{i \in I(m|n)}$.
For $\la =\sum_{i \in I(m|n)} \la_i \delta_i 
\in X(m|n)$, we define its {\em parity}
\begin{equation}\label{pardef}
\bar\la := \la_{-m}+\dots+\la_{-2}+\la_{-1}
\in \Z_2.
\end{equation}
Let $X^+(m|n)$ be the set of all {\em dominant integral weights}, namely,
the $\lambda = \sum_{i \in I(m|n)} 
\lambda_i \delta_i \in X(m|n)$ such that
$\lambda_{-m} \geq \dots \geq \lambda_{-1}$ and 
$\lambda_1 \geq \dots \geq \lambda_n$.
Define a bijection
\begin{equation}\label{shift}
X(m|n) \rightarrow \Z^{m|n}, \quad\lambda \mapsto f_\lambda
\end{equation}
where $f_\lambda \in \Z^{m|n}$ is the function 
defined by $f_\lambda(i) = (\lambda + \rho | \delta_i)$
for $i \in I(m|n)$.
Under this bijection,
$X^+(m|n)$ maps onto $\Z_+^{m|n}$, see \ref{thespace}.
Also $f_{x \cdot \lambda} = f_\lambda \cdot x^{-1}$
for each $x \in W = S_{m|n}$, i.e. the dot action of $W$
on $X(m|n)$ corresponds to the action of $S_{m|n}$ on $\Z^{m|n}$
introduced in \ref{s1}.
Now we lift all the remaining combinatorial definitions 
involving $\Z^{m|n}$ directly to $X(m|n)$. 
For instance, recalling (\ref{dat}), we define the degree of atypicality
$\#\lambda$ of $\lambda \in X(m|n)$ by
$\#\lambda := \# f_\lambda$; this is the same notion 
as in \cite[(1.1)]{Serg}.
Similarly, let $\wt(\la) := \wt(f_\la)$,
an element of the weight lattice $P$, see (\ref{wd}),
and write $\la \preceq \mu$ if $f_\la \preceq f_\mu$, see \ref{bord}.
This ordering on $X(m|n)$
plays the role of the {\em Bruhat ordering}, see e.g. Theorem~\ref{typthm}(ii) below.
It should not be confused with the dominance ordering $\leq$:
we have that $\la \preceq \mu\Rightarrow \la \leq \mu$ but not
conversely.

We are ready to
introduce two categories of representations of $\mathfrak g$.
All the results summarized in the remainder of this subsection
are taken from \cite[section 7]{Btilt}, 
where they are deduced from
a general framework for representations of graded Lie superalgebras
similar to that of Soergel \cite{so2}.

The first category is denoted $\mathcal O_{m|n}$, and is the
(integral weight) analogue of the \cite{BGG} category $\mathcal O$ 
for a semisimple Lie algebra.
By definition, $\mathcal O_{m|n}$ is the category of all
finitely generated $\mathfrak g$-supermodules $M$ that are
locally finite dimensional over $\mathfrak b$ and satisfy
\begin{equation}\label{intnt}
M = \bigoplus_{\lambda \in X(m|n)} M_\lambda.
\end{equation}
An object $P \in \mathcal O_{m|n}$ is {\em projective}
if every (not necessarily even) morphism from $P$ to a 
quotient of an object $M \in \mathcal O_{m|n}$ 
lifts to a morphism from $P$ to $M$.
By \cite[Lemma 7.3]{Btilt}, the category $\mathcal O_{m|n}$
has enough projectives, i.e. every object is a quotient of
a projective object. Moreover, $\mathcal O_{m|n}$
is finite, i.e. every object has a composition series.
For each $\la \in X(m|n)$, we have the {\em Verma module}
\begin{equation}
M(\lambda) := U(\mathfrak g) \otimes_{U(\mathfrak b)} \C_\lambda,
\end{equation}
where $\C_\lambda$ is the one dimensional $\mathfrak b$-module
of weight $\la$
{\em concentrated in degree $\bar \la$}.
The significance of the choice of 
parity here will be explained in \ref{parsec} below. 
As usual, $M(\la)$ has a unique irreducible quotient denoted $L(\la)$, and
$\{L(\lambda)\}_{\lambda \in X(m|n)}$ is a complete set
of pairwise non-isomorphic irreducibles in ${\mathcal O_{m|n}}$.

We say that an object $M \in \mathcal O_{m|n}$ has a 
{\em Verma flag} if it has a filtration
$0 = M_0 < \dots < M_r = M$ such that each $M_i / M_{i-1}$
is $\cong M(\la_i)$ for some $\la_i \in X(m|n)$.
If $M$ has a Verma flag and $\mu \in X(m|n)$, we let
\begin{equation}\label{vm}
(M:M(\mu)) = \dim \hom_{\mathcal O_{m|n}}(M, M(\mu)^\tau).
\end{equation}
By \cite[(6.1)]{Btilt}, this computes the number of subquotients
of a Verma flag of $M$ that are $\cong M(\mu)$.
There is an obvious refinement of these multiplicities: for $p \in \Z_2$,
\begin{equation}\label{vm2}
(M:M(\mu))_p := \dim \hom_{\mathcal O_{m|n}}(M, M(\mu)^\tau)_p
\end{equation}
counts the number of subquotients of a Verma flag of $M$
that are $\simeq \Pi^p M(\mu)$.

By \cite[Theorem 6.3]{Btilt}, there is for each $\la \in X(m|n)$
a unique (up to even isomorphism) indecomposable module
$T(\la) \in \mathcal O_{m|n}$
satisfying the following properties:
\begin{itemize}
\item[(T1)] $T(\lambda)$ has a Verma flag starting with 
$M(\lambda)$ at the bottom;
\item[(T2)] $\ext^1_{\mathcal O_{m|n}}(M(\mu), T(\lambda)) = 0$ for all
$\mu \in X(m|n)$.
\end{itemize}
Moreover, by \cite[(7.4)]{Btilt},
the multiplicity of $M(\mu)$ in a Verma flag of $T(\la)$
is equal to the composition multiplicity
of $L(-\la-2\rho)$ in $M(-\mu-2\rho)$, i.e.
\begin{equation}\label{mdt}
(T(\lambda): M(\mu)) = [M(-\mu - 2\rho):L(-\lambda-2\rho)],
\end{equation}
for $\lambda, \mu \in X(m|n)$.
In particular, 
$(T(\la):M(\la)) = 1$ and 
$(T(\la):M(\mu)) = 0$ unless $\mu \leq \la$.
Consequently, we call $T(\la)$ the 
{\em infinite dimensional tilting module
of highest weight $\la$}.
Finally, note that for every $\la \in X(m|n)$,
\begin{equation}\label{mdt2}
f_{-\lambda-2\rho} = -f_\lambda,
\end{equation}
so the involution $\la \mapsto -\la-2\rho$ on $X(m|n)$ 
appearing in the formula (\ref{mdt})
corresponds
to the involution $f \mapsto -f$ on $\Z^{m|n}$ in $\rho$-shifted notation.

The second category we shall consider is the category
$\mathcal F_{m|n}$ of all finite dimensional
$\mathfrak g$-supermodules satisfying (\ref{intnt}).
Again, this is finite and has enough projectives.
As is well-known, the irreducible finite dimensional 
$\mathfrak g_{\0}$-supermodules with integral highest weights
are parametrized by the set $X^+(m|n)$.
For $\lambda \in X^+(m|n)$, let us write $L'(\lambda)$ 
for the corresponding irreducible highest weight representation of 
$\mathfrak g_{\0}$ {\em concentrated in degree $\bar\la$}.
Then, for each $\la \in X^+(m|n)$, we have the {\em Kac module}
\begin{equation}
K(\lambda) := U(\mathfrak g) \otimes_{U(\mathfrak p)} 
L'(\lambda),
\end{equation}
where we are viewing $L'(\lambda)$ here as a $\mathfrak p$-supermodule
with elements of 
$\mathfrak p_{\0} = \mathfrak g_{\0}$ acting as given and elements of 
$\mathfrak p_{\1}$ acting trivially.
For $\lambda \in X^+(m|n)$, the irreducible module
$L(\lambda)$ defined earlier can also be realized as the unique irreducible
quotient of $K(\lambda)$, and 
$\{L(\lambda)\}_{\lambda \in X^+(m|n)}$ is a complete set of pairwise
non-isomorphic irreducibles in $\mathcal F_{m|n}$.

When working in $\mathcal F_{m|n}$, 
we will talk about {\em Kac flags} in place of Verma flags.
If $M$ has a Kac flag, the number of subquotients of a Kac flag of $M$
that are $\cong K(\mu)$ is denoted $(M:K(\mu))$, and can be computed by
\begin{equation}
(M:K(\mu)) = \dim \hom_{\mathcal F_{m|n}}(M, K(\mu)^\tau).
\end{equation}
Like in (\ref{vm2}), there is a refinement 
denoted $(M:K(\mu))_p$ for $p \in \Z_2$, counting the number of
subquotients of a Kac flag of $M$ that are $\simeq\Pi^p K(\mu)$.
By \cite[(7.6)]{Btilt} or \cite[Proposition 2.5]{Zou},
the projective cover $P(\la)$ of $L(\la)$ in the category $\mathcal F_{m|n}$
has a Kac flag with $K(\la)$ appearing at the top, 
satisfying the BGG reciprocity
\begin{equation}\label{beegeegee}
(P(\la):K(\mu)) = [K(\mu):L(\la)].
\end{equation}
There are also indecomposable tilting modules
in category $\mathcal F_{m|n}$, denoted $U(\la)$ 
for $\lambda \in X^+(m|n)$. Here, by \cite[Theorem 6.3]{Btilt}, 
$U(\la) \in \mathcal F_{m|n}$ is the unique (up to even isomorphism) 
indecomposable object such that
\begin{itemize}
\item[(U1)] $U(\lambda)$ has a Kac flag starting with $K(\lambda)$ at
the bottom;
\item[(U2)] $\ext^1_{\mathcal F_{m|n}}(K(\mu), U(\lambda)) = 0$ for all
$\mu \in X^+(m|n)$.
\end{itemize}
Let $\beta = n(\delta_{-m}+\dots+\delta_{-1}) - m (\delta_1+\dots+\delta_n)$
be the sum of the positive odd roots.
Then, by \cite[(7.7)--(7.8)]{Btilt} and parity considerations,
we have that
\begin{align}
K(\la)^\star &\simeq K(\beta - w_0\la),\label{kdual}\\
U(\la)^\star &\simeq P(\beta - w_0\la).\label{udual}
\end{align}
Note (\ref{beegeegee}), 
(\ref{kdual}) and (\ref{udual}) together imply
\begin{equation}\label{kdt} 
(U(\lambda): K(\mu)) = [K(\beta-w_0 \mu):
L(\beta-w_0 \lambda)],
\end{equation}
for $\lambda, \mu \in X^+(m|n)$.
In particular, $(U(\la):K(\la)) = 1$ and
$(U(\la):K(\mu)) = 0$ unless $\mu \leq \la$.
Accordingly, we will call $U(\la)$ the {\em finite dimensional tilting module
of highest weight $\la$}.
We  remark finally that
\begin{equation}\label{kdt2}
f_{\beta-w_0\la} = -f_\la \cdot w_0 - (m+n+1)\bid,
\end{equation}
where $\bid \in \Z^{m|n}$ is the constant function $i \mapsto 1$.
Thus, up to a constant shift which can usually be ignored, 
the involution $\la \mapsto \beta - w_0\la$ 
on $X^+(m|n)$ appearing in the formula (\ref{kdt}) corresponds
to the involution $f \mapsto -f\cdot w_0$ 
on $\Z_+^{m|n}$.

\Point{Translation functors}\label{stf}
We need some basic facts about central characters.
Let $Z$ be the (even) center of $U(\mathfrak g)$.
The fixed choices of $\mathfrak h \subset \mathfrak b$ 
determine a Harish-Chandra homomorphism
$\varphi:Z \rightarrow U(\mathfrak h),$
see \cite[7.4.3]{Dix}.
Each $\lambda \in \mathfrak h^*$ yields a central character
$\chi_\lambda$ defined by $\chi_\lambda(z) = \lambda(\varphi(z))$.
To parametrize the {\em integral central characters}, i.e. the
$\chi_\lambda$ for $\lambda \in X(m|n)$,
we use the following 
consequence of results of Sergeev \cite{sergeev}, 
see \cite[Corollary 1.9]{Serg}:

\begin{Lemma}\label{hc} Given $\lambda,\mu \in X(m|n)$, we have that
$\chi_\lambda = \chi_\mu$ if and only if $\wt(\lambda) = \wt(\mu)$
(where $\wt(\la) = \wt(f_\la)$, see (\ref{wd})).
\end{Lemma}

For each central character $\chi$, let
$\mathcal O_\chi$ denote the full subcategory of $\mathcal O_{m|n}$
consisting of the modules all of whose composition factors
have central character $\chi$. We have the {\em block decomposition}
\begin{equation*}
\mathcal O_{m|n} = \bigoplus_{\chi} \mathcal O_\chi
\end{equation*}
as $\chi$ runs over all integral central characters.
Lemma~\ref{hc} shows that we can parametrize the
integral characters $\chi$ 
instead by the weights $\gamma \in P$ arising non-trivially in the 
tensor space $\T^{m|n}$ of \ref{s2}.
Let us introduce some notation to do this formally.
Suppose that $\gamma \in P$. Let $\mathcal O_\gamma = \{0\}$ 
if $\gamma$ is not a weight of $\T^{m|n}$, else let $\mathcal 
O_\gamma = 
\mathcal O_{\chi_\lambda}$ where $\lambda \in X(m|n)$ is such that
$\wt(\lambda) = \gamma$.
Then, we can rewrite the above block decomposition as
\begin{equation}
\mathcal O_{m|n} = \bigoplus_{\gamma \in P} \mathcal O_\gamma,
\end{equation}
where $\mathcal O_\gamma$ is non-zero if and only if $\gamma$ is a weight
of $\T^{m|n}$. We let $\pr_\gamma:\mathcal O_{m|n} \rightarrow \mathcal O_\gamma$ be the natural projection functor.
In an entirely similar way, we have the block decomposition of
$\mathcal F_{m|n}$,
\begin{equation}
\mathcal F_{m|n} = \bigoplus_{\gamma \in P} \mathcal F_\gamma,
\end{equation}
where this time $\mathcal F_{m|n}$ is non-zero if and only if 
$\gamma$ is a weight of $\E^{m|n}$, see \ref{thespace}.

Let $V$ be the natural $\mathfrak g$-supermodule.
So, $V$ is the vector superspace on basis $\{v_i\}_{i \in I(m|n)}$,
where $\bar v_i := \bar i$, and the action of the matrix unit 
$e_{i,j} \in \mathfrak g$ is given by 
$e_{i,j} v_k = \delta_{j,k} v_i.$
For $r \geq 0$, let $S^r V$ be the $r$th supersymmetric power of $V$,
a finite dimensional  irreducible representation of $\mathfrak g$.
Let $S^r V^\star = S^r (V^\star) \simeq (S^r V)^\star$.
For $a \in \Z$ and $r \geq 0$, we define additive functors
$F_a^{(r)}, E_a^{(r)} : \mathcal O_{m|n} \rightarrow \mathcal O_{m|n}$
as follows.
It suffices by additivity to define them on objects belonging to
$\mathcal O_\gamma$ for each $\gamma \in P$.
So if $M \in \mathcal O_\gamma$, we let
\begin{align}
F_a^{(r)} M &:= \pr_{\gamma - r (\eps_a - \eps_{a+1})} 
(M \otimes S^r V),\\
E_a^{(r)} M &:= \pr_{\gamma + r (\eps_a - \eps_{a+1})} 
(M \otimes S^r V^\star).
\end{align}
On a morphism $\theta: M \rightarrow N$, 
$F_a^{(r)} \theta$ and $E_a^{(r)} \theta$ 
are defined simply to be the restrictions
of the natural maps $\theta \otimes \operatorname{id}$.
Clearly, the restrictions of $F_a^{(r)}$ and $E_a^{(r)}$ to
$\mathcal F_{m|n}$ give functors
$F_a^{(r)}, E_a^{(r)} : \mathcal F_{m|n} \rightarrow \mathcal F_{m|n}$
too.
The first well-known lemma gives the elementary properties.

\begin{Lemma}\label{tfprops}
On either category $\mathcal O_{m|n}$ or
$\mathcal F_{m|n}$, 
$F_a^{(r)}$ and $E_a^{(r)}$ are exact functors, they
commute with the $\tau$-duality, 
and are both left and right adjoint to each other.
\end{Lemma}

The next lemma is also quite standard, though
we have included a proof since we wish to keep track of parity information.

\begin{Lemma}\label{vf}
Let $\nu_1,\dots,\nu_N$ be the set of weights of $S^r V$
ordered so that $\nu_i > \nu_j \Rightarrow i < j$.
Let $\lambda \in X(m|n)$.
\begin{itemize}
\item[(i)] 
$M(\lambda) \otimes S^r V$
has a multiplicity-free Verma flag with subquotients 
$\simeq
M(\lambda+\nu_1), \dots, M(\lambda+\nu_N)$ in order from
bottom to top.
\item[(ii)]
$M(\lambda) \otimes S^r V^\star$
has a multiplicity-free Verma flag with subquotients 
$\simeq M(\lambda-\nu_1), \dots, M(\lambda-\nu_N)$ in order from
top to bottom.
\end{itemize}
\end{Lemma}

\begin{proof}
We prove (i), (ii) being entirely similar.
By the tensor identity,
$$
M(\la) \otimes S^r V 
=
(U(\mathfrak g) \otimes_{U(\mathfrak b)} \C_\la) \otimes S^r V 
\simeq U(\mathfrak g) \otimes_{U(\mathfrak b)}
(\C_\la \otimes S^r V).
$$
So it suffices by exactness of the functor
$U(\mathfrak g) \otimes_{U(\mathfrak b)} ?$ to show that
$M := \C_\la \otimes S^r V$ has a filtration 
$0 = M_0 < M_1 < \dots < M_N = M$ as a $\mathfrak b$-module
with $M_i / M_{i-1} \simeq \C_{\la+\nu_i}$.
Let $x_1,\dots, x_N$ be a basis for $S^r V$, where
$x_i$ is of weight $\nu_i$.
Then $1 \otimes x_i \in \C_\la \otimes S^r V$ is
of weight $\la+\nu_i$ and degree $\bar\la + \bar x_i
=\bar\la+\bar\nu_i = \overline{\la+\nu_i}$ (recall
(\ref{pardef})).
So taking $M_i$ to be the subspace spanned by 
$1\otimes x_1,\dots,1\otimes x_i$
gives the required filtration.
\end{proof}

\begin{Corollary}\label{cvf}
Let $\lambda \in X(m|n)$ and $a \in \Z$.
Let $(\sigma_{-m},\dots,\sigma_{-1},\sigma_1,\dots,\sigma_n)$
be the $a$-signature of $f_\lambda$,
see (\ref{sigdef}).
\begin{itemize}
\item[(i)] $F_a^{(r)} M(\lambda)$ has a multiplicity-free Verma flag with
subquotients $\simeq M(\lambda + \delta_{i_1} + \dots + \delta_{i_r})$
for all distinct $i_1, \dots, i_r \in I(m|n)$
such that $\sigma_{i_1}= \dots = \sigma_{i_r} = +$.
\item[(ii)]
$E_a^{(r)} M(\lambda)$ has a multiplicity-free Verma flag with
subquotients $\simeq M(\lambda - \delta_{j_1} - \dots - \delta_{j_r})$
for all distinct $j_1, \dots, j_r \in I(m|n)$
such that $\sigma_{j_1} = \dots = \sigma_{j_r} = -$.
\end{itemize}
In both (i) and (ii), the Verma flag can be chosen so that subquotients appear
in order refining dominance, most dominant at the bottom.
\end{Corollary}

\begin{proof}
The Verma module $M(\lambda)$ has central character
$\chi_\lambda$ so belongs to $\mathcal O_{\wt(\lambda)}$
by Lemma~\ref{hc}. 
Applying the exact functor $\pr_{\wt(\la) - 
r (\eps_a - \eps_{a+1})}$ to the filtration in Lemma~\ref{vf}(i),
we deduce that $F_a^{(r)} M(\la)$ has a Verma flag with
subquotients  being the 
$M(\lambda + \nu_i)$ such that $\wt(\lambda + \nu_i) 
= \wt(\lambda) - r (\eps_a - \eps_{a+1})$.
This implies that 
$\nu_i = \delta_{i_1} + \dots + \delta_{i_s}$ for 
distinct $i_1, \dots, i_s \in I(m|n)$
such that $\sigma_{i_1} = \dots = \sigma_{i_s} = +$, giving
(i). Part (ii) is similar.
\end{proof}

There is an analogous statement in the finite dimensional setting.

\begin{Corollary}\label{ckf}
Let $\lambda \in X^+(m|n)$ and $a \in \Z$.
Let $(\sigma_{-m},\dots,\sigma_{-1},\sigma_1,\dots,\sigma_n)$
be the $a$-signature of $f_\lambda$,
see (\ref{sigdef}).
\begin{itemize}
\item[(i)] $F_a^{(r)} K(\lambda)$ has a multiplicity-free Kac flag with
subquotients $\simeq K(\lambda + \delta_{i_1} + \dots + \delta_{i_r})$
for all distinct $i_1, \dots, i_r \in I(m|n)$
such that 
$\lambda + \delta_{i_1} + \dots + \delta_{i_r} \in X^+(m|n)$ and
$\sigma_{i_1} = \dots = \sigma_{i_r} = +$.
\item[(ii)]
$E_a^{(r)} K(\lambda)$ has a multiplicity-free Kac flag with
subquotients $\simeq K(\lambda - \delta_{j_1} - \dots - \delta_{j_r})$
for all distinct $j_1, \dots, j_r \in I(m|n)$
such that 
$\lambda - \delta_{j_1} - \dots - \delta_{j_r} \in X^+(m|n)$ and
$\sigma_{j_1} = \dots = \sigma_{j_r} = -$.
\end{itemize}
In both (i) and (ii), the Kac flag can be chosen so that subquotients appear
in order refining dominance, most dominant at the bottom.
\end{Corollary}

\begin{proof}
We prove (i).
By universal properties,
$K(\la)$ is the largest finite dimensional
quotient of $M(\la)$.
So since $F_a^{(r)}$ is exact,
$F_a^{(r)} K(\la)$ is a quotient of $F_a^{(r)} M(\la)$ and
Corollary~\ref{cvf} implies that $F_a^{(r)} K(\la)$
has a filtration with subquotients being
finite dimensional
quotients of $M(\lambda + \delta_{i_1} + \dots + \delta_{i_r})$
for all distinct $i_1, \dots, i_r \in I(m|n)$
such that $\sigma_{i_1}= \dots = \sigma_{i_r} = +$.
But such a quotient is zero unless $\la+\delta_{i_1}+\dots+\delta_{i_r}
\in X^+(m|n)$. Hence, $F_a^{(r)} K(\la)$ has a filtration with
subquotients being quotients of the Kac modules
$K(\la+\delta_{i_1}+\dots+\delta_{i_r})$ for
all distinct $i_1,\dots,i_r \in I(m|n)$
such that $\la+\delta_{i_1}+\dots+\delta_{i_r} \in X^+(m|n)$
and $\sigma_{i_1} = \dots = \sigma_{i_r} = +$.
Finally the fact that each factor is actually isomorphic
to the corresponding Kac module, rather than a proper quotient,
follows by a character calculation using the Kac character formula
for $K(\mu)$, the Pieri formulae \cite[(5.16),(5.17)]{Mac}
and Lemma~\ref{hc}.
\end{proof}

\begin{Corollary}\label{ttot}
Let $a \in \Z$ and $r \geq 1$.
\begin{itemize}
\item[(i)] For each $\la \in X(m|n)$, 
each indecomposable summand of
$F_a^{(r)} T(\la)$ or of $E_a^{(r)} T(\la)$ 
is $\simeq T(\mu)$ for $\mu \in X(m|n)$.
\item[(ii)] For each $\la \in X^+(m|n)$, 
each indecomposable summand of 
$F_a^{(r)} U(\la)$ or of $E_a^{(r)} U(\la)$ is $\simeq U(\mu)$ for 
$\mu \in X^+(m|n)$.
\end{itemize}
\end{Corollary}

\begin{proof} We prove (i) for $E_a^{(r)}$, the other cases being similar. 
Let $T$ be an indecomposable
summand of $E_a^{(r)} T(\la)$. We need to
show that it has a Verma flag with subquotients $\simeq M(\nu)$
for various $\nu \in X(m|n)$, and that
$\ext^1_{\mathcal O_{m|n}}(M(\mu), T) = 0$ for all $\mu \in X(m|n)$.
The first statement is immediate since
$E_a^{(r)} T(\la)$ has such a Verma flag by 
Corollary~\ref{cvf}, and summands of modules with a Verma flag
also have a Verma flag, see 
\cite[Corollary 4.3]{Btilt}.
For the second statement, Lemma~\ref{tfprops} and a standard argument,
see e.g. \cite[I.4.4]{Jantzen}, shows that
$\ext^1_{\mathcal O_{m|n}}(M(\mu), E_a^{(r)} T(\la))
\simeq \ext^1_{\mathcal O_{m|n}}(F_a^{(r)} M(\mu), T(\la)).$
To see that the right hand side is zero, note that
$F_a^{(r)} M(\mu)$ has a Verma flag by Corollary~\ref{cvf}.
By induction on length using the long exact sequence and
the defining property (T2) of $T(\la)$,
$\ext^1_{\mathcal O_{m|n}}(M, T(\la)) = 0$ 
for every $M \in \mathcal O_{m|n}$ with a Verma flag.
\end{proof}

Let $\mathcal O_{m|n}^\Delta$ be the full subcategory of
$\mathcal O_{m|n}$ consisting of all modules possessing a Verma flag.
Let $K(\mathcal O_{m|n}^\Delta)$ denote the Grothendieck group of
the superadditive 
category $\mathcal O_{m|n}^\Delta$ in the sense of \cite[$\S$2-c]{BK}.
Note $K(\mathcal O_{m|n}^\Delta)$ is a free $\Z$-module on basis
$\{[M(\lambda)]\}_{\lambda \in X(m|n)}$.
Similarly, let $\mathcal F_{m|n}^\Delta$ be the full subcategory
of $\mathcal F_{m|n}$ consisting of all modules possessing
a Kac flag, and let $K(\mathcal F_{m|n}^\Delta)$ denote its Grothendieck group.
So $K(\mathcal F_{m|n}^\Delta)$ is a free $\Z$-module on basis
$\{[K(\lambda)]\}_{\lambda \in X^+(m|n)}$.
In view of Corollary~\ref{cvf} and ~\ref{ckf}, the functors
$F_a^{(r)}$ and $E_a^{(r)}$ map objects in
$\mathcal O_{m|n}^\Delta$ resp. $\mathcal F_{m|n}^\Delta$
to objects in $\mathcal O_{m|n}^\Delta$ resp. $\mathcal F_{m|n}^\Delta$.
Moreover, they preserve short exact sequences in 
$\mathcal O_{m|n}^\Delta$ resp. $\mathcal F_{m|n}^\Delta$.
Hence they induce $\Z$-linear operators on
$K(\mathcal O_{m|n}^\Delta)$
and on $K(\mathcal F_{m|n}^\Delta)$.

Now we make the connection to the 
modules $\T^{m|n}$ and $\E^{m|n}$ from sections 2 and 3 of the article.
Actually we need to specialize these modules at $q = 1$.
So let $\T^{m|n}_{\Z[q,q^{-1}]}$ be the $\Z[q,q^{-1}]$-lattice
in $\T^{m|n}$ spanned by $\{M_f\}_{f \in \Z^{m|n}}$, in the notation
of \ref{s2}.
Let $\E^{m|n}_{\Z[q,q^{-1}]}$ be the $\Z[q,q^{-1}]$-lattice
in $\E^{m|n}$ spanned by $\{K_f\}_{f \in \Z_+^{m|n}}$,
in the notation of \ref{thespace}.
Viewing $\Z$ as a $\Z[q,q^{-1}]$-module so that $q$ acts as $1$, we
define
\begin{align*}
\T^{m|n}_\Z &:= \Z \otimes_{\Z[q,q^{-1}]} \T^{m|n}_{\Z[q,q^{-1}]},\\
\E^{m|n}_\Z &:= \Z \otimes_{\Z[q,q^{-1}]} \E^{m|n}_{\Z[q,q^{-1}]}.
\end{align*}
We write $M_f(1)$ (resp. $K_f(1)$) for the basis element $1 \otimes M_f$
of $\T^{m|n}_{\Z}$ (resp. $1 \otimes K_f$
of $\E^{m|n}_{\Z}$).
Similarly, we define $T_f(1) = 1 \otimes T_f$ and $U_f(1) = 1 \otimes U_f$
(in case of $T_f(1)$, recall that as
a consequence of Conjecture~\ref{pc}
we expect it is a finite sum of $K_g(1)$'s so belongs to 
$\T^{m|n}_{\Z}$, but without this we mean
here to 
work in the completion $\widehat\T^{m|n}_\Z$ constructed as in \ref{s2}).

Note the generators $E_a^{(r)}$ and $F_a^{(r)}$ of
$\U = U_q(\mathfrak{gl}_\infty)$ specialize at $q = 1$
to the usual divided powers $E_a^r / r!$ and $F_a^r / r!$ in the 
Chevalley generators of the Lie algebra $\mathfrak{gl}_\infty$,
so we can view $\T^{m|n}_\Z$ resp. $\E^{m|n}_\Z$ as 
modules over the Kostant $\Z$-form $\U_\Z$ for the universal enveloping algebra
$U(\mathfrak{gl}_\infty)$.

\begin{Theorem}\label{gid}
Identify
$K(\mathcal O_{m|n}^\Delta)$ with $\T^{m|n}_{\Z}$ via the
$\Z$-module isomorphism
$$
i: K(\mathcal O_{m|n}^\Delta) \rightarrow \T^{m|n}_{\Z},
\quad [M(\lambda)] \mapsto M_{f_\lambda}(1).
$$
Then, the representation theoretically defined
operators $F_a^{(r)}, E_a^{(r)}$
act in the same way as 
the Chevalley generators $F_a^{(r)}, E_a^{(r)}$ of $\U_\Z$.
\end{Theorem}

\begin{proof}
Corollary~\ref{cvf} shows that the operators induced by the
functors $F_a^{(r)}, E_a^{(r)}$ act on $[M(\lambda)] \in K(\mathcal O_{m|n}^\Delta)$
in exactly the same way
on $K(\mathcal O_{m|n}^\Delta)$
as $F_a^{(r)}, E_a^{(r)} \in \U_\Z$ act on 
$M_{f_\lambda}(1)\in\T^{m|n}_{\Z}$.
\end{proof}

An entirely similar argument, using Corollary~\ref{ckf} instead, gives
the analogous theorem for category $\mathcal F_{m|n}^\De$:

\begin{Theorem}\label{fdede}
Identify
$K(\mathcal F_{m|n}^\Delta)$ with $\E^{m|n}_{\Z}$ via the
$\Z$-module isomorphism
$$
j: K(\mathcal F_{m|n}^\Delta) \rightarrow \E^{m|n}_{\Z},
\quad [K(\lambda)] \mapsto K_{f_\lambda}(1).
$$
Then, the representation theoretically defined
operators $F_a^{(r)}, E_a^{(r)}$
act in the same way as 
the Chevalley generators $F_a^{(r)}, E_a^{(r)}$ of $\U_\Z$.
\end{Theorem}

\Point{\boldmath Tilting modules in category $\mathcal O_{m|n}$}\label{s10}
We proceed to prove some results and formulate some conjectures
about the infinite dimensional tilting modules $T(\lambda)$.
For $\lambda \in X(m|n)$, write
$$
M'(\lambda) := U(\mathfrak g_{\0}) \otimes_{\mathfrak b_{\0}} \C_\lambda
$$
for the purely even Verma module for $\mathfrak g_{\0}$
concentrated in degree $\bar\la$,
and $L'(\lambda)$ for its unique irreducible quotient.
We will need the following result of Kac \cite[Proposition 2.9]{Kacnote}.
Actually in {\em loc. cit.}, Kac is only concerned with finite dimensional
representations, but the same argument works for the general
case stated here.

\begin{Lemma}\label{kac}
If $\lambda \in X(m|n)$ is typical, then
$L(\lambda) \simeq U(\mathfrak g) \otimes_{U(\mathfrak p)}
L'(\lambda)$.
\end{Lemma}

Recall the definition of the polynomials $t_{g,f}(q)$ and
$l_{g,f}(q)$ from (\ref{tde}).
We use the bijection (\ref{shift}) to shift notation, letting
$t_{\mu,\la}(q) := t_{f_\mu,f_\la}(q)$ and $l_{\mu,\la}(q) := 
l_{f_\mu,f_\la}(q)$.
The first part of the following theorem
is a reformulation of the
Kazhdan-Lusztig conjecture \cite{KL} 
for 
$\mathfrak{gl}(m)\oplus \mathfrak{gl}(n)$,
proved in \cite{BB, BrK}.

\begin{Theorem}\label{typthm} Let $\lambda\in X(m|n)$.
\begin{itemize}
\item[(i)] If $\la$ is typical then  $(T(\lambda): M(\mu)) = 
t_{\mu, \lambda}(1)$ for each $\mu\in X(m|n)$.
\item[(ii)] For arbitrary $\la$, each subquotient of a Verma flag
of $T(\la)$ is $\simeq M(\mu)$ for 
$\mu\preceq\la$.
\end{itemize}
\end{Theorem}

\begin{proof}
(i) 
For the proof, we will assume instead 
that $\lambda \in X(m|n)$ is typical 
with $\lambda + \rho \in X^+(m|n)$.
Let $W_\lambda$
be the stabilizer in $W\cong S_{m|n}$ 
of $\lambda$ under the dot action, and 
let $D^\lambda$ be the set of all maximal length $W / W_\lambda$-coset
representatives.
Let $w_\lambda$ be the longest element of $W_\lambda$.
By the Kazhdan-Lusztig conjecture for the Lie algebra
$\mathfrak g_{\0}$
proved in \cite{BB, BrK} combined with the translation principle
\cite{Jan}, see also \cite[Theorem 3.11.4]{BGS},
we have that
$$
[M'(x \cdot \la): L'(y \cdot \lambda)]
=
P_{x,y}(1)
$$
for arbitrary $x,y \in D^\lambda$.
Here, $P_{x, y}(1)$ denotes the usual Kazhdan-Lusztig polynomial
associated to $x, y \in W$ evaluated at $1$, see \cite{KL}.

We claim that 
$P_{x,y}(1) = t_{-x\cdot \la - 2\rho, -y\cdot \la - 2\rho}(1)$
for all $x,y \in D^{\lambda}$.
To see this, let $f := - f_\lambda$, which is antidominant in the sense
of \ref{s1}.
Define $S_f, D_f$ as in 
\ref{s5}. We will use the fact that the map
$D^\lambda \rightarrow D_f, x \mapsto w_\lambda x^{-1}$ is a bijection.
Observe using (\ref{mdt2}) that  
$f_{-y \cdot \lambda - 2 \rho}
= f \cdot w_\lambda y^{-1}$.
So by Lemma~\ref{tc}, $t_{-x\cdot\la-2\rho,-y\cdot\la-2\rho}(1) = 
t_{f\cdot w_\lambda x^{-1},f\cdot w_\lambda y^{-1}}(1)
=m_{w_\lambda x^{-1}, w_\lambda y^{-1}}^{(f)}(1)$.
Noting $m_{w_\lambda x^{-1}, w_\lambda y^{-1}}^{(f)}(1)$
is the same as the element with the same name
in \cite{soergel},
\cite[Remark 2.6]{soergel} and \cite[Proposition 3.4]{soergel}
show that
$m_{w_\lambda x^{-1}, w_\lambda y^{-1}}^{(f)}(1) = P_{x^{-1}, y^{-1}}(1) = 
P_{x, y}(1)$.
This proves the claim.

Now $M(x \cdot \la) \simeq U(\mathfrak g) \otimes_{U(\mathfrak p)} 
M'(x\cdot\la)$ by associativity of tensor product,
while Lemma~\ref{kac} shows that 
$L(y\cdot \la) \simeq U(\mathfrak g) \otimes_{U(\mathfrak p)} 
L'(y \cdot\la)$.
So as the functor $U(\mathfrak g) \otimes_{U(\mathfrak p)} ?$ is exact,
we deduce from the previous two paragraphs that
$$
[M(x \cdot \la):L(y \cdot \lambda)]
=
t_{-x\cdot\la-2\rho,
-y\cdot\la-2\rho}(1).
$$
Note moreover that 
this argument shows that every subquotient
of $M(x \cdot \la)$ that is $\cong L(y\cdot \lambda)$
is actually $\simeq L(y \cdot \la)$.
Finally applying (\ref{kdt}) gives that
$$
(T(-y\cdot\la-2\rho):M(-x\cdot\la-2\rho))
=
t_{-x\cdot\la-2\rho,
-y\cdot\la-2\rho}(1).
$$
Part (i) of the theorem follows easily from this and central character
considerations.
Moreover, by an obvious refinement of (\ref{kdt}) keeping 
track of parity information too,
we see that every subquotient of a Verma flag
of $T(-y\cdot\la-2\rho)$ that is $\cong M(-x\cdot\la-2\rho)$
is actually $\simeq M(-x\cdot\la-2\rho)$.

(ii)
We proceed by induction on $\#\lambda$. The case that $\lambda$ is typical
follows from (i). So suppose that
$\#\lambda > 0$ and the theorem has been proved for all
$\mu$ with $\#\mu < \#\lambda$.
Let $i: K(\mathcal O_{m|n}^\Delta) \rightarrow \T^{m|n}_{\Z}$
be the map defined in Theorem~\ref{gid}.
Apply the algorithm explained in \ref{s5}
to $f = f_\lambda$ to construct $h = f_\nu$ 
for $\nu \in X(m|n)$ with $\# \nu < \#\lambda$ and 
a sequence $X_1, \dots, X_N$ of monomials in $E_a^{(r)}$ and
$F_a^{(r)}$. 
Let $M := X_N \dots X_1 T(\nu)$.
Note by Corollary~\ref{cvf} that $M$ has a Verma flag, and each
subquotient of a Verma flag of $M$ is $\simeq M(\mu)$ for
some $\mu \in X(m|n)$.
By the induction hypothesis, 
$i([T(\nu)])$ equals $M_h(1)$ plus a linear combination of
$M_g(1)$'s with $g \prec h$. 
By Lemmas~\ref{enasty} and \ref{fnasty}, we deduce that
$X_N \dots X_1 i([T(\nu)])$ equals $M_f(1)$ plus a linear combination
of $M_g(1)$'s with $g \prec f$.
So by Theorem~\ref{gid},
$$
[M] = [M(\lambda)] + (\hbox{a linear combination of
$[M(\mu)]$'s with $\mu \prec \lambda$}).
$$
By Corollary~\ref{ttot}(i),
$T(\la)$ is a summand of $M$, and the result follows.
\end{proof}

Motivated by the theorem, we formulate the following conjecture.

\begin{Conjecture}\label{mc}
Let $i: K(\mathcal O_{m|n}^\Delta) \rightarrow \T^{m|n}_{\Z}$
be the map defined in Theorem~\ref{gid}.
Then, 
$i([T(\lambda)]) = T_{f_\lambda}(1)$
for each $\lambda \in X(m|n)$.
\end{Conjecture}

In view of (\ref{mdt}) this conjecture is equivalent to either of
the statements
\begin{align}
(T(\lambda) : M(\mu)) &= t_{\mu, \lambda}(1),\\
[M(\la):L(\mu)] &= t_{-\la-2\rho,-\mu-2\rho}(1)\label{mmult}
\end{align}
for all $\lambda,\mu \in X(m|n)$.
By Corollary~\ref{lcornew} and (\ref{mdt2}),
the unitriangular matrices 
$(l_{\mu,\lambda}(1))_{\mu,\lambda \in X^(m|n)}$
and $(t_{-\lambda-2\rho, -\mu-2\rho}(1))_{\mu,\lambda\in X^+(m|n)}$
are inverse to each other,
so inverting (\ref{mmult}) also gives that
\begin{equation}
\ch L(\lambda) = {\sum_{\lambda \in X(m|n)}}
l_{\mu, \lambda}(1) \ch M(\mu).
\end{equation}
Although the summation is infinite here, it involves only finitely
many non-zero contributions to the dimensions of each 
fixed weight space of $L(\la)$, thus it 
can be viewed as a conjectural
character formula for irreducibles in $\mathcal O_{m|n}$.

Further evidence for Conjecture~\ref{mc} is given by
the main theorem in the next subsection. We finally mention one other
result which is in keeping with the conjecture, compare in particular
with Theorem~\ref{cthm}. 
Recall the definition of the dual crystal operators
$\tilde E_a^*, \tilde F_a^*, \eps_a^*$ and $\phi_a^*$ 
from \ref{s6}. Again, we lift these directly to $X(m|n)$
via the bijection (\ref{shift}).

\begin{Theorem} [Kujawa]\label{kuj}
Let $\lambda \in X(m|n)$ and $a \in \Z$.
\begin{itemize}
\item[(i)] 
$F_a L(\lambda) \neq 0$ if and only if 
$\phi_a^*(\lambda) \neq 0$, in which case it is a $\tau$-self-dual indecomposable
module with irreducible socle and cosocle $\simeq L(\tilde F_a^*(\la))$.
Moreover, $F_a L(\lambda)$ is irreducible if and only if
$\phi_a^*(\lambda) = 1$.
\item[(ii)] 
$E_a L(\lambda) \neq 0$ if and only if 
$\eps_a^*(\lambda) \neq 0$, in which case it is a $\tau$-self-dual 
indecomposable
module with irreducible socle and cosocle $\simeq L(\tilde E_a^*(\la))$. 
Moreover, $E_a L(\lambda)$ is irreducible if and only if
$\eps_a^*(\lambda) = 1$.
\end{itemize}
\end{Theorem}

\vspace{1mm}

Theorem~\ref{kuj} is a result of Jon Kujawa that will form part of
his PhD thesis \cite{kthesis}. 
The proof, which will hopefully appear elsewhere, is similar to the proof
given in \cite{JWB:branching}
of Kleshchev's modular branching rules from \cite{KbrII}. 
It involves some explicit calculations 
with certain lowering operators in $U(\mathfrak g)$.

\Point{\boldmath Tilting modules in category $\mathcal F_{m|n}$}\label{s11}
Now we study the finite dimensional tilting modules $U(\la)$.
Lift the crystal operators
$\tilde E_a, \tilde F_a, \eps_a, \varphi_a$ from \ref{xs} 
to $X^+(m|n)$ through the bijection (\ref{shift}), as well as the
mutually inverse
bijections $\LL$ and $\RR$ from (\ref{LLdef}).

\begin{Theorem}\label{MT}
Let $j:K(\mathcal F_{m|n}^\Delta) \rightarrow
\E^{m|n}_{\Z}$ be the map defined in Theorem~\ref{fdede}.
Then, $j([U(\la)]) = U_{f_\la}(1)$ for each $\la \in X^+(m|n)$.
Moreover:
\begin{itemize}
\item[(i)] Each subquotient of a Kac flag of $U(\la)$ is $\simeq K(\mu)$
for $\LL(\la) \preceq \mu \preceq \la$;
\item[(ii)] $U(\la)\simeq P(\LL(\la))$;
\item[(iii)] $U(\la) \simeq U(\la)^\tau$.
\end{itemize}
\end{Theorem}

\begin{proof}
If $\#\la = 0$ then Lemma~\ref{kac}
implies that $U(\la) = P(\la) = K(\la) = L(\la)$ 
and the theorem follows in this case.
Now suppose that $\#\la > 0$.
Let $f = f_\la$ and define $h = f_\nu$ for $\nu \in X^+(m|n)$
and operators
$X_a \in \{E_a,F_a\}_{a \in \Z}$
according to Procedure~\ref{alg2}.
We may assume by induction that the theorem has been proved for $\nu$.

Consider $X_a U(\nu)$. 
Theorem~\ref{fdede}, Lemma~\ref{props} and the induction hypothesis shows that
$j([X_a U(\nu)]) = X_a U_{f_\nu}(1) = U_{f_\la}(1)$.
So we get from the explicit description of $U_{f_\la}(1)$ in
Theorem~\ref{nm}(i) that $[X_a U(\nu)] = [K(\la)] + (*) + [K(\LL(\la))]$
where $(*)$ is a sum of $[K(\mu)]$'s for $\LL(\la)\prec \mu \prec \la$.
Using Corollary~\ref{ckf}, we deduce from this and the induction
hypothesis that $X_a U(\nu)$
has a Kac flag with subquotients 
$\simeq K(\la)$, $K(\LL(\la))$ and 
all other subquotients
$\simeq K(\mu)$ for $\LL(\la)\prec\mu\prec \la$.
So $X_a U(\nu)$ must have a summand that is $\simeq U(\la)$,
recalling Lemma~\ref{ttot}(ii).
Also $U(\nu) \simeq U(\nu)^\tau$ is projective by the induction hypothesis,
hence $X_a U(\nu) \simeq (X_a U(\nu))^\tau$ is projective
by Lemma~\ref{tfprops}.
So $X_a U(\nu)$ must have a summand that is $\simeq P(\LL(\la))$.
To complete the proof, it just remains to show that
$X_a U(\nu)$ is indecomposable.
For this, we give two different arguments, the first based on
Theorem~\ref{kuj} and the second using instead a fundamental fact proved by
Serganova in \cite{Serg}.

\vspace{1mm}
\noindent{\em Method one.}
Suppose the space
$$
\hom_{\mathcal F_{m|n}}(X_a U(\nu), L(\mu))
\simeq
\hom_{\mathcal F_{m|n}}(U(\nu), Y_a L(\mu))
$$
is non-zero for some $\mu \in X^+(m|n)$.
By the choice of $a$ in Procedure~\ref{alg2}, $\la$ is
{not} at the end of an $a$-string of length $2$ in the crystal graph.
Since we must have that $\wt(\mu) = \wt(\la)$
by Lemma~\ref{hc}, it follows that $\mu$ is
also {not} at the end of an $a$-string of length $2$.
Theorem~\ref{kuj}
now implies immediately that $Y_a L(\mu) \simeq L(\tilde Y_a^*(\mu))$.
By the induction hypothesis,
$U(\nu)$ is the projective cover of $L(\LL(\nu))$, so we deduce from
the non-vanishing of the right hand hom space above that
$\tilde Y_a^*(\mu) = \LL(\nu)$.
Hence, $\mu = \tilde X_a^*(\LL(\nu)) = \LL(\tilde X_a (\nu)) = \LL(\la)$,
using Remark~\ref{lrem}(2) for the penultimate equality.
We have now shown that $\cosoc_{\mathfrak g} (X_a U(\nu)) \simeq
L(\LL(\la))$, so it is indecomposable.

\vspace{1mm}
\noindent{\em Method two.}
Suppose $X_a U(\nu)$ is decomposable.
Then, by what we have shown already, we can write
$X_a U(\nu) = T_1 \oplus T_2$
where $T_1 \simeq U(\la)$ and $T_2 \neq 0$ 
is a direct sum of indecomposable tilting modules.
Note that $Y_a T_i \neq 0$ for each $i$, indeed we have by adjointness that
$$
\hom_{\mathcal F_{m|n}}(U(\nu), Y_a T_i) \simeq
\hom_{\mathcal F_{m|n}}
(X_a U(\nu), T_i) \neq 0.
$$
Recalling Lemma~\ref{props}, we now consider two cases.
First, suppose that $\nu$ is at the end of an $a$-string of length $1$.
Then, we have that $Y_a U_{f_\la}(1) = U_{f_\nu}(1)$, i.e.
$[Y_a X_a U(\nu)] = [U(\nu)]$.
Since $Y_a X_a U(\nu)$ is a direct sum of 
indecomposable tilting modules, we deduce that 
$U(\nu) \cong Y_a X_a U(\nu) \cong 
Y_a T_1 \oplus Y_a T_2$, a contradiction since $U(\nu)$ 
is indecomposable.
Otherwise, we have that $\nu$ is at the end of an $a$-string of length $2$,
and $[Y_a X_a U(\nu)] = 2[U(\nu)]$.
Hence this time we must have that $Y_a T_1 \cong Y_a T_2 \cong U(\nu)$.
In particular, we get that $[Y_a U(\la):L(\nu)] = 1$.
We now show that $[Y_a U(\la):L(\nu)] \geq 2$, to get the 
desired contradiction.

Let $\mu = \tilde X_a^*(\nu)$, so $\mu = \la - \alpha$
for some $\alpha \in R_{\1}^+$ with $(\la+\rho|\alpha) = 0$.
By  \cite[Theorem 5.5]{Serg} and (\ref{kdt}), we have that
$[K(\la):L(\mu)] \geq 1$ and that
$(U(\la):K(\mu)) = 
[K(\beta-w_0\mu):L(\beta-w_0\la)] \geq 1$.
Hence, $[U(\la):L(\mu)] \geq 2$, since it has a Kac flag involving
both $K(\la)$ and $K(\mu)$, 
each of which have $L(\mu)$ as a composition factor.
Now $X_a K(\nu)$ has a two-step filtration with 
$K(\la)$ at the bottom
and $K(\mu)$ at the top, hence
$$
\hom_{\mathcal F_{m|n}}(K(\nu), Y_a L(\mu))
\simeq
\hom_{\mathcal F_{m|n}}(X_a K(\nu), L(\mu)) \neq 0.
$$
This shows that $[Y_a L(\mu):L(\nu)] \geq 1$.
Finally applying the exact functor $Y_a$ to $U(\la)$ 
and combining our two facts
$[U(\la):L(\mu)] \geq 2$ and $[Y_a L(\mu):L(\nu)] \geq 1$
gives that $[Y_a U(\la):L(\nu)] \geq 2$ as required.
\end{proof}

Now recall the definition of the polynomials
$u_{g,f}(q)$ and $l_{g,f}(q)$ from (\ref{tde2}).
As usual we shift notation, writing
$u_{\mu,\la}(q) := u_{f_\mu,f_\la}(q)$ and $l_{\mu,\la}(q) := 
l_{f_\mu,f_\la}(q)$.
Combining the theorem with
(\ref{kdt}), we get that
\begin{align}\label{cdcd}
(U(\lambda) : K(\mu)) &= u_{\mu, \lambda}(1),\\
[K(\la):L(\mu)] &= u_{\beta-w_0\la,\beta-w_0\mu}(1).\label{lie}
\end{align}
The Main Theorem stated in the introduction follows immediately from the
second of these formulae
and Corollary~\ref{lfomr}(i), since in view of (\ref{kdt2}) 
and Corollary~\ref{lcornew2}
that gives an explicit
formula for $u_{\beta-w_0\la,\beta-w_0\mu}(1)$.
In particular, $[K(\lambda):L(\mu)] \leq 1$ for all $\lambda,\mu\in X^+(m|n)$,
as was conjectured in \cite[Conjecture 7.2]{HKJ}, and
$L(\mu)$ appears as a composition factor 
in exactly $2^{\#\mu}$ different Kac modules $K(\lambda)$,
as was conjectured in \cite[Corollary 7.3]{HKJ}.

By Corollary~\ref{lcornew2} and (\ref{kdt2}),
the unitriangular matrix
$(l_{\mu,\la}(1))_{\mu,\la \in X^+(m|n)}$ is the inverse
of $(u_{\beta - w_0\la, \beta-w_0\mu}(1))_{\mu,\la \in X^+(m|n)}$.
So on inverting (\ref{lie}), we also get that
\begin{equation}\label{chfo}
\ch L(\lambda) = \sum_{\mu \in X^+(m|n)}
l_{\mu, \lambda}(1) \ch K(\mu).
\end{equation}
This can be viewed as a character formula for
the finite dimensional irreducible $\mathfrak{gl}(m|n)$-supermodules
with integral highest weight.
The explicit description of the coefficients
$l_{\mu, \lambda}(1)$ given by Corollary~\ref{clearer}(ii)
seems to be quite different from
the explicit description given by Serganova \cite[Theorem 2.3]{Serg}, and
I have been unable
to prove combinatorially that they are equivalent.

To conclude the subsection, let us record one more consequence of
Theorem~\ref{MT}.

\begin{Corollary}\label{ndty}
For $\la \in X^+(m|n)$, $L(\la)^\star \simeq L(\beta - w_0 \RR(\la))$.
\end{Corollary}

\begin{proof}
By (\ref{udual}) and Theorem~\ref{MT}, 
$P(\beta - w_0 \RR(\la)) \simeq
U(\RR(\la))^\star
\simeq P(\la)^\star$ and it is self-dual under the duality $\tau$.
Hence $L(\beta - w_0 \RR(\la))
\simeq \soc_{\mathfrak g} P(\beta - w_0 \RR(\la))
\simeq \soc_{\mathfrak g} P(\la)^\star \simeq
(\cosoc_{\mathfrak g} P(\la))^\star
\simeq L(\la)^\star$.
\end{proof}

\begin{Remark}\rm
A different description of the highest weight 
of $L(\la)^\star$ can be derived using
Serganova's odd reflections, see \cite{sergthesis},
\cite[Lemma 0.3]{PSI} and \cite[Theorem 4.5]{BKu}.
In view of (\ref{lie}) and (\ref{kdual}), Corollary~\ref{ndty}
implies (indeed is equivalent to) the equality
$u_{\beta - w_0\mu,\beta - w_0\la}(1) = u_{\mu, \RR(\la)}(1)$,
see Corollary~\ref{dty} for a stronger statement.
\end{Remark}

\Point{Highest weight categories}\label{parsec}
Let $\mathcal F_{m|n}^{\0}$ (resp. $\mathcal F_{m|n}^{\1}$)
be the full subcategory
of $\mathcal F_{m|n}$ consisting of the modules all of whose composition
factors are $\simeq L(\la)$ (resp. $\simeq \Pi L(\la)$)
for $\la \in X^+(m|n)$.
Obviously, the parity change functor $\Pi$ defines 
an isomorphism between $\mathcal F_{m|n}^{\0}$ and $\mathcal F_{m|n}^{\1}$.
Since each
$\End_{\mathcal F_{m|n}}(L(\la))$ is concentrated in degree $\0$,
each $\hom_{\mathcal F_{m|n}}(M,N)$ for $M,N \in \mathcal F_{m|n}^{\0}$
is also concentrated in degree $\0$, hence
$\mathcal F_{m|n}^{\0}$ is an abelian category.

\begin{Lemma}\label{mas}
For $\la \in X^+(m|n)$, 
each of the objects $U(\la), P(\la), K(\la)$ and $L(\la)$
belong to $\mathcal F_{m|n}^{\0}$.
Moreover, the dualities $\star$ and $\tau$ and the
functors $F_a^{(r)}$ and $E_a^{(r)}$ map objects in $\mathcal F_{m|n}^{\0}$
to objects in $\mathcal F_{m|n}^{\0}$.
\end{Lemma}

\begin{proof}
By Theorem~\ref{MT}(i), each subquotient of a Kac flag of 
$U(\la)$ is $\simeq K(\mu)$ for some $\mu \in X^+(m|n)$.
We deduce using (\ref{kdual})--(\ref{udual}) that
each subquotient of a Kac flag of $P(\la)$ is $\simeq K(\mu)$
for some $\mu \in X^+(m|n)$.
By the obvious refinement of (\ref{beegeegee}) 
keeping track of parities, it follows that
each composition factor of $K(\mu)$ is $\simeq L(\la)$
for some $\la \in X^+(m|n)$.
Combining these statements shows that all of 
$U(\la), P(\la), K(\la)$ and $L(\la)$ belong to $\mathcal F_{m|n}^{\0}$.
For the remaining statement, we obviously have that $L(\la)^\tau \simeq L(\la)$, hence $\tau$ leaves $\mathcal F_{m|n}^{\0}$ invariant. The
same thing for $\star$ follows from Corollary~\ref{ndty}.
Finally, 
Corollary~\ref{ckf} shows that the exact functors 
$F_a^{(r)}$ and $E_a^{(r)}$
send $K(\la)$ to an object in $\mathcal F_{m|n}^{\0}$, and $L(\la)$
is a quotient of $K(\la)$ so they must also send $L(\la)$
to an object in $\mathcal F_{m|n}^{\0}$.
\end{proof}

\begin{Corollary}\label{indegzero}
For any $M, N \in \mathcal F_{m|n}^{\0}$ and $i \geq 0$,
the space 
$\ext^i_{\mathcal F_{m|n}}(M, N)$
is concentrated in degree $\0$.
\end{Corollary}

\begin{proof}
We have already noted this is the case if $i = 0$.
To get the general case from this, note by the lemma 
that every composition factor of every term of
the obvious minimal projective resolution of $M$
belongs to $\mathcal F_{m|n}^{\0}$.
\end{proof}

It follows easily from the corollary that 
every object $M \in \mathcal F_{m|n}$ decomposes uniquely as
$M = M^{\0} \oplus M^{\1}$ with $M^p \in \mathcal F_{m|n}^p$
for each $p \in \Z_2$.
We deduce that there is a decomposition
 $\mathcal F_{m|n} = \mathcal F_{m|n}^{\0} \oplus \Pi
\mathcal F_{m|n}^{\0}$ allowing us to 
reconstruct the superadditive category
$\mathcal F_{m|n}$ from the additive category
$\mathcal F_{m|n}^{\0}$.
For example,  for $M, N \in \mathcal F_{m|n}$ and $i \geq 0$, we have that
\begin{align}\label{evones}
\ext^i_{\mathcal F_{m|n}}(M,N)_{\0} &= \ext^i_{\mathcal F_{m|n}^{\0}}(M^{\0}, N^{\0})
\oplus \ext^i_{\mathcal F_{m|n}^{\0}}(\Pi M^{\1}, \Pi N^{\1}),\\
\ext^i_{\mathcal F_{m|n}}(M,N)_{\1} &= \ext^i_{\mathcal F^{\0}_{m|n}}(M^{\0}, \Pi N^{\1})
\oplus \ext^i_{\mathcal F_{m|n}^{\0}}(\Pi M^{\1}, N^{\0}).\label{oddones}
\end{align}
At this point, we refer the reader to \cite{CPShw, CPSd}
for the definition of a {\em highest weight category
with duality}.

\begin{Theorem}\label{hwc}
The category $\mathcal F_{m|n}^{\0}$ is a highest weight category
with weight poset $(X^+(m|n), \preceq)$ and duality $\tau$.
For $\la \in X^+(m|n)$, 
$U(\la), P(\la), K(\la)$ and $L(\la)$
are the indecomposable 
tilting, projective, standard and irreducible
modules parametrized by $\la$, respectively.
\end{Theorem}

\begin{proof}
We have seen in Theorem~\ref{MT} that
$(P(\la):K(\mu)) \neq 0 \Rightarrow \la \preceq \mu$.
Given this and (\ref{beegeegee}) 
it is a routine matter to check that $\mathcal F_{m|n}^{\0}$
satisfies the axioms for a highest weight category with duality.
\end{proof}

\begin{Remark}\rm
In an entirely similar fashion, we define
$\mathcal O_{m|n}^{\0}$
to be the full subcategory of $\mathcal O_{m|n}$
consisting of the objects $M$ all of whose composition factors
are $\simeq L(\la)$ for $\la \in X(m|n)$.
Using Theorem~\ref{typthm}(ii) and the refined versions of
 (\ref{mdt}) and BGG reciprocity
\cite[(6.6)]{Btilt} keeping track of parity, one can prove 
analogues of all the results in this subsection for $\mathcal O_{m|n}^{\0}$:
there is a decomposition
$\mathcal O_{m|n} = \mathcal O_{m|n}^{\0} \oplus \Pi\mathcal O_{m|n}^{\0}$,
and $\mathcal O_{m|n}^{\0}$
is a highest weight category with weight poset $(X(m|n), \preceq)$
and duality $\tau$.
\end{Remark}

\Point{Kazhdan-Lusztig polynomials}\label{klpol}
In this subsection, we explain the true significance of the polynomials
$l_{\mu,\la}(q) = l_{f_\mu,f_\la}(q)$ for $\mu, \la \in X^+(m|n)$.

\begin{Lemma}\label{res}
Let $\mu \in X^+(m|n)$.
Then, $K(\mu)$ has a projective resolution
$
\dots \rightarrow P_1(\mu) \rightarrow P_0(\mu) \rightarrow K(\mu) \rightarrow 0$
in $\mathcal F_{m|n}^{\0}$
such that for every $\la \in X^+(m|n)$,
$$
\sum_{i \geq 0}
\dim \hom_{\mathcal F_{m|n}} (P_i(\mu), L(\la)) q^i
=l_{\mu,\la}(-q^{-1}).
$$
\end{Lemma}

\begin{proof}
We first explain how to construct for fixed $d \geq 0$
an exact sequence
$P_d(\mu) \rightarrow\dots\rightarrow P_0(\mu) \rightarrow K(\mu)
\rightarrow 0$
with each $P_i(\mu)$ projective.
In case $\#\mu = 0$, $K(\mu)$ is already projective,
so we can simply take $P_0(\mu) = K(\mu)$ and $P_i(\mu) = 0$ for $i > 0$.
Now suppose $\#\mu > 0$.
Let $g = f_\mu$ and apply Procedure~\ref{alg3} to construct
$h = f_\nu$ and operators $X_a, Y_a \in \{E_a,F_a\}_{a \in \Z}$.
Since Procedure~\ref{alg3} reduces $\mu$ to a typical weight in
finitely many steps, we may assume inductively that we have already
constructed an exact sequence
\begin{equation}\label{indseq}
P_d(\nu)\longrightarrow\dots\longrightarrow P_0(\nu)\longrightarrow K(\nu)
\longrightarrow 0.
\end{equation}
Now we consider two cases.
Suppose first that $\#\nu = \#\mu$.
Then $X_a K(\nu) \simeq K(\mu)$, 
so applying $X_a$ to (\ref{indseq}) gives us
the desired sequence
with $P_i(\mu) = X_a P_i(\nu)$.
In the second case, $\#\nu = \#\mu-1$, and
$X_a K(\nu)$ has a two step filtration 
with $K(\mu)$ at the top
and $K(\tilde X_a (\nu))$ at the bottom.
Applying $X_a$ to (\ref{indseq}) gives us an exact sequence
$X_a P_d(\nu)\rightarrow\dots\rightarrow X_a P_0(\nu)
\rightarrow X_a K(\nu)\rightarrow 0.$
By induction on $d$, we may assume in addition that we have already 
constructed an exact sequence
$P_{d-1}(\tilde X_a (\nu))\rightarrow\dots\rightarrow
P_0(\tilde X_a (\nu))\rightarrow K(\tilde X_a (\nu))
\rightarrow 0.$
Applying the comparison theorem \cite[2.2.6]{Wei} to the embedding
$i:K(\tilde X_a (\nu)) \hookrightarrow X_a K(\nu)$, we get vertical maps
making the diagram commute:
$$
\begin{CD}
\cdots @>>> P_1(\tilde X_a (\nu)) @>>>  P_0(\tilde X_a (\nu)) @>>>
K(\tilde X_a (\nu)) @>>>0\\
&&@VVV@VVV@VViV\\
\cdots @>>>X_a P_1(\nu) @>>>  X_a P_0(\nu) @>>>
X_a K(\nu) @>>>0
\end{CD}
$$
The total complex of this double complex is exact by the acyclic
assembly lemma \cite[2.7.3]{Wei}. 
Factoring out $K(\tilde X_a(\mu))$ yields the required
exact sequence
$$
\dots\longrightarrow
X_a P_1(\nu) \oplus  P_0(\tilde X_a(\nu))
\longrightarrow  X_a P_0(\nu)\longrightarrow
K(\mu) 
\longrightarrow 0
$$
This time, $P_i(\mu) = X_a P_i(\nu) \oplus P_{i-1}(\tilde X_a (\nu))$.

Replacing $d$ by $(d+1)$, the same procedure
constructs an exact sequence
$P_{d+1}(\mu)\rightarrow P_d(\mu) \rightarrow\dots\rightarrow
P_0(\mu)\rightarrow K(\mu)\rightarrow 0,$
where we can always ensure that the first $d$ terms are {\em the same} as
the ones constructed before.
Now letting $d \rightarrow \infty$ we get a projective resolution of $K(\mu)$.
We note moreover by the construction that
whenever $P(\lambda)$ is a summand of $P_i(\mu)$ for some
$i \geq 0$, i.e. $\hom_{\mathcal F_{m|n}}(P_i(\mu), L(\la)) \neq 0$,
we must have that $\mu \preceq \la$.

Finally let $p_{\mu,\la}(q) =
\sum_{i \geq 0} \dim 
\hom_{\mathcal F_{m|n}}(P_i(\mu),L(\la)) q^i$.
To complete the proof, we need to show that $p_{\mu,\la}(q) = 
l_{\mu,\la}(-q^{-1})$
for each $\mu,\la \in X^+(m|n)$. 
For this, we show that 
the polynomials $p_{\mu,\la}(q)$ satisfy the same
relations as the polynomials $l_{\mu,\la}(-q^{-1})$ in Lemma~\ref{secondalg}.
Once this is established, the algorithm explained at the end of
\ref{ta} to compute $l_{\mu,\la}(-q^{-1})$ also computes $p_{\mu,\la}(q)$,
hence $p_{\mu,\la}(q) = l_{\mu,\la}(-q^{-1})$.
So take $\la,\mu \in X^+(m|n)$ with $\#\mu > 0$, where
we may assume that $\wt(\la) = \wt(\mu)$ since otherwise
$p_{\mu,\la}(q) =0$ and the conclusion holds trivially.
Apply Procedure~\ref{alg3} to $g = f_\mu$ to get $h = f_\nu$ and
operators $X_a,Y_a$,
and consider the two cases $\#\nu = \#\mu$ or $\#\nu = \#\mu-1$, just like 
above.
Let us just explain the argument in the second case, 
the first case being easier.
Since $\wt(\la) = \wt(\mu)$ and $\mu$ is not at the end of
an $a$-string of length 2
in the dual crystal graph, 
Theorem~\ref{kuj} shows that $Y_a L(\la)$ equals $L(\tilde Y_a^* (\la))$,
interpreted as $0$ if $\tilde Y_a^* (\la) = \varnothing$.
So we get that
\begin{multline*}
\sum_{i \geq 0} \dim \hom_{\mathcal F_{m|n}}(X_a P_i(\nu),L(\la)) q^i = 
\sum_{i \geq 0} \dim \hom_{\mathcal F_{m|n}}(P_i(\nu),Y_a L(\la)) q^i\\
=\sum_{i \geq 0} \dim \hom_{\mathcal F_{m|n}}(P_i(\nu),L(\tilde Y_a^*(\la))) 
q^i = p_{\nu,\tilde Y_a^* (\la)}(q),
\end{multline*}
interpreted as $0$ in case $\tilde Y_a^* (\la) = \varnothing$.
We noted above that
$P_i(\mu) = X_a P_i(\nu) \oplus P_{i-1} (\tilde X_a (\nu))$,
hence we get that
$p_{\mu,\la}(q) =
p_{\nu,\tilde Y_a^*(\la)}(q)  + q p_{\tilde X_a (\nu),\la}(q),$
which is what we wanted in this case, cf. Lemma~\ref{secondalg}.
\end{proof}

Now choose a length function
on $\Z^{m|n}_+$ as explained in
\ref{lf}, and lift it to $X^+(m|n)$ by setting
$\ell(\la) := \ell(f_\la)$.

\begin{Theorem}\label{parcond}
For $\mu,\la \in X^+(m|n)$,
the superspace
$\ext^i_{\mathcal F_{m|n}}(K(\mu),L(\la))$ is concentrated in degree $\bar 0$,
and
$$
\sum_{i \geq 0} \dim \ext^i_{\mathcal F_{m|n}}(K(\mu),L(\la)) q^i
= l_{\mu,\la}(-q^{-1}).
$$
Hence, 
$\ext^{\bullet}_{{\scriptstyle \mathcal F_{m|n}}}(K(\mu), L(\la)) \neq 0$ if and only
if $\mu \preceq \la$, in which case
\begin{itemize}
\item[(i)]
$\ext^i_{\mathcal F_{m|n}}(K(\mu), L(\la)) \neq 0
\Rightarrow 
i \leq \ell(\la) - \ell(\mu),
i \equiv \ell(\la)-\ell(\mu) \pmod{2}$;
\item[(ii)] $\ext^{\ell(\la)-\ell(\mu)}_{\mathcal F_{m|n}}(K(\mu), L(\la))$
is exactly one dimensional;
\item[(iii)] $\ext^1_{\mathcal F_{m|n}}(K(\mu),L(\la))$
is at most one dimensional.
\end{itemize}
\end{Theorem}

\begin{proof}
Apply the functor $\hom_{\mathcal F_{m|n}}(?, L(\la))$
to the projective resolution constructed in Lemma~\ref{res} and use
Corollaries~\ref{parcor}, \ref{clearer}(ii) and \ref{indegzero}.
\end{proof}

By \cite[Theorem 7.6]{Zou},
Theorem~\ref{parcond} shows that the polynomials
$l_{\mu,\lambda}(-q^{-1})$ defined here agree with 
the Kazhdan-Lusztig polynomials $K_{\lambda,\mu}(q)$
defined by Serganova in \cite{Serg0,Serg}.
(It also proves \cite[Conjecture 4.4]{JZ}, 
and answers a question raised at the end of \cite{Zou}.)
Thus we have a cohomological interpretation of the polynomial
 $l_{\mu,\la}(-q^{-1})$, analogous to
Vogan's interpretation \cite[Conjecture 3.4]{V1}
of Kazhdan-Lusztig polynomials
in category $\mathcal O$ for a semisimple Lie algebra.
The even-odd vanishing
established in Theorem~\ref{parcond}(i) is especially important:
in the language of \cite{CPSkl}, it shows that the highest weight category
$\mathcal F_{m|n}^{\0}$ has a ``Kazhdan-Lusztig theory''.
Applying \cite[Corollary 3.9]{CPSinf} 
(and Corollary~\ref{indegzero} again), we obtain:

\begin{Corollary}\label{extLL}
For $\mu,\la \in X^+(m|n)$,
the superspace
$\ext^i_{\mathcal F_{m|n}}(L(\mu),L(\la))$ is concentrated in degree $\bar 0$,
and
$$
\sum_{i \geq 0} \dim \ext^i_{\mathcal F_{m|n}}(L(\mu),L(\la)) q^i
= \sum_{\nu \in X^+(m|n)} l_{\nu,\mu}(-q^{-1}) l_{\nu,\la}(-q^{-1}).
$$
In particular,
$\ext^i_{\mathcal F_{m|n}}(L(\mu), L(\la)) \neq 0
\Rightarrow i \equiv \ell(\la)-\ell(\mu) \pmod{2}$.
\end{Corollary}

\begin{Example}\rm
Take $\mu =\la = 0$ and let $r = \min(m,n)$.
In this case, the polynomials
$l_{\nu,\mu}(-q^{-1})$
are computed explicitly in Example~\ref{goodie}.
Combining this with Corollary~\ref{extLL}, one deduces
that $\dim \ext^{2i}_{\mathcal F_{m|n}}(\C,\C)$
equals the number of partitions of $i$ with at most $r$ non-zero parts.
Hence:
\begin{equation}
\sum_{i \geq 0} \dim \ext^i_{\mathcal F_{m|n}}(\C,\C) q^i
= 
\frac{1}{(1-q^2) (1-q^4) \dots (1-q^{2r})}.
\end{equation}
\end{Example}

\end{document}